%% file: main.tex
\documentclass[11pt]{article}

\usepackage[a4paper,margin=1in]{geometry}

\usepackage[T1]{fontenc}
\usepackage[utf8]{inputenc}
\usepackage[english]{babel}
\usepackage{bm}

\usepackage{amsmath}
\usepackage{amsthm}
\usepackage{amssymb} 
\usepackage{mathtools}
\usepackage{calc}
\usepackage[short, nocomma]{optidef}
\usepackage{accents}
\usepackage{scalerel}

\usepackage{booktabs}
\usepackage{color}
\usepackage{xcolor}
\usepackage{afterpage}
\usepackage{caption}
\usepackage[labelfont=sf]{subcaption}

\usepackage{eurosym}
\usepackage{siunitx}
\usepackage{algorithm}
\usepackage{algpseudocode}

\usepackage{tikz}
\usepackage{pgfplots}
\pgfplotsset{compat=1.18}
\usetikzlibrary{svg.path, fit, patterns, shapes.geometric, positioning, calc}

\usepackage{layouts}
\usepackage{hyperref}

\providecommand{\mathdefault}[1]{#1}

\newcommand{\ubar}[1]{\underaccent{\bar}{#1}}

\makeatletter
\newcommand{\pushright}[1]{\ifmeasuring@#1\else\omit\hfill$\displaystyle{#1}$\fi\ignorespaces}
\newcommand{\pushleft}[1]{\ifmeasuring@#1\else\omit$\displaystyle{#1}$\hfill\fi\ignorespaces}
\makeatother

\theoremstyle{plain}
\newtheorem{theorem}{Theorem}

\newtheorem{proposition}{Proposition}

\theoremstyle{definition}

\theoremstyle{remark}

\makeatletter
\renewcommand{\proof}[1]{\par\medskip\noindent{\itshape #1}\enspace\ignorespaces}
\renewcommand{\endproof}{\par\medskip}
\makeatother
\newcommand{\Halmos}{\hfill\ensuremath{\square}}

\DeclareMathOperator*{\argmax}{arg\,max}

\newcommand{\up}{\rule{0pt}{2.6ex}}
\newcommand{\down}{\rule[-1.2ex]{0pt}{0pt}}

\makeatletter
\newcommand{\informs@note}[1]{%
  \def\informs@tmp{#1}%
  \ifx\informs@tmp\@empty\else
    \par\smallskip{\footnotesize\noindent\textit{Note.}\enspace #1}%
  \fi}
\newcommand{\FIGURE}[3]{\centering\caption{#2}\par\smallskip#1\par\informs@note{#3}}
\newcommand{\TABLE}[3]{\centering\caption{#1}\par\smallskip#2\par\informs@note{#3}}
\makeatother

\input{common-resources.tex}

\usepackage{natbib}
\bibpunct[, ]{(}{)}{;}{a}{,}{,}%
\newcommand{\spdual}[1]{\;\mathmakebox[\widthof{$[\,\DualOverCLbSym,\DualOverCUbSym\,]$}][l]{#1}}
\newcommand{\doi}[1]{doi: \href{https://doi.org/#1}{\nolinkurl{#1}}}

\addto\extrasenglish{

}

\hypersetup{
  colorlinks = true,
  linkcolor  = black,
  citecolor  = black,
  urlcolor   = blue,
}

\begin{document}

\title{\bfseries Fast Subproblem Solving in Benders Decomposition for
Affine Potential-Based Flow Problems with Topology Switching and Scenarios}

\author{%
  Oliver Gaul\thanks{Corresponding author.}\\
  {\small Chair of Operations Research}\\
  {\small RWTH Aachen University}\\
  {\footnotesize\texttt{gaul@or.rwth-aachen.de}}\\
  {\footnotesize ORCID: \href{https://orcid.org/0000-0002-2131-1911}{0000-0002-2131-1911}}
  \and
  Tim Donkiewicz\\
  {\small Chair of Operations Research}\\
  {\small RWTH Aachen University}\\
  {\footnotesize\texttt{donkiewicz@or.rwth-aachen.de}}\\
  {\footnotesize ORCID: \href{https://orcid.org/0000-0002-5721-3563}{0000-0002-5721-3563}}
  \and
  Marco Lübbecke\\
  {\small Chair of Operations Research}\\
  {\small RWTH Aachen University}\\
  {\footnotesize\texttt{luebbecke@or.rwth-aachen.de}}\\
  {\footnotesize ORCID: \href{https://orcid.org/0000-0002-2635-0522}{0000-0002-2635-0522}}
}
\date{\today}

\maketitle

\begin{abstract}
\noindent
Potential-based flow problems arise for a range of practical problems, like the operation of lossless power networks.
The distribution of flows through the network is governed by the potential equation, which relates arc flows and node potentials following arc potential functions.
We consider a general variant with affine potential functions, topology switching, and scenarios.
It can be solved using a Benders decomposition algorithm, with each scenario corresponding to a subproblem.
Subproblem solving time can then significantly contribute to overall running time.
We present a fast algorithm to obtain Benders optimality and feasibility cuts from fixed master solutions.
Problems with affine potential functions are reduced to equivalent problems where every potential function is the identity function.
Optimal primal and dual solutions as well as unbounded dual rays are determined by systems of linear equations (SLE), resulting in a constant number of SLEs solved per subproblem.
SLEs are solved using LU-factorizations, with the similarity between primal and dual, scenarios, and master solutions exploited by reusing factorizations via low-rank updates to significantly reduce theoretical complexity.
We demonstrate the method on security-constrained optimal transmission switching instances, where subproblem solve time reductions of 97.98 percent compared to a standard LP solver are observed.

\end{abstract}

\section{Introduction}
\input{chapters/0-introduction.tex}

\section{Related Work}
\input{chapters/1-related-work.tex}

\section{Potential-Based Flow Problems}\label{sec:pot_flows}
\input{chapters/2-potential-flows.tex}

\section{A Benders Decomposition Approach}\label{sec:benders}
\input{chapters/3-benders.tex}

\section{Solving Subproblems as Systems of Linear Equations}\label{sec:single_sp}
\input{chapters/4-sle.tex}

\section{Evaluation}\label{sec:evaluation}
\input{chapters/5-evaluation.tex}

\section{Conclusion}
\input{chapters/6-conclusion.tex}

\section*{Acknowledgments}
The authors want to thank Alexander Helber for pointing out the uniqueness
property for the primal subproblems.
This research was supported by the German Federal Ministry for Economic
Affairs and Climate Action (BMWK) under grant 03EI4056D within the InnOpTEM
project (September 2022--August 2025).

\bibliographystyle{plainnat}
\bibliography{literature}

\end{document}

%% file: common-resources.tex
\newcommand{\subto}{\operatorname{s.t.}}
\newcommand{\N}{\ensuremath{\mathbb{N}}}

\newcommand{\R}{\ensuremath{\mathbb{R}}}
\newcommand{\Z}{\ensuremath{\mathbb{Z}}}

\newcommand{\DemandSym}{b}                   
\newcommand{\C}{\ensuremath{C}}              
\newcommand{\DemandParSym}[0]{\ensuremath{b}}  
\newcommand{\PotSym}[0]{\ensuremath{\pi}} 
\newcommand{\CondSym}[0]{\mu}               
\newcommand{\FlowSym}[0]{\ensuremath{x}}     
\newcommand{\FlowParSym}[0]{\ensuremath{x}}  
\newcommand{\FlowOverSym}[0]{\ensuremath{o}}
\newcommand{\Arcs}{A}                        
\newcommand{\LArcs}{\Arcs_L}                 
\newcommand{\SArcs}{\Arcs_S}                 
\newcommand{\SArcsOn}{\SArcs^1}                 
\newcommand{\SArcsOff}{\SArcs^0}                 
\newcommand{\subSym}{W}                      

\newcommand{\TAL}{Z}
\newcommand{\DemandCostsSym}{F}

\newcommand{\fixedSym}[1]{\hat{{#1}}}

\newcommand{\etaFixed}{\fixedSym{\eta}}

\newcommand{\aC}{\ensuremath{\mathcal{C}}} 

\newcommand{\Demand}[1]{\DemandSym_{#1}}
\newcommand{\DemandCosts}[1]{\DemandCostsSym_{#1}}

\newcommand{\DemandMinSym}[0]{\ensuremath{\ubar{\DemandParSym}}}
\newcommand{\DemandMin}[1]{\ensuremath{\DemandMinSym_{{#1}}}}
\newcommand{\DemandMaxSym}[0]{\ensuremath{\bar{\DemandParSym}}}
\newcommand{\DemandMax}[1]{\ensuremath{\DemandMaxSym_{{#1}}}}

\newcommand{\DemandFixedSym}[0]{\ensuremath{\fixedSym{\DemandSym}}}
\newcommand{\DemandFixed}[1]{\ensuremath{\DemandFixedSym_{{#1}}}}

\newcommand{\PotFSym}[0]{\psi}

\newcommand{\PotMinSym}[0]{\ensuremath{\ubar{\PotSym}}}
\newcommand{\PotMin}[1]{\ensuremath{\PotMinSym_{#1}}}
\newcommand{\PotMaxSym}[0]{\ensuremath{\bar{\PotSym}}}
\newcommand{\PotMax}[1]{\ensuremath{\PotMaxSym_{#1}}}
\newcommand{\Pot}[1]{\ensuremath{\PotSym_{#1}}}
\newcommand{\PotCSym}[0]{\ensuremath{\PotSym^{\C}}}
\newcommand{\PotC}[1]{\ensuremath{\PotCSym_{#1}}}
\newcommand{\PotCOptSym}[0]{\ensuremath{\PotSym^{\C*}}}
\newcommand{\PotCOpt}[1]{\ensuremath{\PotCOptSym_{#1}}}

\newcommand{\PotDifMin}[2]{\ensuremath{\ubar{\PotSym}_{({#1},{#2})}}}
\newcommand{\PotDifMax}[2]{\ensuremath{\bar{\PotSym}_{({#1},{#2})}}}

\newcommand{\Cond}[2]{\ensuremath{\CondSym_{({#1},{#2})}}}

\newcommand{\PotF}[2]{\PotFSym_{({#1},{#2})}}
\newcommand{\PotFSlope}[2]{r_{({#1},{#2})}}
\newcommand{\PotFConst}[2]{s_{({#1},{#2})}}

\newcommand{\FlowMinSym}[0]{\ensuremath{\ubar{\FlowParSym}}}
\newcommand{\FlowMin}[2]{\ensuremath{\FlowMinSym_{({#1},{#2})}}}
\newcommand{\FlowMaxSym}[0]{\ensuremath{\bar{\FlowParSym}}}
\newcommand{\FlowMax}[2]{\ensuremath{\FlowMaxSym_{({#1},{#2})}}}
\newcommand{\FlowMinOverSym}[0]{\ensuremath{\underaccent{\tilde}{\FlowParSym}}}
\newcommand{\FlowMinOver}[2]{\ensuremath{\FlowMinOverSym_{({#1},{#2})}}}
\newcommand{\FlowMaxOverSym}[0]{\ensuremath{\tilde{\FlowParSym}}}
\newcommand{\FlowMaxOver}[2]{\ensuremath{\FlowMaxOverSym_{({#1},{#2})}}}

\newcommand{\FlowOverC}[3]{\ensuremath{\FlowOverSym_{({#1},{#2})}^{{#3}}}}
\newcommand{\FlowOverCSym}[0]{\ensuremath{\FlowOverSym^{\C}}}
\newcommand{\FlowOverCOptSym}[0]{\ensuremath{\FlowOverSym^{\C*}}}
\newcommand{\FlowOverCOpt}[2]{\ensuremath{\FlowOverCOptSym_{({#1},{#2})}}}
\newcommand{\Flow}[2]{\ensuremath{\FlowSym_{({#1},{#2})}}}
\newcommand{\FlowOptSym}{\ensuremath{\FlowSym^{*}}}
\newcommand{\FlowOpt}[2]{\ensuremath{\FlowOptSym_{({#1},{#2})}}}
\newcommand{\FlowCSym}[0]{\ensuremath{\FlowSym^{\C}}}
\newcommand{\FlowC}[2]{\ensuremath{\FlowCSym_{({#1},{#2})}}}
\newcommand{\FlowCOptSym}[0]{\ensuremath{\FlowSym^{\C*}}}
\newcommand{\FlowCOpt}[2]{\ensuremath{\FlowCOptSym_{({#1},{#2})}}}

\newcommand{\OverCost}[0]{\ensuremath{\DemandCostsSym^{\FlowOverSym}}}

\newcommand{\zSym}[0]{\ensuremath{z}}
\newcommand{\z}[2]{\ensuremath{\zSym_{({#1},{#2})}}}

\newcommand{\zFixedSym}[0]{\ensuremath{\fixedSym{\zSym}}}
\newcommand{\zFixed}[2]{\ensuremath{\zFixedSym_{({#1},{#2})}}}

\newcommand{\ArcsOut}{\Arcs^\C}

\newcommand{\LArcsOutSym}[1]{\LArcs^{#1}}
\newcommand{\LArcsOut}{\LArcsOutSym{\C}}
\newcommand{\LArcsCSym}[2]{\LArcs^{{#1},{#2}}}
\newcommand{\LArcsC}[1]{\LArcsCSym{\zFixedSym}{#1}}
\newcommand{\LArcsCC}{\LArcsC{\C}}

\newcommand{\SArcsOutSym}[1]{\SArcs^{#1}}
\newcommand{\SArcsOut}{\SArcsOutSym{\C}}
\newcommand{\SArcsCSym}[1]{\SArcs^{{#1,\C}}}
\newcommand{\SArcsz}{\SArcsCSym{\zFixedSym}}
\newcommand{\SArcszNC}{\SArcsOutSym{\zFixedSym}}

\newcommand{\InOutArcs}[1]{\delta({#1}) \setminus \C}
\newcommand{\InArcs}[1]{\delta^{-}({#1}) \setminus \C}
\newcommand{\OutArcs}[1]{\delta^{+}({#1}) \setminus \C}
\newcommand{\InArcsNC}[1]{\delta^{-}({#1})}
\newcommand{\OutArcsNC}[1]{\delta^{+}({#1})}

\newcommand{\Sub}{\subSym}
\newcommand{\OtherSub}{U}
\newcommand{\Subs}{\mathcal{\subSym}^{\zFixedSym,\C}}
\newcommand{\SubsNC}{\mathcal{\subSym}^{\zFixedSym}}
\newcommand{\SubNode}[1]{u({#1})}

\newcommand{\lbsym}[1]{\ubar{{#1}}}
\newcommand{\ubsym}[1]{\bar{{#1}}}
\newcommand{\DualFCSym}[0]{\ensuremath{\varphi}}
\newcommand{\DualFCCSym}[0]{\ensuremath{\DualFCSym^\C}}
\newcommand{\DualFCCOptSym}[0]{\ensuremath{\DualFCSym}^{\C*}}
\newcommand{\DualFCCOpt}[1]{\ensuremath{\DualFCCOptSym_{#1}}}
\newcommand{\DualFCC}[1]{\ensuremath{\DualFCCSym_{#1}}}

\newcommand{\DualPESym}[0]{\ensuremath{\omega}}
\newcommand{\DualPECSym}[0]{\ensuremath{\DualPESym^\C}}
\newcommand{\DualPECOptSym}[0]{\ensuremath{\DualPESym^{\C*}}}
\newcommand{\DualPEC}[2]{\ensuremath{\DualPECSym_{({#1},{#2})}}}
\newcommand{\DualPECLb}[2]{\ensuremath{\lbsym{\DualPESym}^{\C}_{({#1},{#2})}}}
\newcommand{\DualPECUbSym}[0]{\ensuremath{\ubsym{\DualPESym}^{\C}}}
\newcommand{\DualPECUb}[2]{\ensuremath{\DualPECUbSym_{({#1},{#2})}}}

\newcommand{\DualOverSym}[0]{\ensuremath{\sigma}}

\newcommand{\DualOverCOptSym}[0]{\ensuremath{\DualOverSym^{\C*}}}

\newcommand{\DualOverCOpt}[2]{\ensuremath{\DualOverCOptSym_{({#1},{#2})}}}
\newcommand{\DualOverCLbSym}[0]{\ensuremath{\lbsym{\DualOverSym}^{\C}}}
\newcommand{\DualOverCLb}[2]{\ensuremath{\DualOverCLbSym_{({#1},{#2})}}}
\newcommand{\DualOverCLbOptSym}[0]{\ensuremath{\lbsym{\DualOverSym}^{\C*}}}
\newcommand{\DualOverCLbOpt}[2]{\ensuremath{\DualOverCLbOptSym_{({#1},{#2})}}}
\newcommand{\DualOverCUbSym}[0]{\ensuremath{\ubsym{\DualOverSym}^{\C}}}
\newcommand{\DualOverCUb}[2]{\ensuremath{\DualOverCUbSym_{({#1},{#2})}}}
\newcommand{\DualOverCUbOptSym}[0]{\ensuremath{\ubsym{\DualOverSym}^{\C*}}}
\newcommand{\DualOverCUbOpt}[2]{\ensuremath{\DualOverCUbOptSym_{({#1},{#2})}}}

\newcommand{\DualSOnSym}[0]{\ensuremath{\alpha}}
\newcommand{\DualSOnCSym}[0]{\ensuremath{\DualSOnSym^\C}}
\newcommand{\DualSOnC}[2]{\ensuremath{\DualSOnCSym_{({#1},{#2})}}}
\newcommand{\DualSOnCLbSym}[0]{\ensuremath{\lbsym{\DualSOnSym}^{\C}}}
\newcommand{\DualSOnCUbSym}[0]{\ensuremath{\ubsym{\DualSOnSym}^{\C}}}

\newcommand{\DualSOnCUb}[2]{\ensuremath{\DualSOnCUbSym_{({#1},{#2})}}}
\newcommand{\DualSOnCLbOptSym}[0]{\ensuremath{\lbsym{\DualSOnSym}^{\C*}}}
\newcommand{\DualSOnCUbOptSym}[0]{\ensuremath{\ubsym{\DualSOnSym}^{\C*}}}
\newcommand{\DualSOnCLbOpt}[2]{\ensuremath{\DualSOnCLbOptSym_{({#1},{#2})}}}
\newcommand{\DualSOnCUbOpt}[2]{\ensuremath{\DualSOnCUbOptSym_{({#1},{#2})}}}

\newcommand{\DualSOffSym}[0]{\ensuremath{\tau}}

\newcommand{\DualSOffCLbOptSym}[0]{\ensuremath{\lbsym{\DualSOffSym}^{\C*}}}
\newcommand{\DualSOffCLbOpt}[2]{\ensuremath{\DualSOffCLbOptSym_{({#1},{#2})}}}
\newcommand{\DualSOffCUbOptSym}[0]{\ensuremath{\ubsym{\DualSOffSym}^{\C*}}}
\newcommand{\DualSOffCUbOpt}[2]{\ensuremath{\DualSOffCUbOptSym_{({#1},{#2})}}}
\newcommand{\DualSOffCLbSym}{\ensuremath{\lbsym{\DualSOffSym}^{\C}}}
\newcommand{\DualSOffCUbSym}{\ensuremath{\ubsym{\DualSOffSym}^{\C}}}

\newcommand{\DualFlowBSym}[0]{\ensuremath{\gamma}}

\newcommand{\DualFlowBCLbSym}[0]{\ensuremath{\lbsym{\DualFlowBSym}^{\C}}}
\newcommand{\DualFlowBCUbSym}[0]{\ensuremath{\ubsym{\DualFlowBSym}^{\C}}}

\newcommand{\DualFlowBCUb}[2]{\ensuremath{\DualFlowBCUbSym_{({#1},{#2})}}}
\newcommand{\DualFlowBCLbOpt}[2]{\ensuremath{\lbsym{\DualFlowBSym}^{\C*}_{({#1},{#2})}}}
\newcommand{\DualFlowBCUbOpt}[2]{\ensuremath{\ubsym{\DualFlowBSym}^{\C*}_{({#1},{#2})}}}

\newcommand{\DualPotBSym}[0]{\ensuremath{\xi}}

\newcommand{\DualPotBCLbSym}[0]{\ensuremath{\lbsym{\DualPotBSym}^{\C}}}
\newcommand{\DualPotBCUbSym}[0]{\ensuremath{\ubsym{\DualPotBSym}^{\C}}}

\newcommand{\DualPotBCLbOpt}[1]{\ensuremath{\lbsym{\DualPotBSym}^{\C*}_{{#1}}}}
\newcommand{\DualPotBCUb}[1]{\ensuremath{\DualPotBCUbSym_{{#1}}}}
\newcommand{\DualPotBCUbOpt}[1]{\ensuremath{\ubsym{\DualPotBSym}^{\C*}_{{#1}}}}

\newcommand{\subP}{\ensuremath{S_{\C}(\DemandFixedSym, \zFixedSym, \etaFixed_\C)}}

\newcommand{\lesP}{\ensuremath{L^{\zFixedSym, \C, p}}}
\newcommand{\lesD}{\ensuremath{L^{\zFixedSym, \C, d}}}

\newcommand{\ie}{i.e.}
\newcommand{\eg}{e.g.}

\definecolor{rwth-blue}{RGB}{0,84,159}
\definecolor{rwth-lblue}{RGB}{142,186,229}
\definecolor{rwth-purple}{RGB}{122,111,172}
\definecolor{rwth-violet}{RGB}{97,33,88}
\definecolor{rwth-carmine}{RGB}{161,16,53}
\definecolor{rwth-red}{RGB}{204,7,30}
\definecolor{rwth-magenta}{RGB}{227,0,102}
\definecolor{rwth-orange}{RGB}{246,168,0}
\definecolor{rwth-yellow}{RGB}{255,237,0}
\definecolor{rwth-grass}{RGB}{189,171,39}
\definecolor{rwth-green}{RGB}{87,171,39}
\definecolor{rwth-cyan}{RGB}{0,152,161}
\definecolor{rwth-teal}{RGB}{0,97,101}
\definecolor{rwth150red}{HTML}{a11035}

\DeclareDocumentCommand\colorsOn{}{
  \definecolor{visual-red}{RGB}{255,0,0}
  \definecolor{visual-blue}{RGB}{0,0,255}
  \definecolor{visual-green}{RGB}{0,128,0}
  \definecolor{visual-purple}{RGB}{220,0,220}
  \definecolor{visual-orange}{RGB}{255,128,0}
  \definecolor{visual-cyan}{RGB}{0,200,200}
}

\colorsOn

\newcommand{\autocref}[1]{Constraint \eqref{#1}}
\newcommand{\autocsref}[1]{Constraints \eqref{#1}}
\newcommand{\autoeqref}[1]{Equation \eqref{#1}}

\definecolor{ibm-blue}{RGB}{100,143,255}
\definecolor{ibm-purple}{RGB}{120,94,240}
\definecolor{ibm-red}{RGB}{220,38,127}
\definecolor{ibm-orange}{RGB}{254,97,0}
\definecolor{ibm-yellow}{RGB}{255,176,0}

\colorlet{color1}{ibm-blue}
\colorlet{color2}{ibm-orange}
\colorlet{color3}{ibm-red}
\colorlet{color4}{ibm-purple}
\colorlet{color5}{ibm-yellow}

%% file: chapters/0-introduction.tex
Infrastructure networks are an important part of modern supply systems.
They allow for the transportation of critical goods, such as water, electricity or gas from producers to household or industrial consumers.
Due to the scale of operations, efficient planning can result in significant cost savings.
This has become more pronounced in recent years, as the ongoing green energy transition and the resulting transitory systems prove particularly challenging in certain areas.
For example, the German transmission system is experiencing a shift from relatively few, centralized power plants to many decentralized renewable energy producers \citep{redispatchGermany}, as well as an imbalanced distribution of large-scale consumers and producers \citep{redispatchGermany,redispatch}.
These challenges necessitate the expansion and adaptation of existing grid infrastructure, with the current expansion efforts lagging behind the necessary capacities \citep{redispatchGermany,redispatch}.
This has caused an increase in costly operational measures such as the so-called re-dispatch, where already planned generation is changed in relatively short time periods, thus incurring higher costs, to redirect flows within the network \citep{redispatch}.
In Germany alone, re-dispatch cost has increased from \euro113.3 million in 2013 to \euro2,689.2 million in 2022 \citep{redispatchGermany}.

As such, there is significant potential for the application of mathematical optimization in this domain.
Classical flow problem formulations, where each arc of a network carries a flow to satisfy node demands, assume that the flow can be directly controlled, for example via trucks.
However, oftentimes flows can only passively be controlled by adjusting network and demand parameters.
Distribution through the network then depends on arc parameters and the relation between the ``pressure'' of adjacent nodes.
The behavior of these flows can be modeled through the \emph{potential-based flow} framework \citep[\eg,][]{gross2019algorithmic}.
There, a flow problem is extended by variable node \emph{potentials}, measuring the ``pressure'' at each node of the network, as well as constant arc \emph{conductance} and \emph{potential function} to capture the physical properties of the flow.
Flows and potentials then need to satisfy the \emph{potential equation}, which links the incident potentials with conductance and flow following the potential function.
Next to these so-called \emph{line arcs}, \emph{switch arcs} freely carry flow like arcs in a classical flow problem, but can be switched on or off to disrupt the synchronization of potentials of incident nodes.
Based on the choice of potential function, this framework unifies different applications including lossless electrical power systems, and gas or water networks \citep{hendrickson1984common,pfetsch_potential-based_2026,rockafellar1998network}.

An important additional consideration is the robustness of the network against certain contingency scenarios, like the possible outage of a set of lines.
An explicit analysis might even be legally required \citep{n1Germany}.
This requires the simultaneous analysis of the effects of shared decisions on the individual flows and potentials of all scenarios.
Structurally, this naturally lends itself to Benders decomposition \citep{benders1962partitioning}: the master problem determines the shared switching and demand variables, while subproblems verify flow feasibility and solution quality for each scenario.
The application of Benders decomposition to potential-based flow problems has been thoroughly studied in the literature \citep{binato2001new,humpola2017gas,khodaei2010transmission,ruiz2015robust,wu2010accelerating,yang2017modeling}.
However, these works usually focus on improvements to the master problem, for example via cuts, to reduce master problem linear relaxation solving time and branch-and-bound tree size.
The efficiency of the subproblem solution algorithm remains a relatively underexplored subject in the Benders context, even though it can contribute significantly to overall running time.

%% file: chapters/1-related-work.tex
Early works explicitly discussing potential-based flow problems include those by \cite{birkhoff1956potential}, who already show the uniqueness of flows, and \cite{rockafellar1998network}.
Recently, the problems have found increased attention. For example, \cite{borner2025valid} discuss valid inequalities based on potential bounds, \cite{gross2019algorithmic} study the hardness of certain related algorithmic problems, and
\cite{pfetsch_potential-based_2026} provide an overview over methods and results used in the context of mathematical optimization, including also application-specific work.
These applications include lossless DC power flow \citep{fisher2008optimal,hedman_optimal_2009}, which can be modeled with affine potential functions, as well as gas \citep{correa2014gas,koch2015evaluating,schmidt2022nonlinear} and water networks \citep{burgschweiger2009optimization, mozafari2020robust}.
As mentioned by \cite{birkhoff1956potential}, related problems include, for example, traffic assignment \citep{dafermos1969traffic,leblanc1975efficient}, where arcs are less attractive the more they are used.

Switch arcs correspond to different mechanisms, like valve switching and compressor activation control network configuration in gas networks \citep{correa2014gas}, the opening and closing of switches in electrical systems \citep{fisher2008optimal,heidarifar2015network}, or configuring valves to en- or disable and control the direction of flow in water networks \citep{burgschweiger2009optimization}.

Scenarios occur in the context of robust optimization and uncertainty: a set of decisions have to be taken ``here-and-now'' before one of possibly infinitely many scenarios realizes, after which the ``wait-and-see'' decisions can be adjusted.
They can, for example, cover changes in loads \citep{mozafari2020robust,thurauf2025adjustable} or topology \citep{hedman_optimal_2009,li2018minimax}.
Both, an aggregate view over the scenarios \citep{hedman_optimal_2009} and a worst-case view \citep{li2018minimax,mozafari2020robust,thurauf2025adjustable} can be considered.

The Benders decomposition algorithm \citep{benders1962partitioning} is a widely used exact algorithm where complicating variables are temporarily fixed to reduce overall computational burden.
\cite{rahmaniani2017benders} provide a survey of modern approaches.
For affine (linear) potential-based flow problems, this is a commonly used method, with research focusing on formulation or cut design.
A common application are linear power flow problems.
For example,
\cite{binato2001new} solve a power transmission network design problem with Benders decomposition, using big-$M$ tightening and Gomory cuts.
\cite{khodaei2010transmission} determine generation schedules over multiple security scenarios for a security-constrained unit commitment problem.
\cite{wu2010accelerating} propose strengthened Benders cuts for a network-constrained unit commitment problem.
\cite{ruiz2015robust} solve a three-level robust transmission expansion planning problem, using a column-and-constraint generation method with a specialized cut computation procedure.

For nonlinear problems, Benders decomposition is mostly applied to gas networks.
Here, the generalized Benders decomposition algorithm of \cite{geoffrion1972generalized} is used, which allows for solving nonlinear optimization problems where the subproblem is convex.
For example,
\cite{humpola2017gas} discusses multiple valid inequalities for passive and active gas network optimization problems in a Benders context.
They also propose domain relaxation to reduce the time spent on the nonlinear subproblems.
\cite{yang2017modeling} use a Benders decomposition approach to solve a gas flow problem, solving the nonlinear subproblems via successive linear programming, with the Benders cuts based on  Lagrange multipliers on the nodal balance equations.

Overall, we observe that Benders decomposition is regularly used to solve potential-based flow problems, in particular affine ones.
Most works there focus either on the general design of problem or decomposition, or on cut design.
However, subproblem solving can be a significant contributor to overall running time, particularly when many scenarios need to be considered simultaneously.
While stronger cuts can lead to fewer Benders iterations, and thus fewer subproblems to be solved, a specialized subproblem algorithm may still reduce overall running time.

\section*{Contribution}
In this paper, we present a specialized subproblem solution algorithm as part of a Benders decomposition approach for potential-based flow problems with affine potential functions, switchable arcs and arc-based contingency scenarios.
It is similar in spirit to, \eg, the knapsack algorithm used by \cite{fischetti2017redesigning} or the semi-assignment algorithm by \cite{contreras2011benders}.
In our case, the algorithm can be seen as a specialized simplex algorithm, where we obtain (a sufficiently large part of) an optimal basis and corresponding primal and dual solutions or dual proof of infeasibility in one step \citep[cf.][]{vanderbei1998linear}.
It is split up into two parts.

First, a primal solution is determined to decide if and what kind of cut needs to be generated.
This is based on the observation that the flow of a passive, \ie, switch arc free, affine potential-based flow problem is unique \citep{birkhoff1956potential,koch2015evaluating}
, together with a procedure to convert active potential-based flow problems with switch arcs fixed by the master solution into equivalent passive ones.
The flow is then determined with a system of linear equations, which has the weighted Laplacian of the network as coefficient matrix.

In the second step, the primal solution is used to construct a dual linear system of equations (SLE), which then provides the dual solution.
Using a primal basis, another linear system is constructed.
It shares its coefficient matrix with that of the first step.
Again, this requires special handling of switch arcs.
For optimality cuts, only the right-hand side differs between the two steps, and we show that the resulting dual solution is optimal, and partially unique for the chosen basis.
For feasibility cuts, the matrix is extended by an additional row and column, and the right-hand side changes.
Different feasibility cuts are constructed based on the component ``causing'' the infeasibility.

Similarly to simplex algorithms, our algorithm is extended by a factorization update scheme, where potential similarities between subproblems, rather than bases of the same problem, are exploited via low-rank updates.
These occur between different master solutions and scenarios, but also between the primal and dual systems.
We discuss the theoretical implications of the presented algorithm, concluding that the update scheme results in a significant reduction in theoretical complexity.
We note that factorization update schemes, in particular also involving weighted Laplacian matrices such as the one we use, have been studied in other works as well, see, \eg, \cite{alsac1983sparsity}.

We demonstrate our approach on instances of the lossless DC security-constrained optimal transmission switching problem \citep{hedman_optimal_2009}, where transmission system operators minimize operating costs while considering single-line contingency scenarios.
As we show, any affine potential-based flow problem instance can be transformed into such an instance, so these instances are representative for the larger class.
Computational experiments on 174 instances demonstrate total running time reductions of more than \qty{71.2108}{\percent} in the geometric mean over the standard solving method, which further reduces to \qty{97.9819}{\percent} when considering only the subproblem solve time. 

The paper is organized as follows:
\autoref{sec:pot_flows} provides a detailed problem description of the core potential-based flow problem as well as the additional aspects considered here.
We then outline a branch-and-Benders-cut approach in \autoref{sec:benders}, where each subproblem takes the form of a specific, simpler potential-based flow problem.
For individual subproblems, this structure is exploited in \autoref{sec:single_sp} to arrive at an SLE-based solving algorithm.
The factorization update scheme is detailed in \autoref{sec:sp:comp_fac}, tying together the different subproblems solved as part of the branch-and-Benders-cut algorithm.
Therein, we also briefly discuss the theoretical complexity of the proposed subproblem solving algorithm.
\autoref{sec:evaluation} then provides a computational evaluation.

%% file: chapters/2-potential-flows.tex
\emph{Potential-based flows}, as discussed by \cite{gross2019algorithmic}, describe flow problems that often model real-life physical networks.
We consider directed graphs $(V, \LArcs \cup \SArcs)$ with \emph{nodes} $V$ and \emph{line} $\LArcs \subseteq V \times V$ and \emph{switch arcs} $\SArcs \subseteq V \times V, \LArcs \cap \SArcs = \emptyset$.
The set of all arcs is given by $\Arcs = \LArcs \cup \SArcs$.
Each node $u \in V$ has a variable potential $\Pot{u} \in \R$ with
lower and upper bounds $\PotMin{u} \in \R_{\leq 0}$ and $\PotMax{u}
\in \R_{\geq 0}$, respectively.
Every arc $a = (u,v) \in \LArcs$ has a variable \emph{flow} $\Flow{u}{v} \in \R$ and a \emph{conductance} $\Cond{u}{v} > 0$ which describes its ability to carry flow.
The \emph{potential network} $(G, \CondSym)$ comprises the graph $G$ and arc conductances $\CondSym \in \R^{\LArcs}$.
For each line arc $(u,v) \in \LArcs$, the interdependence between flow $\Flow{u}{v}$ and potentials $\Pot{u}, \Pot{v}$ is given by the \emph{potential equation}
\begin{equation*}
  \Pot{u} - \Pot{v} = \PotF{u}{v}(\Flow{u}{v})\Cond{u}{v}^{-1},
\end{equation*}
where $\PotF{u}{v}: \R \rightarrow \R$ is a \emph{potential function}.
Potential functions are continuous and strictly increasing.
Rearranging yields the flow as a function of the potentials:
\begin{equation*}
  \Flow{u}{v} = \PotF{u}{v}^{-1}((\Pot{u} - \Pot{v}) \Cond{u}{v}).
\end{equation*}
The flow is limited by lower and upper bounds $\FlowMin{u}{v} \in \R_{\leq 0}$ and $\FlowMax{u}{v} \in \R_{\geq 0}$, respectively.

Switch arcs $(u,v)\in\SArcs$ can be switched on or off to connect or disconnect nodes in the graph, denoted by switching state $\z{u}{v} \in \{0,1\}$.
If a switch arc $(u,v) \in \SArcs$ is switched ``on'' ($\z{u}{v} = 1$), the incident node potentials must be equal, $\Pot{u} = \Pot{v}$.
If switched ``off'' ($\z{u}{v} = 0$), the arc flow must be zero, $\Flow{u}{v} = 0$, and the potentials are no longer directly linked.
In the latter case, the difference of the respective potentials is confined to lower and upper limits $\PotDifMin{u}{v} \in \R$ and $\PotDifMax{u}{v} \in \R$, respectively, which are selected such that they can never restrict the choice of potentials.

We refer to switched on and off switch arcs as ``active'' and ``inactive,'' respectively.
We require that there are no two nodes in the same switch arc-connected component in the graph $(V, \SArcs)$ that are incident to the same line arc.
Additionally, for any two such components, there can be at most one line arc that is incident to a node in each component.
\autoref{fig:network_component} shows an example graph with both line and switch arcs.

\begin{figure}[htbp]
  \FIGURE{
  \centering
  \tikzset{
    aux_node/.style = {circle, draw, minimum size=#1},
    aux_node/.default = 11 
  }
  \tikzset{
    coupler_node/.style = {circle, draw, minimum size=#1},
    coupler_node/.default = 11 
  }
  \begin{tikzpicture}
    \node (0,0)                  (center) {};
    \node [left =1.37 of center] (lcenter) {};
    \node [right=1.37 of center] (rcenter) {};
    \node [below=2.532 of center.center] (push) {};
    \node [coupler_node, above=1.2 of center] (coup1) {};
    \node [coupler_node, below=1.2 of center] (coup2) {};
    \node [aux_node, above =.2 of lcenter] (aux1) {};
    \node [aux_node, above =.2 of rcenter] (aux2) {};
    \node [aux_node, below =.2 of lcenter] (aux3) {};
    \node [aux_node, below =.2 of rcenter] (aux4) {};
    \node [aux_node, left  =1 of aux1] (naux1) {};
    \node [aux_node, right =1 of aux2] (naux2) {};
    \node [aux_node, left  =1 of aux3] (naux3) {};
    \node [aux_node, right =1 of aux4] (naux4) {};
    \node [left =.5 of aux1.west] (L1_west) {};
    \node [right=.5 of aux2.east] (L2_east) {};
    \node [left =.5 of aux3.west] (L3_west) {};
    \node [right=.5 of aux4.east] (L4_east) {};
    \draw [<-]  (naux1) -- (aux1);
    \draw [<-] (naux2) -- (aux2);
    \draw [<-] (naux3) -- (aux3);
    \draw [<-] (naux4) -- (aux4);
    \draw [<-, dashed] (aux1)  -- (coup1);
    \draw [<-, dotted] (aux3)  -- (coup1);
    \draw [<-, dotted] (aux1)  -- (coup2);
    \draw [<-, dashed] (aux3)  -- (coup2);
    \draw [<-, dashed] (aux2)  -- (coup1);
    \draw [<-, dotted] (aux4)  -- (coup1);
    \draw [<-, dotted] (aux2)  -- (coup2);
    \draw [<-, dashed] (aux4)  -- (coup2);
    \draw [->, dotted] (coup1) -- (coup2);
    \draw [->, dotted] (aux3)  -- (aux4);
    \node [above=.03 of L1_west.east] {$a^{\ell}_1$};
    \node [above=.03 of L2_east.west] {$a^{\ell}_2$};
    \node [above=.03 of L3_west.east] {$a^{\ell}_3$};
    \node [above=.03 of L4_east.west] {$a^{\ell}_4$};
    \node [above=.01 of coup1.north] {$u_1$};
    \node [below=.01 of coup2.south] {$u_2$};
    \draw (-0.85, 1.48  ) node {$a^s_1$};
    \draw (-0.85, -1.48  ) node {$a^s_2$};
    \node [above=.01 of aux1.north] {$u_3$};
  \end{tikzpicture}
  }
  {Network topology.\label{fig:network_component}}
  {Line arcs are represented by solid lines, active and inactive switch arcs by dashed and dotted lines, respectively. The active switch arcs and incident nodes form switch arc components.}
\end{figure}

Each node has an adjustable \emph{demand} $\Demand{u} \in [\DemandMin{u}, \DemandMax{u}]$, $\DemandMin{u}, \DemandMax{u} \in \R$, $\DemandMin{u} \leq \DemandMax{u}$, with per-unit cost $\DemandCosts{u} \in \R$.
The demands are always balanced, \ie,
\begin{equation*}
    \sum_{u \in V} b_u = 0,
\end{equation*}
due to the flow conservation constraints.
If demands are fixed to a value and there are no switch arcs, $\SArcs = \emptyset$, the network is called \emph{passive}.
Otherwise, it is \emph{active}.

A potential-based flow $(\FlowSym, \PotSym, \zSym, \DemandSym)$ then comprises arc flows $\FlowSym \in \R^{\Arcs}$, node potentials $\PotSym \in \R^{V}$, switching states $\zSym \in \{0,1\}^{\SArcs}$ and demands $\DemandSym \in \R^{V}$.
It is considered feasible if it fulfills \eqref{eq:core:switch_off}--\eqref{eq:core:pot_domain}:

\begin{subequations}
  \begin{alignat}{2}
    \z{u}{v} \FlowMin{u}{v} \leq \Flow{u}{v} & \leq \z{u}{v} \FlowMax{u}{v} && \quad \forall (u,v) \in \SArcs \label{eq:core:switch_off}\\
    (1 - \z{u}{v})\PotDifMin{u}{v} \leq \Pot{u} - \Pot{v} & \leq (1 - \z{u}{v}) \PotDifMax{u}{v} && \quad \forall (u,v) \in \SArcs \label{eq:core:switch_on}\\
    \sum_{(v,u) \in \delta^-(u)} \Flow{v}{u} - \sum_{(u,v) \in \delta^+(u)} \Flow{u}{v} & = \Demand{u} && \quad \forall u \in V \label{eq:core:flow_cons}\\
    \PotF{u}{v}^{-1}((\Pot{u} - \Pot{v}) \Cond{u}{v}) & = \Flow{u}{v} && \quad \forall (u,v) \in \LArcs \label{eq:core:pot_eq}\\
    \Demand{u} & \in [\DemandMin{u}, \DemandMax{u}] && \quad \forall u \in V \label{eq:core:demand_domain}\\
    \z{u}{v} & \in \{0,1\} && \quad \forall (u,v) \in \SArcs \label{eq:core:z_domain}\\
    \Flow{u}{v} & \in [\FlowMin{u}{v}, \FlowMax{u}{v}] && \quad \forall (u,v) \in \Arcs \label{eq:core:flow_domain}\\
    \Pot{u} & \in [\PotMin{u}, \PotMax{u}] && \quad \forall u \in V \label{eq:core:pot_domain}
  \end{alignat}
\end{subequations}
  
\autocsref{eq:core:switch_off} and \eqref{eq:core:switch_on} ensure that switch arcs carry no flow when off and have equal incident potentials when on.
Here, for all switch arcs $(u,v) \in \SArcs$, the ``big-$M$'' parameters $\FlowMin{u}{v}, \FlowMax{u}{v} \in \R$ are chosen such that the bound can never be met while the arc is on, and analogously for $\PotDifMin{u}{v}, \PotDifMax{u}{v} \in \R$.
\autocref{eq:core:flow_cons} enforces flow conservation: the net flow at each node meets its demand.
\autocref{eq:core:pot_eq} is the potential equation relating flows and potentials.
\autocsref{eq:core:demand_domain}--\eqref{eq:core:pot_domain} specify the variables' domains.

\paragraph{Extended potential-based flow problems.}
In addition to this basic po\-ten\-tial-based flow problem, we consider scenario-based robustness, line arc overloading, and switching limits.
Robustness is enforced via a set of \emph{contingencies} $\aC \subseteq 2^{\Arcs}$, where each contingency $\C \in \aC$ represents a scenario in which all contained arcs fail.
The methods discussed in this work extend naturally to arbitrary contingencies with individual bounds, (fixed) demands or conductances, however we focus on (single) arc contingencies for the sake of simplicity.
The set $\aC$ may include the empty set $\emptyset$, representing the base case without failures.
Switching and demand decisions are made independently of which contingency occurs, while flows and potentials adapt to each contingency.
We denote by $\LArcsOut \coloneq \LArcs \setminus \C$ and $\SArcsOut \coloneq \SArcs \setminus \C$, respectively, the sets of line and switch arcs unaffected by contingency $\C$, and by $\ArcsOut \coloneq \LArcsOut \cup \SArcsOut$ the set of all such arcs.

In each contingency $\C \in \aC$, each line arc $(u,v) \in \LArcs$ can be overloaded by $\FlowOverC{u}{v}{\C} \in \R_{\geq 0}$ beyond its soft bounds $\FlowMinOver{u}{v}, \FlowMaxOver{u}{v} \in \R$ at cost $\OverCost \in \R_{\geq 0}$ per unit, within the already introduced hard bounds $\FlowMin{u}{v}$ and $\FlowMax{u}{v}$, \ie, $\FlowMin{u}{v} \leq \FlowMinOver{u}{v} \leq \FlowMaxOver{u}{v} \leq \FlowMax{u}{v}$.
Overloaded arcs contribute to the objective function for each contingency individually.
Alternatively, ways to aggregate subproblem costs, \eg, the maximum over all subproblems, could also be considered.

Each switch arc $(u,v) \in \SArcs$ has a default state.
We collect arcs switched on by default in set $\SArcsOn\subseteq\SArcs$, and in set $\SArcsOff=\SArcs\setminus\SArcsOn$ those which are inactive in the default state.
At most $\TAL \in \N$ switch arcs may be switched away from their defaults.
The graph has to be ``islanding-free'' in all contingencies. That is, all nodes $u \in V$ where demands are not fixed to zero, \ie, $\DemandMin{u} \neq 0$ or $\DemandMax{u} \neq 0$,  have to be part of the same connected component induced by the active arcs (arcs that are switched on and not part of the contingency).
This is a looser requirement than full connectedness in that it allows for parts of the graph to be disconnected by switching or contingencies, as long as only one of them is relevant to the distribution of demands.

All these parameters define a
potential-based flow problem $\mathcal{P}=(G,\CondSym,\PotFSym,\aC,H)$, where 
\begin{equation*}
  H=(\DemandMinSym,\DemandMaxSym,\FlowMinSym,\FlowMaxSym,\FlowMinOverSym,\FlowMaxOverSym,\PotMinSym,\PotMaxSym,\SArcsOff,\TAL,\DemandCostsSym,\OverCost)  
\end{equation*}
collects bounds, default state, and costs.
Should the contingencies be irrelevant in a statement, we omit them from the description of $\mathcal{P}$.

We are ready to state the complete extended model (omitting the connectedness constraints) in \autocsref{eq:pOrig:obj}--\eqref{eq:pOrig:flowover_domain}. 
Arc flows and node potentials now materialize in every scenario.
The corresponding variables are thus superscripted with the respective contingency.

\begin{mini!}
  {}{\sum_{u \in V} \DemandCosts{u} \Demand{u} + \sum_{\C \in \aC} \sum_{(u,v)\in \LArcsOut} \OverCost \FlowOverC{u}{v}{\C} \label{eq:pOrig:obj}}{}{}
  \addConstraint{\forall \C \in \aC:}\nonumber
  \addConstraint{\sum_{(u,v) \in \SArcsOn} \left( 1 - \z{u}{v} \right) + \sum_{(u,v) \in \SArcsOff} \z{u}{v}}{\leq Z \label{eq:pOrig:TAL}}
  \addConstraint{\z{u}{v} \FlowMin{u}{v} \leq \FlowC{u}{v}}{\leq \z{u}{v} \FlowMax{u}{v}}{\quad \forall (u,v) \in \SArcsOut \label{eq:pOrig:switch_off}}
  \addConstraint{(1 - \z{u}{v}) \PotDifMin{u}{v} \leq \PotC{u} - \PotC{v}}{\leq (1 - \z{u}{v}) \PotDifMax{u}{v}}{\quad \forall (u,v) \in \SArcsOut \label{eq:pOrig:switch_on}}
  \addConstraint{\sum_{(v,u) \in \InArcs{u}} \FlowC{v}{u} - \sum_{(u,v) \in \OutArcs{u}} \FlowC{u}{v}}{= \Demand{u}}{\quad \forall u \in V \label{eq:pOrig:flow_cons}}
  \addConstraint{\PotF{u}{v}^{-1}((\PotC{u} - \PotC{v}) \Cond{u}{v})}{= \FlowC{u}{v}}{\quad \forall (u,v) \in \LArcsOut \label{eq:pOrig:pot_eq}}
  \addConstraint{\FlowMinOver{u}{v} - \FlowOverC{u}{v}{\C} \leq \FlowC{u}{v}}{\leq \FlowMaxOver{u}{v} + \FlowOverC{u}{v}{\C}}{\quad \forall (u,v) \in \LArcsOut \label{eq:pOrig:overload}}
  \addConstraint{\FlowMin{u}{v} \leq \FlowC{u}{v}}{\leq \FlowMax{u}{v}}{\quad \forall (u,v) \in \LArcsOut \label{eq:pOrig:flow_bound}}
  \addConstraint{\PotMin{u} \leq \PotC{u}}{\leq \PotMax{u}}{\quad \forall u \in V \label{eq:pOrig:pot_bound}}
  \addConstraint{\Demand{u}}{\in [\DemandMin{u}, \DemandMax{u}]}{\quad \forall u \in V \label{eq:pOrig:demand_domain}}
  \addConstraint{\z{u}{v}}{\in \{0,1\}}{\quad \forall (u,v) \in \SArcs}
  \addConstraint{\FlowC{u}{v}}{\in \R}{\quad \forall (u,v) \in \ArcsOut}
  \addConstraint{\PotC{u}}{\in \R}{\quad \forall u \in V}
  \addConstraint{\FlowOverC{u}{v}{\C}}{\geq 0}{\quad \forall (u,v) \in \LArcsOut \label{eq:pOrig:flowover_domain}}
\end{mini!}

The objective function minimizes total demand and overload costs across all contingencies.
\autocref{eq:pOrig:TAL} enforces the switching limit on deviations from default settings.
The core potential-based flow \autocsref{eq:core:switch_off}--\eqref{eq:core:pot_eq} are replicated for each contingency via \autocsref{eq:pOrig:switch_off}--\eqref{eq:pOrig:pot_eq}.
\autocref{eq:pOrig:overload} tracks line arc overloading beyond soft bounds.
\autocsref{eq:pOrig:flow_bound} and \eqref{eq:pOrig:pot_bound} enforce the bounds on flows and potentials for each contingency.
\autocsref{eq:pOrig:demand_domain}--\eqref{eq:pOrig:flowover_domain} specify bounds and domains of the variables.

\paragraph{Affine Potential Functions.}
\label{sec:aff_pot_func}
In this work, we focus on instances with affine potential functions with positive slope:
\begin{equation*}
    \PotF{u}{v}(x) = \PotFSlope{u}{v} x + \PotFConst{u}{v},\, \forall x \in \R, (u,v) \in \LArcs,
\end{equation*}
where $\PotFSlope{u}{v} \in \R_{>0}$ and $\PotFConst{u}{v} \in \R$.
By standard transformation, any such function can be replaced by the identity function:

\begin{proposition}
  Let $\mathcal{P}=(G,\CondSym,\PotFSym,H)$ be a potential-based flow problem with affine potential functions $\PotFSym$.
  Then, there exists an equivalent potential-based flow problem $\mathcal{P}'=(G,\CondSym',\PotFSym',H')$ where every potential function is the identity function $\text{id}:\R \rightarrow\R$, $\text{id}(x) = x$ such that any feasible solution $(\FlowSym', \PotSym, \zSym, \DemandSym')$ to $\mathcal{P}'$ can be transformed into a feasible solution $(\FlowSym, \PotSym, \zSym, \DemandSym)$ to $\mathcal{P}$ with the same objective function value, and vice versa.\label{the:aff_id}
\end{proposition}
\proof{Proof.}
  We transform the original instance via arc parameter rescaling and demand adjustment.
  For each $(u,v) \in \LArcs$, inverting the potential function in \autocref{eq:core:pot_eq} yields
  \begin{align*}
      \Flow{u}{v} &= \PotF{u}{v}^{-1}((\Pot{u} - \Pot{v}) \Cond{u}{v})\\
                  &= \frac{(\Pot{u} - \Pot{v}) \Cond{u}{v} - \PotFConst{u}{v}}{\PotFSlope{u}{v}}\\
                  &= (\Pot{u} - \Pot{v}) \frac{\Cond{u}{v}}{\PotFSlope{u}{v}} - \frac{\PotFConst{u}{v}}{\PotFSlope{u}{v}}.
  \end{align*}
  The rescaled conductance is
  \begin{equation*}
    \Cond{u}{v}' \coloneq \frac{\Cond{u}{v}}{\PotFSlope{u}{v}},
  \end{equation*}
  yielding
  \begin{align*}
      \Flow{u}{v}' \coloneq& (\Pot{u} - \Pot{v}) \Cond{u}{v}' \\
                   =& (\Pot{u} - \Pot{v})  \frac{\Cond{u}{v}}{\PotFSlope{u}{v}} - \frac{\PotFConst{u}{v}}{\PotFSlope{u}{v}} + \frac{\PotFConst{u}{v}}{\PotFSlope{u}{v}} \\
                   =& \Flow{u}{v} + \frac{\PotFConst{u}{v}}{\PotFSlope{u}{v}}.
  \end{align*}
  The shift $\PotFConst{u}{v}/\PotFSlope{u}{v}$ is absorbed into all flow bounds, for instance, the transformed lower bound becomes $\FlowMin{u}{v}' = \FlowMin{u}{v} + \PotFConst{u}{v}/\PotFSlope{u}{v}$.
  These shifts also affect flow conservation \autocsref{eq:core:flow_cons}. Substituting
  \begin{equation}
    \Flow{u}{v} = \Flow{u}{v}' - \frac{\PotFConst{u}{v}}{\PotFSlope{u}{v}} \label{eq:shifted_flow}
  \end{equation}
  gives
  \begin{align*}
      \Demand{u} &= \sum_{(v,u) \in \delta^-(u)} \Flow{v}{u} - \sum_{(u,v) \in \delta^+(u)}  \Flow{u}{v}\\
      &= \sum_{(v,u) \in \delta^-(u)}  \Flow{v}{u}' - \sum_{(u,v) \in \delta^+(u)}  \Flow{u}{v}'\\
      &\quad - \sum_{(v,u) \in \delta^-(u)} \frac{\PotFConst{v}{u}}{\PotFSlope{v}{u}} + \sum_{(u,v) \in \delta^+(u)} \frac{\PotFConst{u}{v}}{\PotFSlope{u}{v}}.
  \end{align*}
  Rearranging yields
  \begin{equation*}
    \sum_{(v,u) \in \delta^-(u)}  \Flow{v}{u}' - \sum_{(u,v) \in \delta^+(u)}  \Flow{u}{v}' = \Demand{u}',
  \end{equation*}
  with transformed demand
  \begin{equation}
    \Demand{u}' \coloneq \Demand{u} + \sum_{(v,u) \in \delta^-(u)} \frac{\PotFConst{v}{u}}{\PotFSlope{v}{u}} - \sum_{(u,v) \in \delta^+(u)} \frac{\PotFConst{u}{v}}{\PotFSlope{u}{v}} \label{eq:shifted_demand}
  \end{equation}
  for all $u \in V$.
  We shift the demand bounds accordingly
  \begin{equation*}
    \DemandMin{u}' = \DemandMin{u} + \sum_{(u,v) \in \delta^-(u)} \frac{\PotFConst{u}{v}}{\PotFSlope{u}{v}} - \sum_{(u,v) \in \delta^+(u)} \frac{\PotFConst{u}{v}}{\PotFSlope{u}{v}},
  \end{equation*}
  analogously for $\DemandMax{u}'$, and add the term
  \begin{equation*}
    -\DemandCosts{u} \cdot \left( \sum_{(v,u) \in \delta^-(u)} \frac{\PotFConst{v}{u}}{\PotFSlope{v}{u}} - \sum_{(u,v) \in \delta^+(u)} \frac{\PotFConst{u}{v}}{\PotFSlope{u}{v}} \right)
  \end{equation*}
  to the objective function.

  To summarize, the transformed instance is obtained from the original instance by,
  for all line arcs $(u,v) \in \LArcs$,
  replacing the conductance $\Cond{u}{v}$ by the rescaled $\Cond{u}{v}'$,
  shifting flow bounds $\FlowMin{u}{v}$, $\FlowMax{u}{v}$, $\FlowMinOver{u}{v}$, and $\FlowMaxOver{u}{v}$ by $\PotFConst{u}{v}/\PotFSlope{u}{v}$, and
  replacing its potential function with the identity potential function
  \begin{equation*}
      \Flow{u}{v}' = \Cond{u}{v}' (\Pot{u} - \Pot{v}) \quad \forall (u,v) \in \LArcs;
  \end{equation*}
  and for all nodes $u \in V$,
  shifting its demand bounds $\DemandMin{u}, \DemandMax{u}$ and adding a term to the objective to match the adjusted demand $\Demand{u}$.
  The transformed instance is equivalent to the original instance: following Equations \eqref{eq:shifted_flow} and \eqref{eq:shifted_demand}, any feasible solution to the transformed instance can be transformed into a feasible solution to the original instance with the same objective function value, and vice versa.
  \Halmos
\endproof

In the following, we can thus assume all potential functions to be the identity function. 
As such, we can substitute all flow variables on line arcs by a term over potential variables following the potential equation.
For ease of notation, we still refer to flows where helpful, but line arc flows and potential equations are not explicitly part of the following models.

Problem \eqref{eq:core:switch_off}--\eqref{eq:core:pot_domain} is NP-hard if potential functions are the identity function \citep{lehmann2014complexity}, and thus following \autoref{the:aff_id} also if they are affine.

%% file: chapters/3-benders.tex
We decompose the problem \autocsref{eq:pOrig:obj}--\eqref{eq:pOrig:flowover_domain} into a master problem (MP) and one subproblem (SP) per contingency $\C \in \aC$.
The MP initially consists only of demand $\DemandSym$ and switching $\zSym$ variables as well as the constraint on switching actions \eqref{eq:pOrig:TAL} and constraint
\begin{equation*}
  \sum_{u \in V} \Demand{u} = 0  
\end{equation*}
to enforce global demand balance.
Additionally, it has one Benders recourse variable $\eta_\C \geq 0$ for each contingency $\C \in \aC$ that is added to the objective function, reflecting a lower bound estimate of that contingency's objective function contribution for any master solution.

Subproblems contain the remaining constraints and variables, with master variables as fixed parameters.
Let $(\DemandFixedSym, \zFixedSym, \etaFixed_\C)$, $\DemandFixedSym \in \R^{V}$, $\zFixedSym \in \{0,1\}^{\SArcs}$, $\etaFixed_\C \in \R$ denote a master solution with values for demands, topology, and the Benders variable of contingency $\C \in \aC$, respectively.
We denote the corresponding $\C$-th subproblem as $\subP$, and its optimal objective function value, if feasible, as $\subP^* \in \R$.
The subproblem $\subP$ is a potential-based flow problem where contingency arcs $(u,v) \in \C$ are removed, and the master problem decisions are fixed. 
We call these networks with fixed demands and switch arc states \emph{semi-passive}.
The subproblem model for $\C \in \aC$ is given in Model \eqref{eq:sp:obj}--\eqref{eq:sp:flowover_domain}.
For reference, we indicate the dual variables (or ``Farkas multipliers,'' respectively) corresponding to the respective constraints in square brackets.

\begin{subequations}
  \begin{alignat}{4}
    & \mathrlap{\subP^* =} \nonumber\\
    \min\quad & \mathrlap{\sum_{(u,v)\in \LArcsOut} \OverCost \FlowOverC{u}{v}{\C}} &&&&&& \label{eq:sp:obj}\\
    \subto\quad & & \zFixed{u}{v} \FlowMin{u}{v} \leq \FlowC{u}{v} & \leq \zFixed{u}{v} \FlowMax{u}{v} &\quad& [\DualSOffCLbSym, \DualSOffCUbSym] &\quad& \forall (u,v) \in \SArcsOut \label{eq:sp:switch_off}\\
     & & (1 - \zFixed{u}{v}) \PotDifMin{u}{v} \leq \PotC{u} - \PotC{v} & \leq (1 - \zFixed{u}{v}) \PotDifMax{u}{v} &\quad& [\DualSOnCLbSym,\DualSOnCUbSym] &\quad& \forall (u,v) \in \SArcsOut \label{eq:sp:switch_on}\\
     & & \sum_{(v,u) \in \delta^-(u) \cap \LArcsOut} (\PotC{v} - \PotC{u}) \Cond{v}{u} & \mathrel{\phantom{=}} &\quad&  &\quad&   \nonumber\\
     & & {} - \sum_{(u,v) \in \delta^+(u) \cap \LArcsOut} (\PotC{u} - \PotC{v}) \Cond{u}{v} & \mathrel{\phantom{=}} &\quad&  &\quad&   \nonumber\\
     & & {} + \sum_{(v,u) \in \delta^-(u) \cap \SArcsOut} \FlowC{v}{u} & \mathrel{\phantom{=}} &\quad&  &\quad&   \nonumber\\
     & & {} - \sum_{(u,v) \in \delta^+(u) \cap \SArcsOut} \FlowC{u}{v} & = \DemandFixed{u} &\quad& [\DualFCCSym] &\quad& \forall u \in V \label{eq:sp:flow_cons}\\
     & & \FlowMinOver{u}{v} - \FlowOverC{u}{v}{\C} \leq (\PotC{u} - \PotC{v}) \Cond{u}{v} & \leq \FlowMaxOver{u}{v} + \FlowOverC{u}{v}{\C} &\quad& [\DualOverCLbSym,\DualOverCUbSym] &\quad& \forall (u,v) \in \LArcsOut \label{eq:sp:overload}\\
     & & \FlowMin{u}{v} \leq (\PotC{u} - \PotC{v}) \Cond{u}{v} & \leq \FlowMax{u}{v} &\quad& [\DualFlowBCLbSym,\DualFlowBCUbSym] &\quad& \forall (u,v) \in \LArcsOut \label{eq:sp:flow_bound}\\
     & & \PotMin{u} \leq \PotC{u} & \leq \PotMax{u} &\quad& [\DualPotBCLbSym,\DualPotBCUbSym] &\quad& \forall u \in V \label{eq:sp:pot_bound}\\
     & & \FlowC{u}{v} & \in \R &\quad&  &\quad& \forall (u,v) \in \SArcsOut\\
     & & \PotC{u} & \in \R &\quad&  &\quad& \forall u \in V\\
     & & \FlowOverC{u}{v}{\C} & \in \R_{\geq 0} &\quad&  &\quad& \forall (u,v) \in \LArcsOut \label{eq:sp:flowover_domain}
  \end{alignat}
\end{subequations}

The decomposed problem is solved using a branch-and-Benders-cut algorithm \citep[cf.][]{rahmaniani2017benders}.
The MP is solved as usual using a branch-and-bound algorithm.
Every time an integer feasible MP solution is obtained during the algorithm, it needs to be verified with all subproblems, and is only valid if it is ``accepted'' by all of them.
For this, each subproblem is solved with the master solution's values as parameters.
If an SP is infeasible, a feasibility cut is added to the MP.
This feasibility cut is based on the values of a dual direction of unboundedness (``Farkas multipliers''), and cuts off at least the master solution causing the infeasibility.
With $\DualSOffCLbOptSym$ and $\DualSOffCUbOptSym$ denoting dual direction values of the lower and upper bound, respectively, of \autocref{eq:sp:switch_off}, and analogously for the other constraints, a feasibility cut is
\begin{equation}
  \begin{aligned}
    & \sum_{(u,v) \in \SArcsOut} \Bigl( \z{u}{v} \FlowMax{u}{v} \DualSOffCUbOpt{u}{v} - \z{u}{v} \FlowMin{u}{v} \DualSOffCLbOpt{u}{v} \\
    & \hphantom{\sum_{(u,v) \in \SArcsOut} \Bigl(} + (1 - \z{u}{v}) \PotDifMax{u}{v} \DualSOnCUbOpt{u}{v} - (1 - \z{u}{v}) \PotDifMin{u}{v} \DualSOnCLbOpt{u}{v} \Bigr) \\
    & + \sum_{(u,v) \in \LArcsOut} \left( \FlowMaxOver{u}{v} \DualOverCUbOpt{u}{v} - \FlowMinOver{u}{v} \DualOverCLbOpt{u}{v} + \FlowMax{u}{v} \DualFlowBCUbOpt{u}{v} - \FlowMin{u}{v} \DualFlowBCLbOpt{u}{v}\right) \\
    & + \sum_{u \in V} \left( \Demand{u} \DualFCCOpt{u} + \PotMin{u} \DualPotBCLbOpt{u} + \PotMax{u} \DualPotBCUbOpt{u}\right) \leq 0.
  \end{aligned}\label{eq:feas_cut}
\end{equation}
Otherwise, if an SP is feasible, and the optimal objective function value of $\subP$ is strictly larger than the master solution value for $\eta_\C$, an optimality cut is added to the MP.
This cut is based on the values of an optimal dual solution.
Denoting these the same way as for the feasibility cut, the cut inequality only differs in the right-hand side, here $\eta_\C$, \ie, the Benders recourse variable.
By adding these cuts, changes in switching actions and demands in the MP reducing line arc overloads in contingencies are thereby incentivized.
If neither cut is added in any subproblem, the master solution is accepted as a valid solution in the current node of the branch-and-bound tree. 

%% file: chapters/4-sle.tex
In this section, we focus on solving a single subproblem $\subP$ after obtaining the master solution, \ie, fixed values for demand $\DemandSym$, topology $\zSym$ and $\eta_\C$.
We discuss how a constant number of systems of linear equations (SLE) are sufficient to solve subproblems, \ie, check for and compute optimality and feasibility cuts.
The next section then integrates this as part of the larger branch-and-Benders-cut algorithm, where multiple subproblems are solved in succession.

We assume that the graph is islanding-free.
That is, all nodes with demands not fixed to zero by their bounds are connected in the graph induced by the active arcs, \ie, those switched-on and not part of a given contingency.
This can be enforced through, \eg, a separate separation procedure.
We verify connectedness by computing a minimum cut on the contracted graph (see upcoming section) with line arcs not in the contingency contracted as well.
Additionally, all nodes that are not part of the connected component containing these demands are removed with their incident arcs, as we can assume that they do not carry flow (or potentials) in an optimal solution.

The subproblem solving algorithm is motivated by the comparatively high computational effort of the simplex algorithm---it requires solving a potentially very high number of linear equation systems, whereas our algorithm requires solving only one or two of equal or smaller size.
It is roughly divided into two parts:
first we solve a reduced version of the primal subproblem, free of soft and hard variable bounds, the objective function, and switch arcs.
From this, we obtain flow and potential variable values for all arcs and nodes.
The potential differences as well as relevant flow values of a switch arc-free network are unique for a given master solution \citep{birkhoff1956potential}.
With these, we can thus check whether the subproblem is feasible, and, if so, compute its objective function value.
This determines whether the subproblem yields an optimality or a feasibility cut, or neither.
If it does yield a cut, we solve a reformulated version of the dual subproblem as an SLE to obtain all necessary dual variable values to construct the respective cut.
While similar, this dual SLE differs between optimality and feasibility cuts, as detailed later.

The remainder of this section is structured as follows:
in Subsection \ref{sec:sp:primal_sb}, we discuss an arc contraction procedure for fixed switch arcs to obtain a switch arc-free graph yielding an equivalent line arc flow.
Subsection \ref{sec:sp:primal} discusses the first part of the algorithm, \ie, how to construct the primal SLE from a given master solution.
Subsections \ref{sec:sp:opt_cut} and \ref{sec:sp:feas_cut} then explain the second part of the algorithm, \ie, how to construct the dual SLE for an optimality or feasibility cut, respectively, from a given master and primal SLE solution.
The resulting improvements in computational complexity are discussed jointly in the next section.

\subsection{Switch Arc Contraction}\label{sec:sp:primal_sb}
With all switching variables fixed in the subproblem, and thus a semi-passive network, we can simplify the underlying network by substituting all switch arcs in a way that preserves the flow on line arcs.
This is done by ``contracting'' switch arcs iteratively, eventually resulting in a graph with only line arcs.
The uniqueness of flows with regard to a passive network was observed by, \eg, \cite{birkhoff1956potential}, so the contraction procedure can be seen as a precondition for the uniqueness of flows in our subproblems.
Similarly, contracted graphs are used for constructing the dual SLEs. 
For both primal and dual solutions, we first determine most variable values on the contracted graph, and then any remaining values as necessary, as discussed in the following sections.

We first show that we can disregard inactive switch arcs.
For ease of readability, we omit any notational reference to the contingency.
Then, we describe the contraction procedure that eliminates active switch arcs.
\begin{proposition}
  Let $\mathcal{P}=((V, \LArcs \cup \SArcs),\CondSym,\PotFSym,H)$ be a potential-based flow problem with a semi-passive network fixed by $\DemandFixedSym  \in \R^V$ and $\zFixedSym \in \SArcs$. Let $\SArcszNC = \{(u,v) \in \SArcs \mid \zFixed{u}{v} = 1\}$ be the set of active switch arcs. Any optimal solution to the identical problem without inactive switch arcs $\mathcal{P}'=((V, \LArcs \cup \SArcszNC),\CondSym,\PotFSym,H)$ can be transformed into an optimal solution to $\mathcal{P}$. \label{the:inactive_switch_arcs}
\end{proposition}
\proof{Proof.}
  Let $(\PotSym^*, \FlowSym^{*\prime})$, $\PotSym \in \R^V$, $\FlowSym^{*\prime} \in \R^{\SArcszNC}$, be an optimal solution to  $\mathcal{P}'$. Then, $(\PotSym^*, \FlowSym^*)$, where
  \begin{equation*}
    \Flow{u}{v}^* = \begin{cases}
                        \Flow{u}{v}^{*\prime}, & (u,v) \in \SArcszNC\\
                        0, & \text{else}
                    \end{cases}
  \end{equation*}
  is a feasible solution to $\mathcal{P}$:
  all constraints except for \autocsref{eq:sp:switch_off} and \eqref{eq:sp:switch_on} for all $(u,v) \notin \SArcszNC$ as well as \autocsref{eq:sp:flow_cons} of the incident nodes are unaffected and thus trivially fulfilled.
  Let $(u,v) \mathrel{\notin} \SArcszNC$ be an inactive arc.
  Since it is inactive, \ie, $\Flow{u}{v}^* = 0$, \autocsref{eq:sp:switch_off} are fulfilled with equality, and \autocsref{eq:sp:switch_on} are always fulfilled by construction.
  \autocsref{eq:sp:flow_cons} for the incident nodes $u, v \in V$ sum over the additional variable $\Flow{u}{v}$, but since $\Flow{u}{v}^* = 0$, the result does not change. 

  The solution $(\PotSym^*, \FlowSym^*)$ is optimal for $\mathcal{P}$, otherwise, there must be a solution $(\PotSym^{*\prime\prime}, \FlowSym^{*\prime\prime})$ for $\mathcal{P}$ with improved objective function value and $\Flow{u}{v}^{*\prime\prime} = 0$ for all $(u,v) \notin \SArcszNC$.
  But then, $(\PotSym^{*\prime\prime}, \FlowSym^{*\prime\prime})$ reduced to arcs in $\SArcszNC$ is clearly feasible for $\mathcal{P}'$, so $(\PotSym^*, \FlowSym^{*\prime})$ can not be optimal for $\mathcal{P}'$.
  \Halmos
\endproof

\begin{figure*}[htbp]
  \FIGURE{
  \subcaptionbox{After contracting arc $a^s_1$.}{
  \centering
  \tikzstyle{aux_node} = [diamond, draw]
  \tikzset{
    coupler_node/.style = {circle, draw, minimum size=#1},
    coupler_node/.default = 11 
  }
  \begin{tikzpicture}
    \node (0,0)                  (center) {};
    \node [left =1.37 of center] (lcenter) {};
    \node [right=1.37 of center] (rcenter) {};
    \node [below=2.532 of center.center] (push) {};
    \node [coupler_node, above=1.2 of center] (coup1) {};
    \node [coupler_node, below=1.2 of center] (coup2) {};
    \node [coupler_node, draw=white, above =.2 of lcenter] (aux1) {};
    \node [coupler_node, above =.2 of rcenter] (aux2) {};
    \node [coupler_node, below =.2 of lcenter] (aux3) {};
    \node [coupler_node, below =.2 of rcenter] (aux4) {};
    \node [left =-.07 of coup1.center] (mlaux1) {};
    \node [left =-.07 of coup2.center] (mlaux2) {};
    \node [right =-.07 of coup1.center] (mraux1) {};
    \node [right =-.07 of coup2.center] (mraux2) {};
    \node [coupler_node, left  =1 of aux1] (naux1) {};
    \node [coupler_node, right =1 of aux2] (naux2) {};
    \node [coupler_node, left  =1 of aux3] (naux3) {};
    \node [coupler_node, right =1 of aux4] (naux4) {};
    \node [left =.5 of aux1.west] (L1_west) {};
    \node [right=.5 of aux2.east] (L2_east) {};
    \node [left =.5 of aux3.west] (L3_west) {};
    \node [right=.5 of aux4.east] (L4_east) {};  
    \draw [<-] (naux1) -- (coup1);
    \draw [<-] (naux2) -- (aux2);
    \draw [<-] (naux3) -- (aux3);
    \draw [<-] (naux4) -- (aux4);
    \draw [<-, dashed] (aux3)  -- (coup2);
    \draw [<-, dashed] (aux2)  -- (coup1);
    \draw [<-, dashed] (aux4)  -- (coup2);
    \node [coupler_node, above=1.2 of center, fill=white] (coup1) {};
    \node [coupler_node, below=1.2 of center, fill=white] (coup2) {};
    \node [above=.01 of L2_east.west] {$a^{\ell}_2$};
    \node [above=.01 of L3_west.east] {$a^{\ell}_3$};
    \node [above=.01 of L4_east.west] {$a^{\ell}_4$};
    \node [above=.01 of coup1.north] {$u_1$};
    \node [below=.01 of coup2.south] {$u_2$};
    \draw (-1.65, 1.48) node {$a^{\ell}_1$};
  \end{tikzpicture}
  }\hspace{1.5cm}
  \subcaptionbox{Full switch arc contraction.}{
  \centering
  \tikzstyle{aux_node} = [diamond, draw]
  \tikzset{
    coupler_node/.style = {circle, draw, minimum size=#1},
    coupler_node/.default = 11 
  }
  \begin{tikzpicture}
    \node (0,0)                  (center) {};
    \node [left =1.37 of center] (lcenter) {};
    \node [right=1.37 of center] (rcenter) {};
    \node [below=2.532 of center.center] (push) {};
    \node [coupler_node, above=1.2 of center] (coup1) {};
    \node [coupler_node, below=1.2 of center] (coup2) {};
    \node [coupler_node, draw=white, above =.2 of lcenter] (aux1) {};
    \node [coupler_node, draw=white, above =.2 of rcenter] (aux2) {};
    \node [coupler_node, draw=white, below =.2 of lcenter] (aux3) {};
    \node [coupler_node, draw=white, below =.2 of rcenter] (aux4) {};
    \node [coupler_node, left  =1 of aux1] (naux1) {};
    \node [coupler_node, right =1 of aux2] (naux2) {};
    \node [coupler_node, left  =1 of aux3] (naux3) {};
    \node [coupler_node, right =1 of aux4] (naux4) {};
    \draw [<-] (naux1) -- (coup1);
    \draw [<-] (naux2) -- (coup1);
    \draw [<-] (naux3) -- (coup2);
    \draw [<-] (naux4) -- (coup2);
    \draw (-1.65, 1.48) node {$a^{\ell}_1$};
    \draw (1.65, 1.48) node {$a^{\ell}_2$};
    \draw (-1.65, -1.48) node {$a^{\ell}_3$};
    \draw (1.65, -1.48) node {$a^{\ell}_4$};
    \node [above=.01 of coup1.north] {$u_1$};
    \node [below=.01 of coup2.south] {$u_2$};
  \end{tikzpicture}
  }
  }
  {Switch arc component contraction. \label{fig:busbar_contract}}
  {Contraction of switch arc components from \autoref{fig:network_component} after removing inactive switch arcs. Solid lines represent line arcs, dashed lines active switch arcs. In (a), the arc $a^\ell_1$ is now directly connected to $u_1$. In (b), the components have fully been merged into the representative nodes $u_1$ and $u_2$.}
\end{figure*}

For a specific contingency $C$ denote by $\SArcsz \subseteq \{(u,v) \in \SArcsOut \mid \zFixed{u}{v} = 1\}$ the set of remaining, active switch arcs in the network.
Let $\Subs \subseteq 2^V$ denote the set of connected components in the subgraph $(V, \SArcsz)$, the \emph{active switch arc components}.
By \autocsref{eq:sp:switch_on}, the nodes of any such component share the same potential.
To determine the potential of all nodes, it thus suffices to consider each connected component in the graph as a single representative node, the potential of which represents the potential of all members of the component.
We iteratively replace each arc $(u,v) \in \SArcsz$ and its incident nodes in the component $\Sub \in \Subs$ with a single contracted node:
first, we remove $(u,v)$ from $\SArcsz$.
Then, we replace all arcs $(u',v), (v,v') \in \LArcsOut \cup \SArcsz$ that are incident to $v$ with $(u',u)$ and $(u,v')$, respectively (while keeping their other properties the same).
If a node is adjacent to both $u$ and $v$, \ie, replacing the arc incident to $v$ would result in two parallel arcs, we do not add the second arc.
As no line arc can be incident to two nodes of the same switch arc-connected component, such an arc can only be a switch arc.
Lastly, we remove $v$ from $V$ and set the demand of the contracted node to $\DemandFixed{u} + \DemandFixed{v}$.

We denote the remaining node of component $\Sub$ by $\SubNode{\Sub}$.
The set of (line) arcs in the contracted graph is denoted by $\LArcsCC$.
Refer to \autoref{fig:busbar_contract} for examples.

Solving the contracted problem is sufficient to obtain an optimal primal solution to the original problem.
As before, our notation does not refer to the contingency in the result and proof.

\begin{theorem}
  Let $\mathcal{P}=((V, \LArcs \cup \SArcs),\CondSym,\PotFSym,H)$ be a potential-based flow problem with a semi-passive network fixed by $\DemandFixedSym  \in \R^V$ and $\zFixedSym \in \SArcs$. Let $\mathcal{P}'=((V', \LArcs' \cup \SArcs'),\CondSym',\PotFSym',H')$, $\SArcs' = \emptyset$, be the potential-based flow problem on the corresponding passive contracted network. Any optimal solution $(\PotSym^{*\prime}, \FlowSym^{*\prime})$ for $\mathcal{P}'$ can be transformed directly into an optimal solution $(\PotSym^*, \FlowSym^*)$ for $\mathcal{P}$, and vice versa.\label{the:contraction}
\end{theorem}
\proof{Proof.}
  Following \autoref{the:inactive_switch_arcs}, we can disregard all inactive switch arcs of $\mathcal{P}$.
  For the transformed solution, for all switch arc components $\Sub \in \SubsNC$, we set the potential of all contained nodes to that of the representative node,
  \begin{equation*}
    \Pot{u}^* = \Pot{u(\Sub)}^{*\prime} \quad \forall u \in \Sub, \Sub \in \SubsNC.
  \end{equation*}
  The switch arc flows $\FlowOptSym$ need to be assigned values such that the flow conservation \autocsref{eq:sp:flow_cons} are fulfilled.
  Let again $\SArcszNC = \{(u,v) \in \SArcs \mid \zFixed{u}{v} = 1\}$ denote the set of active switch arcs.
  We determine a spanning tree on the subgraph $(\Sub, \SArcszNC \cap (\Sub \times \Sub))$ for each component $\Sub \in \SubsNC$.
  All arcs of that subgraph not in the spanning tree are assigned a flow of zero, in analogy to non-basic variables in a network simplex algorithm.
  Then, starting from the leaves of the spanning tree until all variables are determined, there always is a node in $\Sub$ where the flow conservation equation has only one undetermined variable, the flow of an incident switch arc.
  It is easy to see that this solution is feasible and has the same objective function value in $\mathcal{P}$, the flows on line arcs remain unchanged.

  Conversely, we show that for any feasible solution of $\mathcal{P}$, the contraction procedure preserves the flow on line arcs, \ie, that a solution with its variables adapted as described remains feasible in $\mathcal{P}'$. 
  Let $\FlowOptSym \in \R^{\LArcs \cup \SArcszNC}$ denote an optimal flow solution, with line arc flows obtained following the potential equations.
  Let $(u,v) \in \SArcszNC$ be the first arc we contract (into node $u$).
  By \autocsref{eq:sp:flow_cons}, it holds that
  \begin{equation}
    \smashoperator[r]{\sum_{(v,u) \in \InArcsNC{u}}} \FlowOpt{v}{u} - \smashoperator[r]{\sum_{(u,v) \in \OutArcsNC{u}}} \FlowOpt{u}{v} = \DemandFixed{u}\label{eq:subst_nbu}
  \end{equation}
  and
  \begin{equation}
    \sum_{(u,v) \in \InArcsNC{v}} \FlowOpt{u}{v} - \sum_{(v,u) \in \OutArcsNC{v}} \FlowOpt{v}{u} = \DemandFixed{v}.\label{eq:subst_nbv}
  \end{equation}
  The adapted \autocref{eq:sp:flow_cons} at the contracted node $u$ is
  \begin{equation*}
    \begin{aligned}
      &\sum_{(v,u) \in \InArcsNC{u}} \Flow{v}{u} - \sum_{(u,v) \in \OutArcsNC{u}} \Flow{u}{v}\\
      +&\sum_{(u,v) \in \InArcsNC{v}} \Flow{u}{v} - \sum_{(v,u) \in \OutArcsNC{v}} \Flow{v}{u}\\ &\hspace{2.7cm}= \DemandFixed{u} + \DemandFixed{v},
    \end{aligned}
  \end{equation*}
  which is the sum of Equations (\ref{eq:subst_nbu}) and (\ref{eq:subst_nbv}), and thus the solution fulfills the new constraint.
  It is clear that all other constraints remain fulfilled, and that the adapted solution is still optimal.
  By induction, this holds after repeating the procedure on the remaining arcs in $\SArcszNC$ as well.
  \Halmos
\endproof

Any semi-passive network can thus be simplified to a passive network, and the results by \citet{birkhoff1956potential} apply.

\subsection{Primal Solving}\label{sec:sp:primal}
In this subsection, we construct the ``primal'' SLE, and show how it is used to obtain a full primal solution.
We first show how this SLE is constructed from a master solution, and provide a proof showing that it always has a unique solution \citep{birkhoff1956potential}.
This unique solution then allows us to determine if and which cut needs to be generated.
Additionally, certain of the resulting primal variable solution values are needed to construct the dual SLEs in the following sections.

Following \autoref{the:inactive_switch_arcs} and \autoref{the:contraction}, we are able to omit all switch arcs. 
For ease of notation, we assume without loss of generality that the graph is already switch arc free.
Then, we can simplify the formulation given in \autocsref{eq:sp:obj}--\eqref{eq:sp:flowover_domain} to the one in \eqref{eq:sp_simp:obj}--\eqref{eq:sp_simp:flowover_domain} by removing \autocsref{eq:sp:switch_off} and \eqref{eq:sp:switch_on} and variables $\FlowC{u}{v}$ for $(u,v) \in \SArcsOut$.

\begin{subequations}
  \begin{alignat}{4}
    & \mathrlap{\subP^* =} \nonumber\\
    \min\quad & \mathrlap{\sum_{(u,v)\in \LArcsOut} \OverCost \FlowOverC{u}{v}{\C}} &&&&&& \label{eq:sp_simp:obj}\\
    \subto\quad & & \sum_{(v,u) \in \InArcs{u}} (\PotC{v} - \PotC{u}) \Cond{v}{u} & \mathrel{\phantom{=}} &\quad&  &\quad&   \nonumber\\
     & & {} - \sum_{(u,v) \in \OutArcs{u}} (\PotC{u} - \PotC{v}) \Cond{u}{v} & = \DemandFixed{u} &\quad& [\DualFCCSym] &\quad& \forall u \in V \label{eq:sp_simp:flow_cons}\\
     & & \FlowMinOver{u}{v} - \FlowOverC{u}{v}{\C} \leq (\PotC{u} - \PotC{v}) \Cond{u}{v} & \leq \FlowMaxOver{u}{v} + \FlowOverC{u}{v}{\C} &\quad& [\DualOverCLbSym,\DualOverCUbSym] &\quad& \forall (u,v) \in \LArcsOut \label{eq:sp_simp:overload}\\
     & & \FlowMin{u}{v} \leq (\PotC{u} - \PotC{v}) \Cond{u}{v} & \leq \FlowMax{u}{v} &\quad& [\DualFlowBCLbSym,\DualFlowBCUbSym] &\quad& \forall (u,v) \in \LArcsOut \label{eq:sp_simp:flow_bound}\\
     & & \PotMin{u} \leq \PotC{u} & \leq \PotMax{u} &\quad& [\DualPotBCLbSym,\DualPotBCUbSym] &\quad& \forall u \in V \label{eq:sp_simp:pot_bound}\\
     & & \PotC{u} & \in \R &\quad&  &\quad& \forall u \in V\\
     & & \FlowOverC{u}{v}{\C} & \in \R_{\geq 0} &\quad&  &\quad& \forall (u,v) \in \LArcsOut \label{eq:sp_simp:flowover_domain}
  \end{alignat}
\end{subequations}
\autocsref{eq:sp_simp:flow_cons} of the simplified model then yield a system of linear equations on the potential variables.
We call it the \emph{primal} SLE and denote it by $\lesP$.
With $\InOutArcs{u}$ as the set of in- and outgoing arcs of node $u \in V$, the coefficient matrix of the primal SLE is given by the (negated) weighted Laplacian matrix of the graph, \ie,
\begin{equation}
  D^{\zFixedSym, \C}_{u,v} =  \begin{cases}
                  -\sum_{(u,v) \in \InOutArcs{u}} \Cond{u}{v},& u=v,\\
                  \Cond{u}{v},& (u,v) \in \LArcsOut,\\
                  \Cond{v}{u},& (v,u) \in \LArcsOut,\\
                  0,& \text{else},
                \end{cases}\label{eq:lesp_matrix}
\end{equation}
for all $u, v \in V$:
for \autocsref{eq:sp_simp:flow_cons} for node $u \in V$, variable $\Pot{u}$ occurs in the first sum for each incoming arc $(u,v) \in \InArcs{u}$ with coefficient $-\Cond{u}{v}$, and in the second sum for each outgoing arc $(u,v) \in \OutArcs{u}$ with coefficient $\Cond{u}{v}$, but the sum itself has a negative coefficient.
Similarly, each node $v \in V$ adjacent to $u$ occurs exactly once with the conductance of the connecting arc as positive coefficient.

For passive potential-based flow problems such as this one, it is known that the flow, obtained from the solution of the SLE following the potential equation, is unique \citep{birkhoff1956potential}.
As stated by \cite{pfetsch_potential-based_2026}, we can thus first obtain a primal solution to $\subP$ by solving the SLE, and then check its adherence to the bounds afterward.
This covers the soft bounds given by \autocsref{eq:sp_simp:overload} that are penalized in the objective function, and the hard bounds given by \autocsref{eq:sp_simp:flow_bound} and \eqref{eq:sp_simp:pot_bound}.
We show a proof for the uniqueness here as we give similar proofs later on, and as it is integral to our work---otherwise we might derive a suboptimal solution or falsely conclude that the subproblem is infeasible.

\begin{theorem}\label{the:unique_flows}
  Any solution to the primal SLE $D^{\zFixedSym, \C} \PotCSym = \DemandFixedSym$ yields the same flow. 
\end{theorem}
\proof{Proof.}
  We select as \emph{reference node} an arbitrary node, say node $w \in V$, such that removing it does not split the graph into disconnected components.
  Such a node always exists, for example any leaf of a spanning tree of the graph.
  We then set its potential $\PotC{w}$ to some arbitrary value $\PotSym^{\text{ref}} \in~\R$.
  This value can indeed be chosen arbitrarily at first, as, due to the linear potential function, the flows are unaffected by linear shifts in potentials, \ie, $\Pot{1} - \Pot{2} = (\Pot{1} + \varepsilon) - (\Pot{2} + \varepsilon)$ for any $\Pot{1}, \Pot{2}, \varepsilon \in \R$.

  We then remove the entries corresponding to the reference node in the SLE to obtain the remaining coefficient matrix $D^{\zFixedSym, \C, \text{ref}}$, with
  \begin{equation*}
    D^{\zFixedSym, \C, \text{ref}}_{u,v} = D^{\zFixedSym, \C}_{u,v} \qquad \forall u,v \in V \setminus \{w\}.
  \end{equation*}
  The matrix is identical to the original matrix, but without the row and column corresponding to the reference node.
  In particular, the remaining diagonal entries remain unchanged.

  From standard linear algebra, 
  we know that any SLE $D x = b$, $D \in \R^{m \times n}, b \in \R^m, m \in \N$ with $n \in \N$ variables has a unique solution exactly if $\text{rank}(D) = n$.
  In our remaining SLE, 
  $D^{\zFixedSym, \C, \text{ref}}$ needs to have rank $|V|-1$.
  This is the case if
  \begin{enumerate}
    \item it is weakly diagonally dominant, \ie,
      \begin{equation}
        \label{eq:weakly_diagonally_dominant}
        \left|D^{\zFixedSym, \C, \text{ref}}_{u,u}\right| \geq \sum_{v \in V \setminus \{u, w\}} \left|D^{\zFixedSym, \C, \text{ref}}_{u,v}\right| \quad \forall u \in V \setminus \{w\},
      \end{equation} 
    \item 
      \autoeqref{eq:weakly_diagonally_dominant} holds with strict inequality for at least one $u \in V \setminus \{w\}$, and \label{text:unq_dds}
    \item $D^{\zFixedSym, \C, \text{ref}}$ is irreducible \citep{JohnsonMarijuánPisonero+2023}. 
  \end{enumerate}
  All of these hold:
  \begin{enumerate}
    \item This is true due to the matrix structure, and since all conductances are positive.
          For all  $u, v\in V \setminus \{w\}$:
          \begin{equation*}
            D^{\zFixedSym, \C, \text{ref}}_{u,v} =  \begin{cases}
                            -\sum_{(u,v) \in \InOutArcs{u}} \Cond{u}{v},& u=v,\\
                            \Cond{u}{v},& (u,v) \in \LArcsOut,\\
                            \Cond{v}{u},& (v,u) \in \LArcsOut,\\
                            0,& \text{else}.
                          \end{cases}
          \end{equation*}
          Then for all $u \in V \setminus \{w\}$:
          \begin{align}
            \begin{split}
              \left|D^{\zFixedSym, \C, \text{ref}}_{u,u}\right| & =   \left|- \sum_{(u,v) \in \InOutArcs{u}}            \Cond{u}{v}\right|\\
                                        & =    \sum_{(u,v) \in \InOutArcs{u}}           \left|\Cond{u}{v}\right|\\
                                        & \geq \sum_{\substack{(u,v) \in \InOutArcs{u}\\ \{u,v\} \not\ni w}} \left|\Cond{u}{v}\right|\\
                                        & =    \sum_{v \in V \setminus \{w, u\}}                    \left|D^{\zFixedSym, \C, \text{ref}}_{u,v}  \right|.
            \end{split}\label{eq:diag_dom}
          \end{align}
    \item Since the graph is connected, there must be at least one node $u \in V$ that is connected to the reference node, so $(u,w) \in \LArcsOut$ or $(w,u) \in \LArcsOut$.
          For each such $u$,
          \begin{equation*}
            (\InOutArcs{u}) \supsetneq (\InOutArcs{u}) \setminus \{(u,w), (w,u)\}.
          \end{equation*}
          Thus, since all conductances are positive, the inequality in \autoeqref{eq:diag_dom} holds strictly.
    \item The matrix is reducible exactly if there exists a partition $U,U'$ of $V$ such that $D^{\zFixedSym, \C, \text{ref}}_{u,u'} = 0$ for all $u \in U$, $u' \in U'$.
          But then the graph would be disconnected (or it is only connected at node $w$, but $w$ was chosen such that removing it does not split the graph).
  \end{enumerate}
  
  \Halmos
\endproof
Since any solution to $\lesP$ determines a unique flow and the objective function solely depends on the flow, this flow must coincide with the optimal flow in $\subP$, provided the solution is feasible.
Consider a solution $\PotCOptSym \in \R^{V}$ to $\lesP$.
From \autocsref{eq:sp_simp:overload}, we have that for any optimal solution,
\begin{align*}
  \FlowOverCOpt{u}{v} = \max\{&- (\PotCOpt{u} - \PotCOpt{v}) \Cond{u}{v} + \FlowMinOver{u}{v},\\
  & \qquad(\PotCOpt{u} - \PotCOpt{v}) \Cond{u}{v} - \FlowMaxOver{u}{v}, 0\}
\end{align*}
for all $(u,v) \in \LArcsOut$, \ie, if overloaded, the overloading variables match the violation of the soft bounds.
The optimal objective function value to $\subP$ is given by
\begin{equation*}
  \subP^* = \sum_{(u,v)\in \LArcsOut} \OverCost \FlowOverCOpt{u}{v}.
\end{equation*}
Adherence to \autocsref{eq:sp_simp:flow_bound} can be immediately determined from the solution.
For \autocsref{eq:sp_simp:pot_bound}, it is based on the choice of reference potential.
They are necessarily violated exactly if
\begin{align*}
    & \max\{0, \max\{\PotCOpt{u} - \PotMax{u} \mid u \in V\}\} \\
  > & \max\{0, \min\{\PotCOpt{u} - \PotMin{u} \mid u \in V\}\}, \text{ or}\\
    & \max\{0, \min\{\PotMax{u} - \PotCOpt{u} \mid u \in V\}\} \\
  < & \max\{0, \max\{\PotMin{u} - \PotCOpt{u} \mid u \in V\}\},
\end{align*}
\ie, if potentials can not be shifted enough to satisfy the largest violation of an upper bound without violating any lower bound, or vice versa.

This is sufficient to decide on the type of cut we construct:
If $\subP^* > \etaFixed_\C$, and the subproblem is feasible, we construct an optimality cut.
The subproblem is infeasible if any of \autocsref{eq:sp_simp:flow_bound} or \eqref{eq:sp_simp:pot_bound} are violated, in which case we construct a feasibility cut.
Otherwise, if the subproblem is feasible, and $\subP^* \leq \etaFixed_\C$, the master solution is accepted by this subproblem.
Note that, unlike when using a simplex algorithm, we do not have a dual solution to the subproblem. This comes next.

\subsection{Optimality Cuts}\label{sec:sp:opt_cut}
In this subsection, we construct a separate SLE, the solution of which allows us to easily construct an optimality cut.
We consider the uncontracted subproblem with switch arcs in this section.
We first determine an optimal primal basis based on our optimal primal solution, and then use that basis to construct the ``dual'' SLE and show that it yields an optimal dual solution.

\paragraph{Optimal primal basis.}
The primal solution potentially induces multiple optimal primal bases for $\subP$ in augmented form, as the potential values can be additively shifted, as discussed, and the flow on cyclic switch arc components is not necessarily uniquely determined.
We deterministically select one of them as follows, and then use it to compute an optimal dual solution.
The basis needs to contain $m \in \N$ variables, where $m$ is the number of constraints of $\subP$, including all variables that are non-zero.
This includes the additional slack variables that we need to introduce to convert the inequality constraints to equality constraints.
The columns of the selected variables need to be linearly independent.
Following the proof of Theorem \ref{the:unique_flows}, we know that this does not hold directly for the potential variables.
As in the proof, we thus substitute the potential variable $\PotC{w}$ of the reference node $w \in V$ by the reference value $\PotSym^\text{ref} \in \R$ and remove the corresponding constraint.
With this, there are $m = 4|\SArcsOut| + 3|V|-1 + 4|\LArcsOut| $ constraints.

For the basis $B$, we select one variable per constraint, always selecting the ones that are non-zero for a solution derived from $\PotCOptSym$:
\begin{itemize}
  \item We select both slack variables of \autocsref{eq:sp:switch_off}.
  \item For \autocsref{eq:sp:switch_on} for arc $(u, v) \in \SArcsOut$, both slack variables are zero if the arc is on, else both are potentially non-zero. 
        Of all the switch arcs that are on, we choose a maximal subset $\SArcs^{\zFixedSym,\C\prime} \subseteq \SArcsz$ such that the underlying undirected subgraph $(V, \{\{u,v\} \mid (u,v) \in \SArcs^{\zFixedSym,\C\prime}\})$ is a forest (\ie, for each switch arc component, we select only arcs corresponding to an arbitrary undirected spanning tree), and select the corresponding flow variables $\FlowCSym$ for the basis.
        We select the upper bound slack variables of the remaining switch arcs, as well as all lower bound slack variables.
  \item Up to $|V|-1$ $\PotCSym$ variables are non-zero as a result of $\lesP$ and \autocsref{eq:sp:flow_cons}.
        We select all $\PotCSym$ variables for the basis (excluding that of the reference node).
  \item For \autocsref{eq:sp:overload} for arc $(u,v) \in \LArcsOut$, either $\FlowOverCOpt{u}{v}>0$, in which case at least one of the slack variables is zero, or $\FlowOverCOpt{u}{v}=0$.
  In the latter case, we select both slack variables, in the former the non-zero one and $\FlowOverC{u}{v}{\C}$.
  \item We select both slack variables of \eqref{eq:sp:flow_bound} and \eqref{eq:sp:pot_bound}, they can be zero or non-zero.
\end{itemize}

These variables form a valid basis for the Model~\eqref{eq:sp:obj}--\eqref{eq:sp:flowover_domain}, derived from our optimal solution.
We select exactly $m$ variables, one for each constraint.
This includes all non-zero variables, as only zero valued variables are discarded.
The corresponding columns are also linearly independent:
the columns of all slack variables, \ie, those for Constraints \eqref{eq:sp:switch_off}, \eqref{eq:sp:flow_bound}, and \eqref{eq:sp:pot_bound}; and those we selected for Constraints \eqref{eq:sp:switch_on} and \eqref{eq:sp:overload} are unit vectors.
This also holds for any selected overload variable $\FlowOverC{u}{v}{\C}, (u,v) \in \LArcs$.
These unit vectors are trivially linearly independent and can be used to remove any non-zero entries in the rows corresponding to the respective constraints. Then, only the coefficients of \autocsref{eq:sp:switch_on} where the arc is switched on and \autocsref{eq:sp:flow_cons} remain.
Together, they form a $(|\SArcsz| + |V|-1) \times (|\SArcsz| + |V|-1)$ coefficient matrix of some SLE with the right-hand side according to that of the constraints.
Its variables are the potential variables of all nodes except for the reference node, as well as the flow variables of all selected switch arcs.
As discussed in Subsection \ref{sec:sp:primal}, the solution of the contracted graph is uniquely determined, and immediately gives us the potentials of all nodes of the original graph. The potential of the nodes that were contracted is equal to that of their representative node.
The flows on all selected switch arcs are then uniquely determined following the nodal balance equations \eqref{eq:sp:flow_cons}, as discussed for the proof of \autoref{the:contraction}:
since the underlying undirected graph of each component is always a spanning tree, the nodal balance equation at each leaf node of such a tree always only has at most one undetermined variable, the flow on an incident switch arc.
The same then holds for at least one additional node, until eventually all flow variables are uniquely determined.
In total, the SLE has a unique solution, and thus the columns of its coefficient matrix must be linearly independent, as desired.
Hence, our selection represents a valid primal basis for our optimal solution.

\paragraph{Dual system of linear equations.}
We use the basis to construct a system of linear equations to obtain an optimal dual solution for Model~\eqref{eq:sp:obj}--\eqref{eq:sp:flowover_domain}.
For this, we first dualize the primal subproblem containing only the variables in our basis.
For brevity, we do not explicitly show the full form of this dual subproblem here, as it is not of further interest.
Instead, we discuss how it can be simplified to obtain the system of linear equations shown in \eqref{eq:sp_dual:pot}--\eqref{eq:sp_dual:fcc_domain}, which then is reduced to the SLE \eqref{eq:lesd_sb_merge}, the SLE that is solved during the algorithm.
We also keep the reference node's variable and constraints for ease of notation, removing them is only needed for uniqueness arguments.

First, all dual constraints of non-slack primal variables can be treated as equality constraints, since we want to obtain an optimal solution:
$\FlowCSym$ and $\PotCSym$ are unrestricted, so their dual constraints are equality constraints immediately from dualization.
The variable $\FlowOverCSym$ is only chosen if positive, so its dual constraint must be an equality constraint following complementary slackness.

Second, we can disregard most of the dual variables, as they are equal to $0$ or not relevant to the remaining solution.
We keep the dual variables for flow conservation \eqref{eq:sp:flow_cons} $\DualFCC{u}$ for all nodes $u \in V$, switch arc potential linking \eqref{eq:sp:switch_on} $\DualSOnC{u}{v}$ for all active switch arcs $(u,v) \in \SArcsz$ (including those that are not in $\SArcs^{\zFixedSym,\C\prime}$), and overloading constraints \eqref{eq:sp:overload} for all overloaded line arcs $\DualOverCLb{u}{v}$, $(u,v) \in \LArcsOut$, $\FlowCOpt{u}{v} < \FlowMinOver{u}{v}$, as well as $\DualOverCUb{u}{v}$, $(u,v) \in \LArcsOut$, $\FlowCOpt{u}{v} > \FlowMaxOver{u}{v}$.

All other dual variables are assumed to be equal to zero:
from complementary slackness, we know that in an optimal solution, this must be the case for any dual variable where the corresponding slack variable is positive.
For \autocsref{eq:sp:switch_off}, if the arc is off, the dual variable has objective function coefficient $0$ and does not appear in any other constraints (the flow variable is not in the basis), so it can be chosen arbitrarily. 
If the arc is on, the slack variables are positive by construction.
For \autocsref{eq:sp:switch_on}, if the arc is off, the slack variables are positive by construction.
For the remaining dual variables of \autocsref{eq:sp:overload}, \eqref{eq:sp:flow_bound}, and \eqref{eq:sp:pot_bound}, some constraints can be fulfilled with equality.
We assume for now that the corresponding dual variables can always be chosen as zero in an optimal solution, which we show in the proof of \autoref{the:dual_optimal}.

For ease of notation, we denote for all line arcs $(u,v) \in \LArcs$
\begin{align*}
  \DualPECLb{u}{v} &\coloneq (\DualFCC{u} - \DualFCC{v} - \DualOverCLb{u}{v})\Cond{u}{v},\\
  \DualPECUb{u}{v} &\coloneq (\DualFCC{u} - \DualFCC{v} + \DualOverCUb{u}{v})\Cond{u}{v},\\
  \DualPEC{u}{v} &\coloneq (\DualFCC{u} - \DualFCC{v})\Cond{u}{v},\\
  \intertext{and}
  \hat{\DualPESym}^\C_{(u, v)} &\coloneq \begin{cases}
                                            \DualPECLb{u}{v}, & \FlowCOpt{v}{u} < \FlowMinOver{v}{u}\\
                                            \DualPECUb{u}{v}, & \FlowCOpt{v}{u} > \FlowMaxOver{v}{u}\\
                                            \DualPEC{u}{v},   & \text{else}.
                                          \end{cases}
\end{align*}
The resulting system of linear equations can then be seen in \eqref{eq:sp_dual:pot}--\eqref{eq:sp_dual:fcc_domain}.

\begin{subequations}
  \begin{alignat}{2}
    \sum_{(v,u) \in \InArcs{u}} \hat{\DualPESym}^\C_{(v, u)} - \sum_{(u,v) \in \OutArcs{u}} \hat{\DualPESym}^\C_{(u, v)} & \nonumber\\
    {} + \sum_{(v,u) \in \delta^-(u) \cap \SArcsz} \DualSOnCUb{v}{u} & \nonumber\\
    {} - \sum_{(u,v) \in \delta^+(u) \cap \SArcsz} \DualSOnCUb{u}{v} & = 0 && \quad\spdual{[\PotCSym]}\quad \forall u \in V \label{eq:sp_dual:pot}\\
    -\DualOverCLb{u}{v} & = \OverCost && \quad\spdual{[\FlowOverCSym]}\quad \forall (u,v) \in \LArcsOut, \FlowCOpt{u}{v} < \FlowMinOver{u}{v} \label{eq:sp_dual:over_lb}\\
    -\DualOverCUb{u}{v} & = \OverCost && \quad\spdual{[\FlowOverCSym]}\quad \forall (u,v) \in \LArcsOut, \FlowCOpt{u}{v} > \FlowMaxOver{u}{v} \label{eq:sp_dual:over_ub}\\
    - \DualFCC{u} + \DualFCC{v} & = 0 && \quad\spdual{[\FlowCSym]}\quad \forall (u,v) \in \SArcs^{\zFixedSym,\C\prime} \label{eq:sp_dual:sarc}\\
    \DualSOnCUb{u}{v} & \in \R && \quad\spdual{}\quad \forall (u,v) \in \SArcsz \\
    \DualOverCLb{u}{v} & \in \R && \quad\spdual{}\quad \forall (u,v) \in \LArcsOut, \FlowCOpt{u}{v} < \FlowMinOver{u}{v}\\
    \DualOverCUb{u}{v} & \in \R && \quad\spdual{}\quad \forall (u,v) \in \LArcsOut, \FlowCOpt{u}{v} > \FlowMaxOver{u}{v}\\
    \DualFCC{u} & \in \R && \quad\spdual{}\quad \forall u \in V \label{eq:sp_dual:fcc_domain}
  \end{alignat}
\end{subequations}

This SLE has $|\SArcsz| - |\SArcs^{\zFixedSym,\C\prime}|$ more variables than equations, so we can not always immediately obtain a unique solution.
Similarly to the last section, we avoid this by first calculating a solution without switch arcs.
We first describe the procedure for this and then show that it always provides an optimal solution.

From \autocsref{eq:sp_dual:over_lb} and \eqref{eq:sp_dual:over_ub}, the values of $\DualOverCLb{u}{v}$ and $\DualOverCUb{u}{v}$ are directly determined for all line arcs:
\begin{align*}
  \DualOverCLbOpt{u}{v} &= \begin{cases}
                            0, & \FlowCOpt{u}{v} \geq \FlowMinOver{u}{v} \\
                            -\OverCost, & \FlowCOpt{u}{v} < \FlowMinOver{u}{v},
                          \end{cases}\\
  \DualOverCUbOpt{u}{v} &= \begin{cases}
                            0, & \FlowCOpt{u}{v} \leq \FlowMaxOver{u}{v} \\
                            -\OverCost, & \FlowCOpt{u}{v} > \FlowMaxOver{u}{v}
                          \end{cases}\\
  \intertext{for all $(u,v) \in \LArcsOut$.
We unify this as}
  \DualOverCOpt{u}{v}   &= \begin{cases}
                            0, & \FlowMinOver{u}{v} \leq \FlowCOpt{u}{v} \leq \FlowMaxOver{u}{v} \\
                            -\OverCost, & \FlowCOpt{u}{v} < \FlowMinOver{u}{v} \\
                            \OverCost, & \FlowCOpt{u}{v} > \FlowMaxOver{u}{v}
                          \end{cases}
\end{align*}
for all $(u,v) \in \LArcsOut$.
Note that this term is fixed based on the primal solution.
Observe that the remaining dual problem is almost a semi-passive potential-based flow with identity potential functions itself---potentials $\DualFCCSym$, conductances $\CondSym$, switch flows $\DualSOnCUbSym$, and demands $\DualOverCOptSym$ with $\CondSym$.
However, the objective function does not fit with our definition of a potential-based flow problem.

Nevertheless, from \autoref{the:contraction}, we know that the contraction procedure preserves ``flows'' for both line $\DualPECSym = (\DualFCC{u} - \DualFCC{v})\Cond{u}{v}$, $(u,v) \in \LArcsOut$, and switch arcs $\DualSOnCUb{u}{v}$, $(u,v) \in \SArcsz$.
As before, let $\Sub \subseteq V$ be a set of all the nodes that are contracted to a common node, \ie, the nodes of a connected component in $(V, \SArcsz)$, and let $\Subs$ denote the set of all such sets.
Summing up Equation (\ref{eq:sp_dual:pot}) for all $u \in \Sub$ allows us to remove the $\DualSOnCUbSym$ variables, it yields
\begin{align}
  &\sum_{u \in \Sub} \left( \sum_{(v,u) \in \InArcs{u}} \DualPEC{v}{u} - \DualOverCOpt{v}{u} \Cond{v}{u} - \sum_{(u,v) \in \OutArcs{u}} \DualPEC{u}{v} - \DualOverCOpt{u}{v} \Cond{u}{v}\right) = 0 \nonumber 
\intertext{
(each switch arc in $\SArcsz \cap (\Sub \times \Sub)$, \ie, that is part of the component, occurs exactly once as incoming and outgoing arc each).
We can move the fixed terms containing $\DualOverCOptSym$ to the right-hand side
}
  &\sum_{u \in \Sub} \left( \sum_{(v,u)\in \InArcs{u}} \DualPEC{v}{u} - \sum_{(u,v)\in \OutArcs{u}} \DualPEC{u}{v} \right) \nonumber\\
  &\qquad = \sum_{u \in \Sub} \left( \sum_{(v,u)\in \InArcs{u}} \DualOverCOpt{v}{u} \Cond{v}{u} - \sum_{(u,v) \in \OutArcs{u}} \DualOverCOpt{u}{v} \Cond{u}{v}\right). \nonumber
  \label{eq:dual_sb_leq}
\end{align}
We denote the right-hand side term as
\begin{align}
  \DualOverCOptSym_\Sub \coloneq \sum_{u \in \Sub} \left( \sum_{(v,u)\in \InArcs{u}} \DualOverCOpt{v}{u} \Cond{v}{u} - \sum_{(u,v)\in \OutArcs{u}} \DualOverCOpt{u}{v} \Cond{u}{v} \right).
\end{align}
This yields the SLE
  \begin{align}
    \sum_{u \in \Sub} \left( \sum_{(v,u)\in \InArcs{u}} \DualPEC{v}{u} - \sum_{(u,v)\in \OutArcs{u}} \DualPEC{u}{v} \right) = \DualOverCOptSym_\Sub \quad \forall \Sub \in \Subs.
  \label{eq:lesd}
  \end{align}
This SLE may still have more variables than equations.
We mitigate this by merging the $\DualFCCSym$ variables of each component onto its representative node in the contracted graph.
This preserves the optimality of our computed dual solution, as discussed for \autoref{the:dual_optimal}.
We denote the representative node in the contracted graph of component $\Sub \in \Subs$ by $\SubNode{\Sub}$.
The resulting dual SLE, denoted $\lesD$ is
  \begin{align}
    \sum_{(\SubNode{\OtherSub},\SubNode{\Sub}) \in \delta^-(\SubNode{\Sub})} \DualPEC{\SubNode{\OtherSub}}{\SubNode{\Sub}} - \sum_{(\SubNode{\Sub},\SubNode{\OtherSub}) \in \delta^+(\SubNode{\Sub})} \DualPEC{\SubNode{\Sub}}{\SubNode{\OtherSub}} = \DualOverCOptSym_{\Sub} \quad \forall \Sub \in \Subs.
  \label{eq:lesd_sb_merge}
  \end{align}
The sets $\delta^-$ and $\delta^+$ refer to the contracted graph here.
Observe that, for any graph, only the right-hand side changes between dual SLE $\lesD$ for the optimality cut and the primal SLE $\lesP$ (cf. \autoeqref{eq:lesp_matrix}).

For switch arcs, it remains to calculate values for the $\DualSOnCUbSym$ variables following \autoeqref{eq:sp_dual:pot}.
This is a separate flow problem for each component, with nodal surpluses derived from the dual line flows given by $\DualFCCOptSym$ and $\CondSym$.
As discussed for the primal switch arc flows, we solve it by constructing a spanning tree for each of the components of a given substation, and then moving the surpluses ``up'' starting from the leaves.
If an (undirected) component contains cycles, the solution values for these dual variables are not unique.
However, their objective function coefficient is $0$, so they do not directly influence the objective function value.
 
It remains to show that the computed dual solution is optimal.
\begin{theorem}\label{the:dual_optimal}
  Any solution to the SLE $\lesD$, extended by the tree-based recovery of the dual variables $\DualSOnCUbSym$ and setting the dual variables corresponding to the remaining bounds to zero, yields an optimal dual solution to the original $\subP$. 
\end{theorem}
\proof{Proof.}
  We need to show that the system produces a unique solution for the dual ``flows'' $\DualPECSym$, and that reference value, contraction procedure, and disregarding part of the slack variables do not lead to a suboptimal solution.

  From \autoref{the:unique_flows}, we know that the SLE yields a unique flow $\DualPECOptSym \in \LArcsCC$ on the contracted graph, with unique values for $\DualFCCOptSym \in \Subs$ given a chosen reference value and node.
  Following \autocsref{eq:sp_dual:sarc}, it holds that $\DualFCCOpt{u} = \DualFCCOpt{\SubNode{\Sub}}$ for all $\Sub \in \Subs, u \in \Sub$, as required by the contraction procedure.

  However, unlike the potential variables, the dual variable for node $u \in V$, including the reference node, appears in the objective function of the dual subproblem with potentially non-zero coefficient $\DemandFixed{u} \in \R$.
  We must hence show that the choice of reference value has no impact on the objective function value.
  We required that the demand in the network is balanced, \ie, that its sum over all nodes is equal to $0$.
  Let $\varepsilon \in \R$.
  Then,
  \begin{equation*}
    \begin{aligned}
      \sum_{u \in V} \DualFCCOpt{u} \cdot \DemandFixed{u} & = \sum_{u \in V} \DualFCCOpt{u} \cdot \DemandFixed{u} + \varepsilon \cdot 0\\
                                                            & = \sum_{u \in V} \DualFCCOpt{u} \cdot \DemandFixed{u} + \varepsilon \cdot \sum_{u \in V} \DemandFixed{u}\\
                                                            & = \sum_{u \in V} (\DualFCCOpt{u} + \varepsilon) \cdot \DemandFixed{u},
    \end{aligned}
  \end{equation*}
  as desired.

  Lastly, assume, for a contradiction, that it is invalid to disregard the dual variables of \autocsref{eq:sp:overload}, \eqref{eq:sp:flow_bound}, and \eqref{eq:sp:pot_bound}.
  That is, there is some constraint where the flow or potential is exactly equal to a (soft or hard) bound, and the dual variable of that constraint cannot be equal to zero in any optimal solution.
  Consider a relaxation that differs from the original LP \eqref{eq:sp:obj}--\eqref{eq:sp:flowover_domain} only in loosening exactly this bound by some $\varepsilon \in \R_{> 0}$ (\ie, increasing the upper bound or decreasing the lower bound, as applicable).
  As discussed, the potentials of the primal problem are uniquely determined independently of this bound; thus, relaxing it impacts neither the feasibility nor the objective function value of the primal solution.
  Since the bound is now no longer active,
  the dual variable vanishes by complementary slackness.
  This dual solution is feasible for the dual of the original LP, because changing a right-hand side of the primal does not affect the dual feasible region.
  Since the changed dual variable is at value zero, its modified objective function coefficient (the loosened bound) does not affect the dual objective function value.
  Therefore, by weak duality, our dual solution (of the same objective function value as the primal) must also be optimal.
  This contradicts our assumption that there is no such solution.
  We conclude that we can set all remaining dual variables corresponding to \autocsref{eq:sp:overload}, \eqref{eq:sp:flow_bound}, and \eqref{eq:sp:pot_bound} to zero.
\Halmos
\endproof 

\subsection{Feasibility Cuts}\label{sec:sp:feas_cut}
In this section, we adapt $\lesD$ for the case that $\subP$ is infeasible, and we thus need to construct a feasibility cut. To do so, we need to find a dual direction of unboundedness, \ie, a direction with positive objective function value, of which we can add any positive multiple to a feasible dual point without violating any constraints.
There are two cases in which the subproblem is infeasible:
an arc flow (resulting from the bound-free primal subproblem) violates its hard bounds in \autocsref{eq:sp_simp:flow_bound}, or a node potential violates its hard bounds in \autocsref{eq:sp_simp:pot_bound}.
Both cases can occur simultaneously, also for multiple arcs or nodes.
We add exactly one cut, prioritizing flows over potentials.

\paragraph{Flow bound violation.}
If a flow violates some hard bound, we know that there must be a line arc $(w, u) \in \LArcsOut$ between nodes $w,u \in V$ and components $\Sub', \OtherSub' \in \Subs$, $w \in \Sub', u \in \OtherSub'$, with $(\PotC{w} - \PotC{u}) \Cond{w}{u} < \FlowMin{w}{u}$ or $(\PotC{w} - \PotC{u}) \Cond{w}{u} > \FlowMax{w}{u}$.
We discuss only the latter case,
the former is analogous.
The dual variable of the violated constraint for the contracted graph is $\DualFlowBCUb{\SubNode{\Sub'}}{\SubNode{\OtherSub'}}$ with objective function coefficient $\FlowMax{\SubNode{\Sub'}}{\SubNode{\OtherSub'}} / \Cond{\SubNode{\Sub'}}{\SubNode{\OtherSub'}}$.
Increasing the solution value of $\DualFlowBCUb{\SubNode{\Sub'}}{\SubNode{\OtherSub'}}$ increases the objective function value of the dual solution.
However, we need to do so in a way that does not impact the feasibility of the solution.
The variable occurs only in \autocsref{eq:sp_dual:pot}, for both endpoints of the arc.
We modify the SLE in \autoeqref{eq:lesd} to
\begin{align}
  &\sum_{(\SubNode{\Sub},\SubNode{\OtherSub}) \in \delta^-(\SubNode{\Sub})} \DualPEC{\SubNode{\Sub}}{\SubNode{\OtherSub}} - \sum_{(\SubNode{\OtherSub}, \SubNode{\Sub}) \in \delta^+(\SubNode{\Sub})} \DualPEC{\SubNode{\OtherSub}}{\SubNode{\Sub}} = 0 \enspace \forall \Sub \in \Subs\setminus\{\Sub',\OtherSub'\}, \nonumber\\
  &\DualFlowBCUb{\SubNode{\Sub'}}{\SubNode{\OtherSub'}} + \sum_{(\SubNode{\Sub'}, \SubNode{\OtherSub}) \in\delta^-(\SubNode{\Sub'})} \DualPEC{\SubNode{\Sub'}}{\SubNode{\OtherSub}} - \sum_{(\SubNode{\OtherSub},\SubNode{\Sub'}) \in\delta^+(\SubNode{\Sub'})} \DualPEC{\SubNode{\OtherSub}}{\SubNode{\Sub'}} = 0, \nonumber\\
  &- \DualFlowBCUb{\SubNode{\Sub'}}{\SubNode{\OtherSub'}} + \sum_{(\SubNode{\OtherSub'},\SubNode{\OtherSub}) \in \delta^-(\SubNode{\OtherSub'})} \DualPEC{\SubNode{\OtherSub'}}{\SubNode{\OtherSub}} - \sum_{(\SubNode{\OtherSub},\SubNode{\OtherSub'}) \in \delta^+(\SubNode{\OtherSub'})} \DualPEC{\SubNode{\OtherSub}}{\SubNode{\OtherSub'}} = 0.
  \label{eq:feas_eqs_flow}
\end{align}
Any solution that satisfies this SLE can be added onto a dual feasible point without leaving the dual feasible region.
It thus is a feasible direction for the dual subproblem.
The direction also needs to have a positive objective function value (here, \eg, $1$), so we add this equation to \eqref{eq:feas_eqs_flow}
\begin{align*}
  &\sum_{\Sub \in \Subs} \DemandFixed{\SubNode{\Sub}} \cdot \DualFCC{\SubNode{\Sub}} + \sum_{(\SubNode{\Sub},\SubNode{\OtherSub}) \in \LArcsCC} \left( (\DualOverCLbOpt{\SubNode{\Sub}}{\SubNode{\OtherSub}} + \DualOverCUbOpt{\SubNode{\Sub}}{\SubNode{\OtherSub}}) \cdot \FlowMaxOver{\SubNode{\Sub}}{\SubNode{\OtherSub}} \right)\\
  &+ \frac{\FlowMax{\SubNode{\Sub'}}{\SubNode{\OtherSub'}}}{\Cond{\SubNode{\Sub'}}{\SubNode{\OtherSub'}}} \cdot \DualFlowBCUb{\SubNode{\Sub'}}{\SubNode{\OtherSub'}} = 1.
\end{align*}
Since the primal is infeasible and the dual is feasible, a dual direction of unboundedness must exist.
As our SLE is constructed to find exactly such a direction, it is guaranteed to be solvable.

\paragraph{Potential bound violation.}
If \autocref{eq:sp_simp:pot_bound} is violated, let $\tilde{w} \in V$, part of component $\Sub' \in \Subs$ be a node with the largest violation (or closest to violation) of its upper bound, \ie,
\begin{equation*}
  \tilde{w} \mathrel{\in} \argmax\{\PotCOpt{u} - \PotMax{u} \mid u \in V\}.
\end{equation*}
We do not need to consider the lower bound case, as we can always shift the potentials such that only upper bounds are violated.
Analogously to the previous discussion, we modify the SLE in \autoeqref{eq:lesd} to
  \begin{align}
    &\sum_{(\SubNode{\Sub}, \SubNode{\OtherSub}) \in \delta^-(\SubNode{\Sub})} \DualPEC{\SubNode{\Sub}}{\SubNode{\OtherSub}} - \sum_{(\SubNode{\OtherSub}, \SubNode{\Sub}) \in \delta^+(\SubNode{\Sub})} \DualPEC{\SubNode{\OtherSub}}{\SubNode{\Sub}} = 0 \enspace \forall \Sub \in \Subs\setminus\{\Sub'\}, \nonumber\\
    &\DualPotBCUb{\SubNode{\Sub'}} + \sum_{(\SubNode{\Sub'},\SubNode{\OtherSub}) \in \delta^-(\SubNode{\Sub'})} \DualPEC{\SubNode{\Sub'}}{\SubNode{\OtherSub}} - \sum_{(\SubNode{\OtherSub},\SubNode{\Sub'}) \in \delta^+(\SubNode{\Sub'})} \DualPEC{\SubNode{\OtherSub}}{\SubNode{\Sub'}} = 0,
  \label{eq:feas_eqs_pot}
  \end{align}
and add the objective function via equation
\begin{align*}
  &\sum_{\Sub \in \Subs} \DemandFixed{\SubNode{\Sub}} \cdot \DualFCC{\SubNode{\Sub}} + \sum_{(\SubNode{\Sub},\SubNode{\OtherSub}) \in \LArcsCC} \left( (\DualOverCLbOpt{\SubNode{\Sub}}{\SubNode{\OtherSub}} + \DualOverCUbOpt{\SubNode{\Sub}}{\SubNode{\OtherSub}}) \cdot \FlowMaxOver{\SubNode{\Sub}}{\SubNode{\OtherSub}} \right)\\
  &+ \PotMax{\SubNode{\Sub'}} \cdot \DualPotBCUb{\SubNode{\Sub'}} = 1.
\end{align*}

As for optimality cuts, solving the respective SLEs yields values for all dual variables, and the desired cut.

\section{Factorization Scheme \& Complexity}\label{sec:sp:comp_fac}
In this section, we discuss the implications of our proposed subproblem solving regarding the overall branch-and-Benders-cut algorithm.
We first discuss how to leverage the similarity of the SLEs within and between subproblems.
Afterward, we analyze the theoretical complexity of the resulting algorithm.

\paragraph{Subproblem similarity.}
During the branch-and-Benders-cut algorithm, the set of active arcs changes between master solutions and between contingencies and their corresponding subproblems.
Each subproblem $\subP$ is solved with a primal and possibly a dual SLE.
Each SLE $Dx = b$ consists of a subproblem-dependent matrix $D \in \R^{n \times n}$ and a right-hand side $b \in \R^n$ representing currently active arcs (cf. \autocsref{eq:sp_simp:flow_cons} and \autoeqref{eq:lesd}), where $n = |\Subs|$ for primal and optimality cut subproblems, and $n = |\Subs| + 1$ in case of a feasibility cut.
We solve these linear equation systems using an $LU$-factorization:
first, $D$ is decomposed into a lower triangular matrix $L$ and an upper triangular matrix $U$, with $L, U \in \R^{|\Subs| \times |\Subs|}$, such that $LU = D$.
Afterward, we solve $Ly = b$ for $y$ and then $Ux = y$ for $x$ to obtain a solution to the original SLE.
Computing the factorization from scratch has as running time complexity of $\mathcal{O}(|\Subs|^3)$ \citep{golub2013matrix}.
This dominates the running time complexity when the matrix is triangular, \eg, for solving an SLE using an existing factorization, $\mathcal{O}(|\Subs|^2)$ \citep{golub2013matrix}.
As discussed below, we can convert the factorization of a matrix into that of a sufficiently similar one in $\mathcal{O}(|\Subs|^2)$ as well.
We thus want to avoid complete factorizations by exploiting the similarity between the matrices of the subproblems that we solve.

The similarities occur between Benders iterations, between contingencies, and between the primal and dual SLEs of the same subproblem.
As observed before, the coefficient matrices of the SLEs for primal $\lesP$ and optimality cuts $\lesD$ are equal.
The similarity of contingencies depends on the contingency size, in our evaluation instances, for example, it is small and independent of the network size.
The switching limit parameter $\TAL$ restricts the number of differing switch arcs between master solutions.

\paragraph{Updating individual arcs.}
Subproblems within or between Benders iterations differ in the underlying contracted graph. 
Switch arc states determine the components of the graph, line arc states determine how components are connected. 
In the simplest case, two subproblems differ in a single line arc that is operational or failed due to a contingency.
A matrix $D^{-(u',v')}$ with arc $(u',v') \in \LArcs$ disconnected differs from the same matrix $D^{+(u',v')}$ with that arc connected only slightly (cf. Equation \eqref{eq:lesp_matrix}):
\begin{equation*}
  D^{-(u',v')} + D^{(u',v')} = D^{+(u',v')},
\end{equation*}
where
\begin{equation*}
  D^{(u',v')}_{(u,v)} = \begin{cases}
                                        -\Cond{u'}{v'},  & \{u, v\} \subsetneq \{u', v'\}\\
                                        \Cond{u'}{v'}, & \{u, v\} =          \{u', v'\}\\
                                        0,                           & \text{else}
                                      \end{cases}
\end{equation*}
for all $u,v \in V$. Let
\begin{equation*}
  q^{(u',v')}_u = \begin{cases}
                                   \Cond{u'}{v'}, & u = u' \\
                                  -\Cond{u'}{v'}, & u = v' \\
                                   0                         , & \text{else,}
                                \end{cases}
\end{equation*}
and
\begin{equation*}
  r^{(u',v')}_u = \begin{cases}
                                  -1, & u = u' \\
                                   1, & u = v' \\
                                   0, & \text{else,}
                                \end{cases}
\end{equation*}
for all $u \in V$.
Then,
\begin{equation*}
  D^{(u',v')} = q^{(u',v')} r^{(u',v')\top}.
\end{equation*}
Given an $LU$-factorization $L^{-(u',v')} U^{-(u',v')} = D^{-(u',v')}$, we can thus perform a rank-$1$ update to arrive at an $LU$-factorization $L^{+(u',v')} U^{+(u',v')} = D^{+(u',v')}$.
This update can be performed in $\mathcal{O}(|V|^2)$ \citep{stange2007efficient}.

\begin{figure*}
  \FIGURE{
  \subcaptionbox{$U$ and $V$ split.}{%
  \centering
  \tikzset{
    aux_node/.style = {circle, draw, minimum size=#1},
    aux_node/.default = 11 
  }
  \tikzset{
    coupler_node/.style = {circle, draw, minimum size=#1},
    coupler_node/.default = 11 
  }
  \begin{tikzpicture}
    \node (0,0)                  (center) {};
    \node [left =1.37 of center] (lcenter) {};
    \node [right=1.37 of center] (rcenter) {};
    \node [coupler_node, above=1.2 of center] (coup1) {};
    \node [coupler_node, below=1.2 of center] (coup2) {};
    \node [aux_node, above =.2 of lcenter] (aux1) {};
    \node [aux_node, above =.2 of rcenter] (aux2) {};
    \node [aux_node, below =.2 of lcenter] (aux3) {};
    \node [aux_node, below =.2 of rcenter] (aux4) {};
    \node [coupler_node, left  =1 of aux1] (naux1) {};
    \node [coupler_node, right =1 of aux2] (naux2) {};
    \node [coupler_node, left  =1 of aux3] (naux3) {};
    \node [coupler_node, right =1 of aux4] (naux4) {};
    \def\coupHalo{0.27cm}
    \path [
      draw=black!50,
      line width=.4pt,
      fill=lightgray,
      fill opacity=.35
    ]
      ($(coup1.center)+(60:\coupHalo)$)
      arc[start angle=60, end angle=120, radius=\coupHalo] --
      ($(aux1.center)+(115:\coupHalo)$)
      arc[start angle=115, end angle=270, radius=\coupHalo] --
      ($(aux2.center)+(270:\coupHalo)$)
      arc[start angle=270, end angle=425, radius=\coupHalo] -- cycle;
    \path [
      draw=black!50,
      line width=.4pt,
      fill=lightgray,
      fill opacity=.35
    ]
    ($(coup2.center)+(240:\coupHalo)$)
    arc[start angle=240, end angle=300, radius=\coupHalo] --
    ($(aux4.center)+(295:\coupHalo)$)
    arc[start angle=295, end angle=450, radius=\coupHalo] --
    ($(aux3.center)+(450:\coupHalo)$)
    arc[start angle=450, end angle=605, radius=\coupHalo] -- cycle;
    \draw [->, dotted] (coup1) -- (coup2) node [midway, anchor=west, right=-.08] {$a_1$};
    \draw [->, dashed] (coup1) -- (aux1)  node [midway, anchor=south east, above left =.19 and -.08] {$U$};
    \draw [->, dashed] (coup1) -- (aux2);
    \draw [->, dotted] (coup1)  -- (aux3);
    \draw [->, dotted] (coup1)  -- (aux4);
    \draw [->, dotted] (coup2)  -- (aux1);
    \draw [->, dotted] (coup2)  -- (aux2);
    \draw [->, dashed] (coup2)  -- (aux3) node [midway, anchor=north east, below left =.19 and -.08] {$W$};
    \draw [->, dashed] (coup2)  -- (aux4);
    \draw [->, dotted] (aux3)   -- (aux4) node [midway, anchor=south, above left= -.1 and .11] {$a_2$};
    \draw [->]         (aux1)   -- (naux1);
    \draw [->]         (aux2)   -- (naux2);
    \draw [->]         (aux3)   -- (naux3);
    \draw [->]         (aux4)   -- (naux4);
    \end{tikzpicture}
  }\hspace{1.5cm}
  \subcaptionbox{$U$ and $V$ merged.}{%
  \centering
  \tikzset{
    aux_node/.style = {circle, draw, minimum size=#1},
    aux_node/.default = 11 
  }
  \tikzset{
    coupler_node/.style = {circle, draw, minimum size=#1},
    coupler_node/.default = 11 
  }
  \begin{tikzpicture}
    \node (center2) {};
    \node [left =1.37 of center2] (lcenter2) {};
    \node [right=1.37 of center2] (rcenter2) {};
    \node [coupler_node, above=1.2 of center2] (coup1_2) {};
    \node [coupler_node, below=1.2 of center2] (coup2_2) {};
    \node [aux_node, above =.2 of lcenter2] (aux1_2) {};
    \node [aux_node, above =.2 of rcenter2] (aux2_2) {};
    \node [aux_node, below =.2 of lcenter2] (aux3_2) {};
    \node [aux_node, below =.2 of rcenter2] (aux4_2) {};
    \node [coupler_node, left  =1 of aux1_2] (naux1) {};
    \node [coupler_node, right =1 of aux2_2] (naux2) {};
    \node [coupler_node, left  =1 of aux3_2] (naux3) {};
    \node [coupler_node, right =1 of aux4_2] (naux4) {};
    \def\coupHalo{0.27cm}
    \path [
      draw=black!50,
      line width=.4pt,
      fill=lightgray,
      fill opacity=.35
    ]
      ($(coup1_2.center)+(60:\coupHalo)$)
      arc[start angle=60, end angle=115, radius=\coupHalo] --
      ($(aux1_2.center)+(115:\coupHalo)$)
      arc[start angle=115, end angle=185, radius=\coupHalo] --
      ($(aux3_2.center)+(185:\coupHalo)$)
      arc[start angle=185, end angle=240, radius=\coupHalo] --
      ($(coup2_2.center)+(-115:\coupHalo)$)
      arc[start angle=-115, end angle=-65, radius=\coupHalo] -- 
      ($(aux4_2.center)+(-65:\coupHalo)$)
      arc[start angle=-65, end angle=5, radius=\coupHalo] --
      ($(aux2_2.center)+(5:\coupHalo)$)
      arc[start angle=5, end angle=65, radius=\coupHalo] -- cycle;
    \draw [->, dashed] (coup1_2) -- (coup2_2) node [midway, anchor=west, right=-.08] {$a_1$};
    \draw [->, dashed] (coup1_2) -- (aux1_2)  node [midway, anchor=south east, above left =.19 and -.08] {$U \cup W$};
    \draw [->, dashed] (coup1_2) -- (aux2_2);
    \draw [->, dotted] (coup1_2) -- (aux3_2);
    \draw [->, dotted] (coup1_2) -- (aux4_2);
    \draw [->, dotted] (coup2_2) -- (aux1_2);
    \draw [->, dotted] (coup2_2) -- (aux2_2);
    \draw [->, dashed] (coup2_2) -- (aux3_2);
    \draw [->, dashed] (coup2_2) -- (aux4_2);
    \draw [->, dashed] (aux3_2)  -- (aux4_2)  node [midway, anchor=south, above left= -.1 and .11] {$a_2$};
    \draw [->]         (aux1_2)  -- (naux1);
    \draw [->]         (aux2_2)  -- (naux2);
    \draw [->]         (aux3_2)  -- (naux3);
    \draw [->]         (aux4_2)  -- (naux4);
  \end{tikzpicture}
  }
  }
  {Component reconfiguration.\label{fig:component_cases}}
  {Reconfiguration line arcs are represented by solid lines, active and inactive switch arcs by dashed and dotted lines, respectively. Activating arc $a_1$ merges or splits components $U$ and $W$ (case 1 and 2). Changing the state of arc $a_2$ has no immediate effect, as its incident nodes are part of the same component regardless (case 3).}
\end{figure*}

Differing switch arc states can lead to a change in the components represented in the matrix.
Recall that the matrix considers only active components, \ie, components that contain at least one node with non-zero demand. 
There are three distinct outcomes regarding switch arc components after changing the state of a single arc:
\begin{itemize}
  \item After an arc is activated, two previously separated components are now connected, resulting in one merged component. If at least one of them is active, we connect it (in the matrix, via rank-$1$ updates) to the line arcs of the other component. If the other component is also active, we disconnect all its line arcs and then delete the corresponding row and column in the matrix. 
  \item The reverse: after deactivating an arc, a component is split into two now distinct components. If at least one of them is active, we disconnect all line arcs of the other component, the matrix entries now represent this first component. If the second component is also active, we add it as a new component to the matrix by adding a row and column to represent the new component in the matrix and then connecting its line arcs.
  \item No change in active components. If the arc is activated, this is because the incident nodes are already part of the same component, or because one of the incident nodes was previously not part of an active component. If the arc is deactivated, this is because the incident nodes are connected via other active switch arcs, or because nodes are removed from the component, but do not form an active component. The matrices stay the same, so we do not need to update the factorization.
\end{itemize}
Refer to \autoref{fig:component_cases} for an example.

A special case occurs if the switch arc is paired with exactly one line arc, \ie, there are nodes $u,v,w \in V$, $(u,v) \in \LArcs$, $(v, w) \in \SArcsOut$, $\delta(v) = \{(u,v), (v,w)\}$.
This can be seen as a combined line and switch arc $(u,w)$, if switched on it has the properties of $(u,v)$ and if switched off, the potentials $u, w$ are disconnected.
The potential at $v$ is not explicitly relevant, it is equal to that of $u$ if the arc is on and to that of $w$ if it is off.
We omit its row and column in the SLE matrices, and treat switching the arc analogously to a contingency including the line arc.

\begin{algorithm*}
\caption{Subproblem Solving with Factorization Updates.}\label{alg:factorization}
\begin{algorithmic}[1]
  \State \textbf{input} network and contingency data $(G,\CondSym,\PotFSym,\aC,H)$; previous factorization $L^{\zFixedSym'}, U^{\zFixedSym'}$; master solution $\zFixedSym$, $\DemandFixedSym$, $\etaFixed$; previous master topology $\zFixedSym'$
  \State $\C' \gets \emptyset$
  \State $L^{\zFixedSym, \C'}, U^{\zFixedSym, \C'}$ $\gets$ rank-$1$ updates to $L^{\zFixedSym'}, U^{\zFixedSym'}$ according to $\zFixedSym$
  \begin{itemize}
    \item[$\bullet$] disconnect all $(u,v) \in \LArcsCSym{\zFixedSym'}{\emptyset} \setminus \LArcsC{\emptyset}$
    \item[$\bullet$] update component representation for all $(u,v) \in \SArcs$, $\zFixed{u}{v} \neq \zFixed{u}{v}'$
    \item[$\bullet$] reconnect all $(u,v) \in \LArcsC{\emptyset}                \setminus \LArcsCSym{\zFixedSym'}{\emptyset}$
  \end{itemize}
  \ForAll{$\C \in \aC$}
    \State $L^{\zFixedSym,\C}, U^{\zFixedSym,\C} \gets$ rank-$1$ updates to $L^{\zFixedSym, \C'}, U^{\zFixedSym, \C'}$
    \begin{itemize}
      \item[\hspace{.65cm}$\bullet$] disconnect all $(u,v) \in \C  \cap \LArcsC{\C'}$
      \item[\hspace{.65cm}$\bullet$] update component representation for all $(u,v) \in ((\C \cup \C') \setminus (\C \cap \C')) \setminus \SArcs^{\bar{\zSym}}$
      \item[\hspace{.65cm}$\bullet$] reconnect  all $(u,v) \in \C' \cap \LArcsCC $
    \end{itemize}
    \State solve primal SLE with $L^{\zFixedSym,\C}, U^{\zFixedSym,\C}$ (cf. Subsection \ref{sec:sp:primal})
    \If{$\subP$ infeasible}
      \State feasibility cut: expand $L^{\zFixedSym,\C}, U^{\zFixedSym,\C}$ and solve feasibility cut SLE (cf. Subsection \ref{sec:sp:feas_cut})
      \begin{itemize}
        \item[\hspace{1.25cm}$\bullet$] for $\DualFCCSym, \DualFlowBCUb{\tilde{u}}{\tilde{v}}$, where $(\tilde{u}, \tilde{v}) \in \LArcsOut$ is a violating arc, or
        \item[\hspace{1.25cm}$\bullet$] for $\DualFCCSym, \DualPotBCUb{\tilde{u}}$, where $\tilde{u} \in V$ is a violating node
      \end{itemize}
      \Else
      \State $\eta_\C \gets \subP^*$
      \If{$\eta_\C > \etaFixed_\C$}
        \State optimality cut: solve optimality cut SLE with $L^{\zFixedSym,\C}, U^{\zFixedSym,\C}$ (cf. Subsection \ref{sec:sp:opt_cut})
      \EndIf
    \EndIf
    \State$\C' \gets \C$
  \EndFor
  \State $L^{\zFixedSym, \C'}, U^{\zFixedSym, \C'} \gets$ rank-$1$ updates to $L^{\zFixedSym, \C'}, U^{\zFixedSym, \C'}$
  \begin{itemize}
    \item[$\bullet$] update component representation for all $(u,v) \in \C \cap \SArcsz$
    \item[$\bullet$] reconnect all line arcs $(u,v) \in \C'$
  \end{itemize}
  \State \Return all found cuts
\end{algorithmic}
\end{algorithm*}

\paragraph{Factorization update algorithm.}
Multiple differing arcs can be resolved by an iterative chain of such updates.
Algorithm \autoref{alg:factorization} is called after a valid integer feasible solution is found during the branch-and-Benders-cut algorithm.
It describes how to perform one subproblem iteration of the Benders algorithm, \ie, solve a subproblem for each contingency and return optimality or feasibility cuts if applicable.
We depart from the factorization of the previous master topology $\zFixedSym' \in \{0,1\}^{|\SArcs|}$, $L^{\zFixedSym'} U^{\zFixedSym'} = D^{\zFixedSym'}$. For the first call, this is the ``base'' factorization $LU = D$ that corresponds to the equation system with all arcs switched according to the default configuration of the network, \ie, $\zFixed{u}{v}' = i$ for all $(u,v) \in \SArcs^i$, $i\in\{0,1\}$. 

From this initial factorization, we want to compute the factorization $L^{\zFixedSym} U^{\zFixedSym} = D^{\zFixedSym}$ of the current master topology $\zFixedSym \in \{0,1\}^{\SArcs}$. For each arc $(u,v) \in \SArcs$, $\zFixed{u}{v} \neq \zFixed{u}{v}'$, that differs between the master solutions, the underlying contracted graph can change.
As discussed before, this can lead to a change in the components represented in the matrix, which is resolved by sequentially activating or deactivating the line arcs incident to either component.

After updating the factorization to the master topology, a subproblem is solved for each contingency.
To go from one contingency $\C' \in \aC$ to the next $\C \in \aC$, we need to disconnect each line arc in $\C \cap \LArcsC{\C'}$, and connect each line arc in $\C' \cap \LArcsCC$.
Each of these requires a rank-$1$ update, if the incident components are active.
Switch arcs $(u,v) \in ((\C \cup \C') \setminus (\C \cap \C')) \setminus \SArcs^{\bar{\zSym}}$ activated or deactivated between the contingencies can lead to a change in matrix components, handled as above.

We then solve the primal SLE, as discussed in Subsection \ref{sec:sp:primal}.
If the subproblem is feasible and its objective function value $\eta_\C$ is smaller than or equal to the ``estimated'' objective function value $\etaFixed$, no cut needs to be added, we proceed with the next subproblem.

If it is feasible and its objective function value is larger than $\etaFixed$, we compute an optimality cut.
As detailed in Subsection \ref{sec:sp:opt_cut}, the matrix of the corresponding SLE is the same as for the primal subproblem, only the right hand side changes.

If it is infeasible, let $(\tilde{u}, \tilde{v}) \in \LArcsCC$ be an arc violating the hard flow bounds, if one exists. As discussed in Subsection \ref{sec:sp:feas_cut}, we need to solve the SLE
\begin{equation*}
  \begin{pmatrix}
    D^{\zFixedSym, \C}             & s^{(\tilde{u},\tilde{v})} \\
    \DemandFixedSym^\top & \FlowMax{\tilde{u}}{\tilde{v}} 
  \end{pmatrix}
  \begin{pmatrix}
    \DualFCCSym \\
    \DualFlowBCUb{\tilde{u}}{\tilde{v}}
  \end{pmatrix}
  =
  \begin{pmatrix}
    0\\
    1
  \end{pmatrix},
\end{equation*}
where $D^{\zFixedSym, \C}$ is the matrix of the primal SLE, $\DemandFixedSym$ is the dual objective function coefficient vector for $\DualFCCSym$, and
\begin{equation*}
  s^{(\tilde{u},\tilde{v})}_{\SubNode{\Sub}} = \begin{cases}
                                   1, & \SubNode{\Sub} = \tilde{u} \\
                                  -1, & \SubNode{\Sub} = \tilde{v} \\
                                   0, & \text{else},
                                \end{cases}               \quad \forall \Sub \in \Subs
\end{equation*}
corresponds to the coefficients of $\DualFlowBCUb{\tilde{u}}{\tilde{v}}$ in Equations \eqref{eq:feas_eqs_flow}.
We already have the factorization for $D^{\zFixedSym, \C} = L^{\zFixedSym,\C} U^{\zFixedSym,\C}$ and want to expand it to cover the remaining entries of the matrix.
It must hold that
\begin{equation*}
  \begin{aligned}
     &
    \begin{pmatrix}
      D^{\zFixedSym, \C}             & s^{(\tilde{u},\tilde{v})} \\
      \DemandFixedSym^\top & \FlowMax{\tilde{u}}{\tilde{v}} 
    \end{pmatrix}\\
    =&
    \begin{pmatrix}
      L^{\zFixedSym,\C}            & 0 \\
      \ell^{\zFixedSym,\C\top} & 1 
    \end{pmatrix}
    \begin{pmatrix}
      U^{\zFixedSym,\C} & u^{\zFixedSym,\C} \\
      0^\top & v^{\zFixedSym,\C}
    \end{pmatrix}\\
    =&
    \begin{pmatrix}
      L^{\zFixedSym,\C} U^{\zFixedSym,\C} & L^{\zFixedSym,\C} u^{\zFixedSym,\C}\\
      \ell^{\zFixedSym,\C\top} U^{\zFixedSym,\C}             & \ell^{\zFixedSym,\C\top} u^{\zFixedSym,\C} + v^{\zFixedSym,\C}
    \end{pmatrix}.
  \end{aligned}
\end{equation*}
We need to compute $\ell^{\zFixedSym,\C}, u^{\zFixedSym,\C} \in \R^{\Subs}$ and $v^{\zFixedSym,\C} \in \R$.
For this, we can solve the SLEs
\begin{equation*}
   U^{\zFixedSym,\C\top} \ell^{\zFixedSym,\C} = \DemandFixedSym
\end{equation*}
and
\begin{equation*}
   L^{\zFixedSym,\C} u^{\zFixedSym,\C} = s^{(\tilde{u},\tilde{v})},
\end{equation*}
and then solve
\begin{equation*}
  v^{\zFixedSym,\C} = \FlowMax{\tilde{u}}{\tilde{v}} - \ell^{\zFixedSym,\C\top} u^{\zFixedSym,\C}.
\end{equation*}
Note that $U^{\zFixedSym,\C\top}$ and $L^{\zFixedSym,\C}$ are already triangular, and these SLEs are thus solvable in $\mathcal{O}(|\Subs|^2)$.

If no violating arc exists, let instead $\tilde{u} \in V$ be a node with the highest or closest to violation, as discussed in Subsection \ref{sec:sp:feas_cut}.
The SLE we need to solve is
\begin{equation*}
  \begin{pmatrix}
    D^{\zFixedSym, \C}             & s^{\tilde{u}} \\
    \DemandFixedSym^\top & \PotMax{\tilde{u}} 
  \end{pmatrix}
  \begin{pmatrix}
    \DualFCCSym \\
    \DualPotBCUb{\tilde{u}}
  \end{pmatrix}
  =
  \begin{pmatrix}
    0\\
    1
  \end{pmatrix},
\end{equation*}
where
\begin{equation*}
  s^{\tilde{u}}_{\SubNode{\Sub}} = \begin{cases}
                                      1, & \SubNode{\Sub} = \tilde{u} \\
                                      0, & \text{else},
                                  \end{cases}               \quad \forall \Sub \in \Subs.
\end{equation*}
We proceed analogously to update the factorization.

Finally, after all subproblems have been solved, we revert the changes to $L^{\zFixedSym,\C}$ and $U^{\zFixedSym,\C}$ caused by the contingencies, and return all cuts that were found.

\paragraph{Theoretical complexity.}
As a last point, we discuss the theoretical complexity of the proposed methods, as a complement to the empirical results in the next section.
We propose a method with which the subproblems can be solved with systems of linear equation, and combine it with a factorization update scheme.
We compare this to solving the subproblems as LPs with a simplex algorithm.
The methods differ only in their way to obtain primal and dual solutions.

While a simplex algorithm can require up to $2^{n-1}$ iterations, in practice it empirically uses only $\mathcal{O}(m)$, with a small constant \citep{BertsimasTsitsiklis1997}, where $n = |V| + |\SArcsOut|$ is the number of variables and $m = 4|\SArcsOut| + 3|V| + 4|\LArcsOut|$ is the number of constraints.
Each of these iterations requires solving up to three SLEs of size $m$.
For solving the SLEs, a factorization update scheme similar to ours can be used, which reduces their overall complexity to $\mathcal{O}(m^2)$ per simplex iteration.
We solve $|\aC|$ different subproblems in one Benders iteration.
In total, we arrive at $\mathcal{O}(|\aC| \cdot m \cdot m^2) = \mathcal{O}(|\aC| \cdot (|V| + |\Arcs|)^3) = \mathcal{O}(|\aC| \cdot |\Arcs|^3)$.
When the number of contingencies is restrained to the number of line arcs by some constant factor (\eg, one contingency per line arc), this reduces to $\mathcal{O}(|\Arcs|^4)$.

To solve each subproblem with our method, we need to solve at most two SLEs.
Additionally, these SLEs contain only at most one variable per node.
Solving all of them without factorization updates thus has a total complexity of $\mathcal{O}(|\aC| \cdot 2 \cdot |V|^3) \subseteq \mathcal{O}(|\aC| \cdot |\Arcs|^3)$, or $\mathcal{O}(|\Arcs|^4)$ if the number of contingencies is restrained by the number of arcs.
This can be reduced when using factorization updates:
Let $\hat{\C} \in \Z$ be the size of the largest contingency.
Computing the initial factorization, or updating from a previous master solution, can be done in $\mathcal{O}(|V|^3)$.
For each subproblem, we then need to revert up to $\hat{\C}$ updates from the previous subproblem, and apply equally many of the current subproblem.
Afterward, we solve the primal and dual SLE via forward and backward substitution.
For feasibility cuts, we require another forward and backward substitution to update the factorization.
We arrive at a total complexity of $\mathcal{O}(|V|^3 + |\aC| \cdot (2 \cdot \hat{\C} \cdot |V|^2 + 4 \cdot |V|^2)) = \mathcal{O}(|V|^3 + |\aC| \cdot \hat{\C} \cdot |V|^2) \subseteq \mathcal{O}(|\Arcs|^3 + |\aC| \cdot \hat{\C} \cdot |\Arcs|^2)$. If the number of contingencies is restrained by the number of arcs, and $\hat{\C}$ is independent of the network size, this reduces to $\mathcal{O}(|\Arcs|^3) \subset \mathcal{O}(|\Arcs|^4)$.
This is, \eg, the case for the so-called $N-1$ criterion, where we have at most one contingency per line arc, and each contingency contains at most one line arc.
Note that our coarse analysis is overly conservative for our approach as in sparse networks where the number of arcs is linear in the number of nodes, whereas in general $|\Arcs| \in \mathcal{O}(|V|^2)$.

It should be remarked that our method can be seen as a specialized simplex algorithm, and many standard acceleration techniques apply.
These are out of the scope of this work.
We refer to, for example, \citet{vanderbei1998linear} and \citet{BertsimasTsitsiklis1997} for details.

%% file: chapters/5-evaluation.tex
\newcommand{\compact}{\texttt{compact}}
\newcommand{\fulllp}{\texttt{GRB}}
\newcommand{\sle}{\texttt{SLE}}
\newcommand{\sleu}{{\texttt{SLE-U}}{}}

\newcolumntype{Y}{>{\arraybackslash}X} 
\newcolumntype{R}{S[table-format=4.2, round-mode=places, round-precision=2, table-space-text-pre={(}, table-space-text-post={)}, detect-weight=true]} 
\newcolumntype{Z}{>{\centering\arraybackslash}S[table-format=3.0, round-mode=places, round-precision=0, table-space-text-pre={(}, table-space-text-post={)}]} 

In this section, we evaluate the performance of the proposed subproblem algorithms: ``classical'' with Gurobi \fulllp{}, with the systems of linear equations solver \sle{}, introduced in Sections \ref{sec:sp:primal}--\ref{sec:sp:feas_cut}, as well as the systems of linear equations solver with factorization updates \sleu{} (Section \ref{sec:sp:comp_fac}).

All algorithms are implemented using Julia 1.12.6 \citep{Julia-2017}.
We use PowerModels 0.21.85 \citep{PowerModels} with JuMP 1.30.0 \citep{Lubin2023} to model the Benders master, as well as the Benders subproblem for \fulllp{}.
They are solved using Gurobi 13.0.0 \citep{gurobi}, with one model object modified between subproblems.
The Benders functionality is implemented via Gurobi's callbacks.
For \sle{} and \sleu{}, subproblems are solved with the Julia package LinearSolve 3.75.0 \citep{LinearSolve}.
The factorization update functionality for \sleu{} is newly implemented, based on the work by \citet{stange2007efficient}.
To avoid numerical issues, factorizations are recomputed in regular intervals, among other things whenever the dimension of the coefficient matrix changes.

\paragraph{Evaluation setting.}
We evaluate the methods on instances of the DC security-constrained optimal transmission switching problem (SCOTS).
For these instances, the graph represents a DC power network, with potentials and flows representing \emph{voltage angles} and \emph{power flows}, respectively.
The potential function is affine, with $r_{(u,v)} = 1$ on all arcs $(u,v) \in \LArcs$.
The constant $s_{(u,v)} \in \R$ can model, \eg, transformers.
Most line arcs in the instances are switchable themselves, which can be modelled as the aforementioned pairing of two correspondingly parametrized line arcs and one switch arc in the potential-based flow setting.
Switch arcs can also represent \emph{substations}, \ie, components that are located so close to each other that their electrical properties can be disregarded.
This corresponds to the switch arc components discussed in the previous sections.
The set of contingencies aligns with the $N-1$-criterion, \ie, one contingency per line arc, as long as it does not split the network.
As discussed, any affine potential based flow instance can be transformed into one with an identity function as potential equation, and thus also a SCOTS instance.

\paragraph{Benchmark dataset.}
We use a custom benchmark set based on the PGLib-OPF \citep{pglib} set, instances of which are commonly used for testing power flow applications \citep{fisher2008optimal,hedman_optimal_2009,heidarifar2015network}.
The original benchmark set itself is insufficient since it contains too few instances that are solvable within sufficient time with our computational setup.
Additionally, many of the instances are infeasible for the SCOTS.
We compute new instances by joining parts of the existing instances onto each other.
The resulting benchmark set contains 174 instances.
The construction of the benchmark set will be detailed in an upcoming publication.

\paragraph{Computational Setup.}
We use machines with $2 \times$ Intel Xeon L5630 quad-core processors at 2.13\,GHz (8~cores, 16\,GB RAM) on Debian GNU/Linux 11 across all experiments.
A three-hour time limit is imposed per instance.

\subsection{Results}
We compare the three methods under scenarios with a reduced number of available switch arcs.
For the remaining discussion in this section, we disregard all 25 instances that are infeasible or where all instances encountered the time limit, leaving a remaining 149 instances.
Where primal solutions were found, the reported gaps for instances exceeding the time limit remain within \qty{1}{\percent} for almost all instances, and differ only marginally between solution methods.
We hence omit gap values from further discussion.

\paragraph{Time distribution.}
Since we exclusively change the method for solving subproblems while leaving the output unchanged, the largest difference between the methods can be expected to occur in the subproblem running time.
\begin{figure*}
    \FIGURE
    {
    \subcaptionbox{Relative.}{%
    \input{data/final/time_decomp_boxplot/detailed_relative.tikz}
    }\hspace{.75cm}
    \subcaptionbox{Absolute.}{%
    \input{data/final/time_decomp_boxplot/detailed.tikz}
    }
    }
    {Mean running time distribution.\label{fig:s:rt_dist}}
    {
    Legend in order from left to right for each method. ``Primal'' is the time spent exclusively solving the subproblem LP with Gurobi or the primal SLE, respectively. ``Dual'' is the equivalent for the dual SLE for \sle{} and \sleu{}. ``Update'' is the factorization update time for \sleu{}. ``Misc. SP'' is any remaining SP time, including connectedness separation and generation of the cut object. ``MP'' is the MP solving time. Each box spans the interquartile range, whiskers extend to 1.5$\times$ IQR, dots are outliers.}
\end{figure*}
\autoref{fig:s:rt_dist} shows the relative distribution of the different parts of the solving process for all three methods in box plots.
We observe that for \fulllp{}, the mean subproblem running time dominates the overall running time at \qty{72.4905}{\percent}, with Gurobi's solving process taking up \qty{49.4566}{\percent} (``Primal'', dual solutions are obtained without relevant additional effort).
As such, noticeable improvements in that solving process also result in noticeable improvements in the overall solving process.
While the absolute time spent on solving the MP remains relatively constant for all methods, its relative running time contribution increases to \qty{75.2525}{\percent} for \sle{} and \qty{82.9347}{\percent} for \sleu{}.
This is due to reductions in the time spent on solving the ``core'' of the SP, \ie, in the LP or SLE solver:
for \sle{}, solving the primal SLEs uses \qty{9.3526}{\percent}, and solving the optimality and feasibility SLEs \qty{4.7457}{\percent} of running time.
For \sleu{}, these are further reduced to \qty{0.6962}{\percent} and \qty{0.6869}{\percent}, respectively, with an additional \qty{0.7662}{\percent} for factorization updates.
The remaining time is taken up by miscellaneous subproblem computations, including cut construction and allocation, checking for disconnected components, updating the LP and SLE objects and checking the solutions against bounds.
Despite effort spent to reduce the time spent on these computations, it remains a noticeable contributor.
Further improvements here are likely possible, but would not influence the SP solving time.

MP and miscellaneous SP times can be seen as relatively static, as they are not directly influenced by the choice of SP solver.
This can be observed particularly for \sle{} and \sleu{}, where they match very closely ($\qty{1099.7319}{\second}$ vs.\ $\qty{1142.3121}{\second}$).
\fulllp{} times differ slightly, this is due to a slightly different set of miscellaneous SP computations, as well as possible differences in dual solutions.
While we often observed equal dual solutions, up to numerical differences, between Gurobi and our algorithms, the subproblems do not necessarily have unique solutions.
In fact, as discussed, while our dual solutions are always feasible and where applicable optimal, a different choice of connected component spanning trees or violated flow or potential for feasibility cuts can yield different solutions for our algorithms as well.
These differences seem to be mildly disadvantageous to our methods (95747.2581 solved SPs for \fulllp{}, 98962.1218 for \sle{}, 101760.1538 for \sleu{} on average).
While this may point to slightly weaker cuts that we produce, it can also be attributed to advanced internal solver mechanisms that are not publicly documented.

\paragraph{Running time comparison.}
\begin{table}
    \TABLE{Aggregate running time statistics for \fulllp{}, \sle{} and \sleu{}.\label{tab:s:aggregate}}
    {
    \centering
    \begin{tabular}{lrrr}
    \toprule
    \up\down& Solved & Mean & Geo. Mean \\
    \midrule
    \up\texttt{GRB} & 141 & 2321.6646 & 958.7539 \\
    \texttt{SLE} & 146 & 986.1129 & 367.8782 \\
    \down\texttt{SLE-U} & 147 & 906.2038 & 276.0177 \\
    \bottomrule
    \end{tabular}
    }
    {``Mean'' is the arithmetic mean, "geo.\ mean" is the shifted geometric mean with shift $10$. Aggregate values cover only the $136$ instances solved by all methods.}
\end{table}
\autoref{tab:s:aggregate} provides an aggregate view on the running time behavior of the methods. In every metric, there is a clear hierarchy, with \sleu{} as the fastest solver, \sle{} in second, and \fulllp{} last.

To investigate the methods as a part of the overall solving process as well as individually, the following time comparisons feature the total solve time (\autoref{fig:s:perf_prof_total}), the total SP time (\autoref{fig:s:perf_prof_SP}), which includes the miscellaneous SP time, and the SP solve time (``Primal,'' ``Dual,'' and ``Update'' in \autoref{fig:s:perf_prof_solve}).
\begin{figure}
    \FIGURE{
    \input{data/final/performance_profiles/total/factor.tikz}
    }
    {Performance profile for total solve time.\label{fig:s:perf_prof_total}}
    {}
\end{figure}
The figures depict performance profiles, where the $x$-axis shows the relative slowdown compared to the best performing method on a given instance, and the $y$-axis the fraction of instances where the slowdown is at most as large as the $x$-value.
For example, in \autoref{fig:s:perf_prof_total} for $\qty{89.9329}{\percent}$ of instances, \fulllp{} is at most $10$ times slower than the best performing method.

For total running time, \autoref{fig:s:perf_prof_total}, \sleu{} is the fastest solver for $\qty{83.2215}{\percent}$ of instances, and less than $\qty{54.7700}{\percent}$ slower for $\qty{95}{\percent}$.
\sle{} is the fastest method for $\qty{14.7651}{\percent}$ of instances.
It is at least $2$ times slower for about $\qty{12.0805}{\percent}$ of instances, and at least $3$ times slower for about $\qty{4.6980}{\percent}$ of instances.
\fulllp{} is by far the slowest method. Only \qty{10.0671}{\percent} of instances are at most $2$ times slower. Half of all instances are at least $3.4924$ times slower, $\qty{10.0671}{\percent}$ of instances are at least $10$ times slower and $\qty{5.3691}{\percent}$ of instances are at least $100$ times slower.

\begin{figure}
    \FIGURE{
    \subcaptionbox{Total time.\label{fig:s:perf_prof_SP}}{%
    \input{data/final/performance_profiles/SP/factor.tikz}
    }\hfill
    \subcaptionbox{Solve time (GRB/SLEs).\label{fig:s:perf_prof_solve}}{%
    \input{data/final/performance_profiles/solve_time/factor.tikz}
    }
    }
    {Performance profiles for SP time.}
    {}
\end{figure}
The magnitude of the differences increases when observing only the SP time.
Here, \sle{} is at least $2$ times slower for $\qty{64.4295}{\percent}$ of instances, and at least $3$ times slower for about $\qty{10.0671}{\percent}$ of instances.
For \fulllp{}, $\qty{69.7987}{\percent}$ of instances are at least $10$ times slower. $\qty{5.3691}{\percent}$ of instances remain at least $100$ times slower.

The SP solve time, the focus of the methods, shows the largest improvements.
\sleu{} is the fastest method for all instances.
\sle{} is at least $2$ times slower for $\qty{97.3154}{\percent}$ of instances, and at least $10$ times slower for about $\qty{47.6510}{\percent}$ of instances.
For \fulllp{}, $\qty{86.5772}{\percent}$ of instances are at least $30$ times slower, and $\qty{16.1074}{\percent}$ of instances are at least $100$ times slower.

In Section \ref{sec:sp:comp_fac}, we discussed the theoretical complexity of the methods, and observed that, for single line arc contingencies, \sle{} has a running time complexity in $\mathcal{O}(|A|\cdot|V|^3)$ and \sleu{} in $\mathcal{O}(|V|^3 + |A| \cdot |V|^2)$.
While $|A|$ would be equal to $|V|^2$ for a fully connected graph, the instances in the benchmark set have an average of only $1.4683$ arcs per vertex, \ie, the underlying graphs are sparsely connected, and thus a relatively low number of subproblems are solved per iteration.
Additionally, the SLE matrices corresponding to a graph with all switch arcs active and without contingencies have an average number of non-zero entries of $\qty{1.80}{\percent}$ over all instances.
SLE solution algorithms specialized on sparse matrices exhibit near-linear running time based on the number of non-zero elements of the matrix \citep{schaub2017sparse,spielman2004nearly}. 
\begin{figure*}
    \FIGURE{
    \subcaptionbox{Performance profile.\label{fig:s:perf_prof_solve_scaled}}{%
    \input{data/final/performance_profiles/solve_time_scaled/factor.tikz}
    }\hfill
    \subcaptionbox{Scatter plot.\label{fig:s:scatter_solve}}{%
    \input{data/final/performance_profiles/solve_time_scaled/scatter_sorted_by_sle_runtime.tikz}
    }
    }
    {Synthesized data for SP solve times (SLEs).}
    {Data for artificial methods $\sle{}^{2/3}$ and $\sle{}^{1/2}$ is synthesized from individual $\sle{}$ data points.}
\end{figure*}

Figures \ref{fig:s:perf_prof_solve_scaled} and \ref{fig:s:scatter_solve} visualize this relationship.
Here, we added artificial methods $\sle{}^{2/3}$ and $\sle{}^{1/2}$, and compare them with \sle{} and \sleu{}.
For each instance, letting $t \in \R$ be the time in seconds needed by $\sle{}$, the time needed for $\sle{}^{2/3}$ and $\sle{}^{1/2}$ is set to be $t^{2/3}$ and $t^{1/2}$, respectively.
This roughly emulates a hypothetical method in a complexity class correspondingly reduced compared to $\sle{}$.
Note that this view ignores hidden constants in $\mathcal{O}$-notation and only serves as an intuition.
In \autoref{fig:s:perf_prof_solve_scaled}, \sleu{} matches closely with $\sle{}^{1/2}$, suggesting that it indeed exhibits a near-linear running time and \sle{} a quadratic one.
Note, however, that the $y$-axis of the performance denotes the fraction of instances instead of specific instances. 
An instance-by-instance comparison is provided by \autoref{fig:s:scatter_solve}.
\sleu{} seems to be closer to $\sle{}^{1/2}$ for easier instances, and slows down to between $\sle{}^{1/2}$ and $\sle{}^{2/3}$ for harder instances.

%% file: data/final/time_decomp_boxplot/detailed_relative.tikz
\begin{tikzpicture}
\providecommand{\mathdefault}[1]{#1}

\definecolor{chocolate213940}{RGB}{213,94,0}
\definecolor{cornflowerblue86180233}{RGB}{86,180,233}
\definecolor{darkcyan0114178}{RGB}{0,114,178}
\definecolor{darkcyan0158115}{RGB}{0,158,115}
\definecolor{darkgray176}{RGB}{176,176,176}
\definecolor{gray}{RGB}{128,128,128}
\definecolor{orange2301590}{RGB}{230,159,0}
\definecolor{whitesmoke242}{RGB}{242,242,242}

\begin{axis}[
height=7.14000cm,
legend cell align={left},
legend style={
  font=\footnotesize,
  fill=none,
  fill opacity=1,
  draw opacity=1,
  text opacity=1,
  at={(0.330645161290323,0.995)},
  anchor=north,
  draw=none
},
tick align=outside,
tick pos=left,
width=7.00000cm,
x grid style={darkgray176},
xmin=-2, xmax=10.4,
xtick style={color=black},
xtick={0,4.2,8.4},
xticklabels={\texttt{GRB},\texttt{SLE},\texttt{SLE-U}},
y grid style={darkgray176},
ylabel={Share of solve time (\%)},
ymin=-3, ymax=103,
ytick style={color=black},
ytick={-20,0,20,40,60,80,100,120},
yticklabels={
  $\mathdefault{\ensuremath{-}20}$,
  $\mathdefault{0}$,
  $\mathdefault{20}$,
  $\mathdefault{40}$,
  $\mathdefault{60}$,
  $\mathdefault{80}$,
  $\mathdefault{100}$,
  $\mathdefault{120}$
}
]

\path [draw=none, fill=whitesmoke242] (axis cs:2.1,-2.748608776844071) rectangle (axis cs:6.300000000000001,102.74860877684407);
\path [draw=black, fill=darkcyan0114178, line width=0.32pt]
(axis cs:-1.775,46.0334050821289)
--(axis cs:-1.225,46.0334050821289)
--(axis cs:-1.225,56.7583553747621)
--(axis cs:-1.775,56.7583553747621)
--(axis cs:-1.775,46.0334050821289)
--cycle;
\addplot [line width=0.32pt, black, forget plot]
table {%
-1.5 46.0334050821289
-1.5 33.8856223327968
};
\addplot [line width=0.32pt, black, forget plot]
table {%
-1.5 56.7583553747621
-1.5 68.7038831207593
};
\addplot [line width=0.32pt, black, forget plot]
table {%
-1.6375 33.8856223327968
-1.3625 33.8856223327968
};
\addplot [line width=0.32pt, black, forget plot]
table {%
-1.6375 68.7038831207593
-1.3625 68.7038831207593
};
\addplot [black, mark=*, mark size=1.5, mark options={solid,fill=gray}, only marks, forget plot]
table {%
-1.5 26.5864657666445
-1.5 81.0318632437103
-1.5 85.9897028492997
-1.5 85.2234041814854
-1.5 77.2118814140415
};
\path [draw=black, fill=orange2301590, line width=0.32pt]
(axis cs:-1.025,0.00423231805110793)
--(axis cs:-0.475,0.00423231805110793)
--(axis cs:-0.475,0.0086298012565366)
--(axis cs:-1.025,0.0086298012565366)
--(axis cs:-1.025,0.00423231805110793)
--cycle;
\addplot [line width=0.32pt, black, forget plot]
table {%
-0.75 0.00423231805110793
-0.75 0.000225964538350421
};
\addplot [line width=0.32pt, black, forget plot]
table {%
-0.75 0.0086298012565366
-0.75 0.0151596635134847
};
\addplot [line width=0.32pt, black, forget plot]
table {%
-0.8875 0.000225964538350421
-0.6125 0.000225964538350421
};
\addplot [line width=0.32pt, black, forget plot]
table {%
-0.8875 0.0151596635134847
-0.6125 0.0151596635134847
};
\addplot [black, mark=*, mark size=1.5, mark options={solid,fill=gray}, only marks, forget plot]
table {%
-0.75 0.0182196339262305
-0.75 0.0229357694010152
-0.75 0.0243488286781141
-0.75 0.0165635004661849
-0.75 0.0205612611561513
-0.75 0.0267942140481358
-0.75 0.0336276506045831
};
\path [draw=black, fill=darkcyan0158115, line width=0.32pt]
(axis cs:-0.275,0)
--(axis cs:0.275,0)
--(axis cs:0.275,0)
--(axis cs:-0.275,0)
--(axis cs:-0.275,0)
--cycle;
\addplot [line width=0.32pt, black, forget plot]
table {%
0 0
0 0
};
\addplot [line width=0.32pt, black, forget plot]
table {%
0 0
0 0
};
\addplot [line width=0.32pt, black, forget plot]
table {%
-0.1375 0
0.1375 0
};
\addplot [line width=0.32pt, black, forget plot]
table {%
-0.1375 0
0.1375 0
};
\path [draw=black, fill=chocolate213940, line width=0.32pt]
(axis cs:0.475,23.1255033084755)
--(axis cs:1.025,23.1255033084755)
--(axis cs:1.025,33.0885376700384)
--(axis cs:0.475,33.0885376700384)
--(axis cs:0.475,23.1255033084755)
--cycle;
\addplot [line width=0.32pt, black, forget plot]
table {%
0.75 23.1255033084755
0.75 12.1829380923432
};
\addplot [line width=0.32pt, black, forget plot]
table {%
0.75 33.0885376700384
0.75 47.5125613791914
};
\addplot [line width=0.32pt, black, forget plot]
table {%
0.6125 12.1829380923432
0.8875 12.1829380923432
};
\addplot [line width=0.32pt, black, forget plot]
table {%
0.6125 47.5125613791914
0.8875 47.5125613791914
};
\addplot [black, mark=*, mark size=1.5, mark options={solid,fill=gray}, only marks, forget plot]
table {%
0.75 49.0753727753202
0.75 57.0903942412157
0.75 48.6702316265611
0.75 48.6350241854666
0.75 53.6044322285554
};
\path [draw=black, fill=cornflowerblue86180233, line width=0.32pt]
(axis cs:1.225,10.7384022396446)
--(axis cs:1.775,10.7384022396446)
--(axis cs:1.775,26.1133845733636)
--(axis cs:1.225,26.1133845733636)
--(axis cs:1.225,10.7384022396446)
--cycle;
\addplot [line width=0.32pt, black, forget plot]
table {%
1.5 10.7384022396446
1.5 0.389685232292582
};
\addplot [line width=0.32pt, black, forget plot]
table {%
1.5 26.1133845733636
1.5 44.8563053197102
};
\addplot [line width=0.32pt, black, forget plot]
table {%
1.3625 0.389685232292582
1.6375 0.389685232292582
};
\addplot [line width=0.32pt, black, forget plot]
table {%
1.3625 44.8563053197102
1.6375 44.8563053197102
};
\path [draw=black, fill=darkcyan0114178, line width=0.32pt]
(axis cs:2.425,8.85067709451099)
--(axis cs:2.975,8.85067709451099)
--(axis cs:2.975,19.8976150851787)
--(axis cs:2.425,19.8976150851787)
--(axis cs:2.425,8.85067709451099)
--cycle;
\addplot [line width=0.32pt, black, forget plot]
table {%
2.7 8.85067709451099
2.7 2.00237680520828
};
\addplot [line width=0.32pt, black, forget plot]
table {%
2.7 19.8976150851787
2.7 31.356945283194
};
\addplot [line width=0.32pt, black, forget plot]
table {%
2.5625 2.00237680520828
2.8375 2.00237680520828
};
\addplot [line width=0.32pt, black, forget plot]
table {%
2.5625 31.356945283194
2.8375 31.356945283194
};
\addplot [black, mark=*, mark size=1.5, mark options={solid,fill=gray}, only marks, forget plot]
table {%
2.7 47.2586948315789
2.7 48.0690528442425
};
\path [draw=black, fill=orange2301590, line width=0.32pt]
(axis cs:3.175,4.55774710623796)
--(axis cs:3.725,4.55774710623796)
--(axis cs:3.725,8.16581346947872)
--(axis cs:3.175,8.16581346947872)
--(axis cs:3.175,4.55774710623796)
--cycle;
\addplot [line width=0.32pt, black, forget plot]
table {%
3.45 4.55774710623796
3.45 0.846175220596752
};
\addplot [line width=0.32pt, black, forget plot]
table {%
3.45 8.16581346947872
3.45 12.2127214200311
};
\addplot [line width=0.32pt, black, forget plot]
table {%
3.3125 0.846175220596752
3.5875 0.846175220596752
};
\addplot [line width=0.32pt, black, forget plot]
table {%
3.3125 12.2127214200311
3.5875 12.2127214200311
};
\addplot [black, mark=*, mark size=1.5, mark options={solid,fill=gray}, only marks, forget plot]
table {%
3.45 15.0235103092453
3.45 13.6127834512119
3.45 17.0250957101731
3.45 17.2078402732797
};
\path [draw=black, fill=darkcyan0158115, line width=0.32pt]
(axis cs:3.925,0)
--(axis cs:4.475,0)
--(axis cs:4.475,0)
--(axis cs:3.925,0)
--(axis cs:3.925,0)
--cycle;
\addplot [line width=0.32pt, black, forget plot]
table {%
4.2 0
4.2 0
};
\addplot [line width=0.32pt, black, forget plot]
table {%
4.2 0
4.2 0
};
\addplot [line width=0.32pt, black, forget plot]
table {%
4.0625 0
4.3375 0
};
\addplot [line width=0.32pt, black, forget plot]
table {%
4.0625 0
4.3375 0
};
\path [draw=black, fill=chocolate213940, line width=0.32pt]
(axis cs:4.675,11.5866553098103)
--(axis cs:5.225,11.5866553098103)
--(axis cs:5.225,26.8981346751315)
--(axis cs:4.675,26.8981346751315)
--(axis cs:4.675,11.5866553098103)
--cycle;
\addplot [line width=0.32pt, black, forget plot]
table {%
4.95 11.5866553098103
4.95 6.31202990368187
};
\addplot [line width=0.32pt, black, forget plot]
table {%
4.95 26.8981346751315
4.95 49.0861917717333
};
\addplot [line width=0.32pt, black, forget plot]
table {%
4.8125 6.31202990368187
5.0875 6.31202990368187
};
\addplot [line width=0.32pt, black, forget plot]
table {%
4.8125 49.0861917717333
5.0875 49.0861917717333
};
\addplot [black, mark=*, mark size=1.5, mark options={solid,fill=gray}, only marks, forget plot]
table {%
4.95 70.0916313151614
4.95 52.0917617173005
4.95 74.3264893984451
4.95 78.8159718667366
4.95 57.158366504365
4.95 56.0624502860934
};
\path [draw=black, fill=cornflowerblue86180233, line width=0.32pt]
(axis cs:5.425,41.5729943299922)
--(axis cs:5.975,41.5729943299922)
--(axis cs:5.975,73.86509569612)
--(axis cs:5.425,73.86509569612)
--(axis cs:5.425,41.5729943299922)
--cycle;
\addplot [line width=0.32pt, black, forget plot]
table {%
5.7 41.5729943299922
5.7 2.51062069347345
};
\addplot [line width=0.32pt, black, forget plot]
table {%
5.7 73.86509569612
5.7 87.6325472611312
};
\addplot [line width=0.32pt, black, forget plot]
table {%
5.5625 2.51062069347345
5.8375 2.51062069347345
};
\addplot [line width=0.32pt, black, forget plot]
table {%
5.5625 87.6325472611312
5.8375 87.6325472611312
};
\path [draw=black, fill=darkcyan0114178, line width=0.32pt]
(axis cs:6.625,0.493285138964015)
--(axis cs:7.175,0.493285138964015)
--(axis cs:7.175,1.45247260388799)
--(axis cs:6.625,1.45247260388799)
--(axis cs:6.625,0.493285138964015)
--cycle;
\addplot [line width=0.32pt, black, forget plot]
table {%
6.9 0.493285138964015
6.9 0.198680429429343
};
\addplot [line width=0.32pt, black, forget plot]
table {%
6.9 1.45247260388799
6.9 2.81622637768372
};
\addplot [line width=0.32pt, black, forget plot]
table {%
6.7625 0.198680429429343
7.0375 0.198680429429343
};
\addplot [line width=0.32pt, black, forget plot]
table {%
6.7625 2.81622637768372
7.0375 2.81622637768372
};
\addplot [black, mark=*, mark size=1.5, mark options={solid,fill=gray}, only marks, forget plot]
table {%
6.9 3.25830340957718
6.9 3.72790817156372
6.9 3.33005028724365
6.9 4.52781985994532
6.9 4.00908266640188
};
\path [draw=black, fill=orange2301590, line width=0.32pt]
(axis cs:7.375,0.411231408515722)
--(axis cs:7.925,0.411231408515722)
--(axis cs:7.925,1.59991939858287)
--(axis cs:7.375,1.59991939858287)
--(axis cs:7.375,0.411231408515722)
--cycle;
\addplot [line width=0.32pt, black, forget plot]
table {%
7.65 0.411231408515722
7.65 0.0723541270882699
};
\addplot [line width=0.32pt, black, forget plot]
table {%
7.65 1.59991939858287
7.65 3.18973987081109
};
\addplot [line width=0.32pt, black, forget plot]
table {%
7.5125 0.0723541270882699
7.7875 0.0723541270882699
};
\addplot [line width=0.32pt, black, forget plot]
table {%
7.5125 3.18973987081109
7.7875 3.18973987081109
};
\addplot [black, mark=*, mark size=1.5, mark options={solid,fill=gray}, only marks, forget plot]
table {%
7.65 5.47766133421448
7.65 4.23675734475796
7.65 7.3389779168665
7.65 3.54035779117574
7.65 7.79554498033154
7.65 4.3466792895697
7.65 6.56971036257556
7.65 3.83860535565547
7.65 4.13682932482755
7.65 4.69733986132714
7.65 15.1543033042899
7.65 5.46721665218084
7.65 10.6219279888003
};
\path [draw=black, fill=darkcyan0158115, line width=0.32pt]
(axis cs:8.125,0.61554629398236)
--(axis cs:8.675,0.61554629398236)
--(axis cs:8.675,1.45872745296296)
--(axis cs:8.125,1.45872745296296)
--(axis cs:8.125,0.61554629398236)
--cycle;
\addplot [line width=0.32pt, black, forget plot]
table {%
8.4 0.61554629398236
8.4 0.14371138454526
};
\addplot [line width=0.32pt, black, forget plot]
table {%
8.4 1.45872745296296
8.4 2.71083600863597
};
\addplot [line width=0.32pt, black, forget plot]
table {%
8.2625 0.14371138454526
8.5375 0.14371138454526
};
\addplot [line width=0.32pt, black, forget plot]
table {%
8.2625 2.71083600863597
8.5375 2.71083600863597
};
\addplot [black, mark=*, mark size=1.5, mark options={solid,fill=gray}, only marks, forget plot]
table {%
8.4 4.18264730758885
8.4 3.52168300396626
8.4 2.7785116255979
8.4 3.14263665610818
8.4 2.85214117340086
8.4 3.32397652785598
8.4 3.0003897498565
8.4 11.6190123554667
8.4 2.90637220420099
8.4 6.0330502317351
8.4 2.74182253907311
8.4 5.33709129964751
};
\path [draw=black, fill=chocolate213940, line width=0.32pt]
(axis cs:8.875,14.3465480150893)
--(axis cs:9.425,14.3465480150893)
--(axis cs:9.425,33.1602143787026)
--(axis cs:8.875,33.1602143787026)
--(axis cs:8.875,14.3465480150893)
--cycle;
\addplot [line width=0.32pt, black, forget plot]
table {%
9.15 14.3465480150893
9.15 7.4244372412304
};
\addplot [line width=0.32pt, black, forget plot]
table {%
9.15 33.1602143787026
9.15 58.1469738316691
};
\addplot [line width=0.32pt, black, forget plot]
table {%
9.0125 7.4244372412304
9.2875 7.4244372412304
};
\addplot [line width=0.32pt, black, forget plot]
table {%
9.0125 58.1469738316691
9.2875 58.1469738316691
};
\addplot [black, mark=*, mark size=1.5, mark options={solid,fill=gray}, only marks, forget plot]
table {%
9.15 64.2784567804033
9.15 70.2977280958687
9.15 67.5835370618448
9.15 79.6400688425032
9.15 80.3135310473272
};
\path [draw=black, fill=cornflowerblue86180233, line width=0.32pt]
(axis cs:9.625,62.5211172198723)
--(axis cs:10.175,62.5211172198723)
--(axis cs:10.175,83.5958597223979)
--(axis cs:9.625,83.5958597223979)
--(axis cs:9.625,62.5211172198723)
--cycle;
\addplot [line width=0.32pt, black, forget plot]
table {%
9.9 62.5211172198723
9.9 31.2097770463417
};
\addplot [line width=0.32pt, black, forget plot]
table {%
9.9 83.5958597223979
9.9 91.9426220483813
};
\addplot [line width=0.32pt, black, forget plot]
table {%
9.7625 31.2097770463417
10.0375 31.2097770463417
};
\addplot [line width=0.32pt, black, forget plot]
table {%
9.7625 91.9426220483813
10.0375 91.9426220483813
};
\addplot [black, mark=*, mark size=1.5, mark options={solid,fill=gray}, only marks, forget plot]
table {%
9.9 29.1246108724397
9.9 10.323083615261
9.9 25.8401335939914
9.9 12.6493503614123
9.9 8.09573620431554
};
\addplot [line width=0.48pt, black, forget plot]
table {%
-1.775 52.183616382918
-1.225 52.183616382918
};
\addplot [line width=0.48pt, black, forget plot]
table {%
-1.025 0.00608563523720818
-0.475 0.00608563523720818
};
\addplot [line width=0.48pt, black, forget plot]
table {%
-0.275 0
0.275 0
};
\addplot [line width=0.48pt, black, forget plot]
table {%
0.475 26.552242161422
1.025 26.552242161422
};
\addplot [line width=0.48pt, black, forget plot]
table {%
1.225 19.3247569123841
1.775 19.3247569123841
};
\addplot [line width=0.48pt, black, forget plot]
table {%
2.425 13.6727823260571
2.975 13.6727823260571
};
\addplot [line width=0.48pt, black, forget plot]
table {%
3.175 6.14833699081876
3.725 6.14833699081876
};
\addplot [line width=0.48pt, black, forget plot]
table {%
3.925 0
4.475 0
};
\addplot [line width=0.48pt, black, forget plot]
table {%
4.675 16.5752339046328
5.225 16.5752339046328
};
\addplot [line width=0.48pt, black, forget plot]
table {%
5.425 62.0442319200945
5.975 62.0442319200945
};
\addplot [line width=0.48pt, black, forget plot]
table {%
6.625 0.846417175057237
7.175 0.846417175057237
};
\addplot [line width=0.48pt, black, forget plot]
table {%
7.375 0.693916346066836
7.925 0.693916346066836
};
\addplot [line width=0.48pt, black, forget plot]
table {%
8.125 0.984077514059556
8.675 0.984077514059556
};
\addplot [line width=0.48pt, black, forget plot]
table {%
8.875 20.6306694193314
9.425 20.6306694193314
};
\addplot [line width=0.48pt, black, forget plot]
table {%
9.625 76.4928517670088
10.175 76.4928517670088
};

\addlegendimage{only marks, mark=square*, mark size=2.5pt, draw=black, fill={rgb,255:red,0;green,114;blue,178}}
\addlegendentry{Primal}
\addlegendimage{only marks, mark=square*, mark size=2.5pt, draw=black, fill={rgb,255:red,230;green,159;blue,0}}
\addlegendentry{Dual}
\addlegendimage{only marks, mark=square*, mark size=2.5pt, draw=black, fill={rgb,255:red,0;green,158;blue,115}}
\addlegendentry{Update}
\addlegendimage{only marks, mark=square*, mark size=2.5pt, draw=black, fill={rgb,255:red,213;green,94;blue,0}}
\addlegendentry{Misc. SP}
\addlegendimage{only marks, mark=square*, mark size=2.5pt, draw=black, fill={rgb,255:red,86;green,180;blue,233}}
\addlegendentry{MP}
\end{axis}

\end{tikzpicture}

%% file: data/final/time_decomp_boxplot/detailed.tikz
\begin{tikzpicture}
\providecommand{\mathdefault}[1]{#1}

\definecolor{chocolate213940}{RGB}{213,94,0}
\definecolor{cornflowerblue86180233}{RGB}{86,180,233}
\definecolor{darkcyan0114178}{RGB}{0,114,178}
\definecolor{darkcyan0158115}{RGB}{0,158,115}
\definecolor{darkgray176}{RGB}{176,176,176}
\definecolor{gray}{RGB}{128,128,128}
\definecolor{orange2301590}{RGB}{230,159,0}
\definecolor{whitesmoke242}{RGB}{242,242,242}

\begin{axis}[
height=7.14000cm,
legend cell align={left},
legend style={
  font=\footnotesize,
  fill=none,
  fill opacity=1,
  draw opacity=1,
  text opacity=1,
  at={(0.330645161290323,0.995)},
  anchor=north,
  draw=none
},
tick align=outside,
tick pos=left,
width=7.00000cm,
x grid style={darkgray176},
xmin=-2, xmax=10.4,
xtick style={color=black},
xtick={0,4.2,8.4},
xticklabels={\texttt{GRB},\texttt{SLE},\texttt{SLE-U}},
y grid style={darkgray176},
ylabel={Time (h)},
ymin=-0.09, ymax=3.09,
ytick style={color=black},
ytick={-0.5,0,0.5,1,1.5,2,2.5,3,3.5},
yticklabels={
  $\mathdefault{\ensuremath{-}0.5}$,
  $\mathdefault{0.0}$,
  $\mathdefault{0.5}$,
  $\mathdefault{1.0}$,
  $\mathdefault{1.5}$,
  $\mathdefault{2.0}$,
  $\mathdefault{2.5}$,
  $\mathdefault{3.0}$,
  $\mathdefault{3.5}$
}
]

\path [draw=none, fill=whitesmoke242] (axis cs:2.1,-0.08245826330532213) rectangle (axis cs:6.300000000000001,3.082458263305322);
\path [draw=black, fill=darkcyan0114178, line width=0.32pt]
(axis cs:-1.775,0.0499784166666667)
--(axis cs:-1.225,0.0499784166666667)
--(axis cs:-1.225,0.649947805555556)
--(axis cs:-1.775,0.649947805555556)
--(axis cs:-1.775,0.0499784166666667)
--cycle;
\addplot [line width=0.32pt, black, forget plot]
table {%
-1.5 0.0499784166666667
-1.5 0.00556669444444444
};
\addplot [line width=0.32pt, black, forget plot]
table {%
-1.5 0.649947805555556
-1.5 1.54362433333333
};
\addplot [line width=0.32pt, black, forget plot]
table {%
-1.6375 0.00556669444444444
-1.3625 0.00556669444444444
};
\addplot [line width=0.32pt, black, forget plot]
table {%
-1.6375 1.54362433333333
-1.3625 1.54362433333333
};
\addplot [black, mark=*, mark size=1.5, mark options={solid,fill=gray}, only marks, forget plot]
table {%
-1.5 1.68636683333333
-1.5 1.57381525
-1.5 2.52483575
-1.5 1.60758841666667
-1.5 1.573654
};
\path [draw=black, fill=orange2301590, line width=0.32pt]
(axis cs:-1.025,6.91666666666667e-06)
--(axis cs:-0.475,6.91666666666667e-06)
--(axis cs:-0.475,5.32777777777778e-05)
--(axis cs:-1.025,5.32777777777778e-05)
--(axis cs:-1.025,6.91666666666667e-06)
--cycle;
\addplot [line width=0.32pt, black, forget plot]
table {%
-0.75 6.91666666666667e-06
-0.75 1e-06
};
\addplot [line width=0.32pt, black, forget plot]
table {%
-0.75 5.32777777777778e-05
-0.75 0.000119444444444444
};
\addplot [line width=0.32pt, black, forget plot]
table {%
-0.8875 1e-06
-0.6125 1e-06
};
\addplot [line width=0.32pt, black, forget plot]
table {%
-0.8875 0.000119444444444444
-0.6125 0.000119444444444444
};
\addplot [black, mark=*, mark size=1.5, mark options={solid,fill=gray}, only marks, forget plot]
table {%
-0.75 0.000147972222222222
-0.75 0.000124
};
\path [draw=black, fill=darkcyan0158115, line width=0.32pt]
(axis cs:-0.275,0)
--(axis cs:0.275,0)
--(axis cs:0.275,0)
--(axis cs:-0.275,0)
--(axis cs:-0.275,0)
--cycle;
\addplot [line width=0.32pt, black, forget plot]
table {%
0 0
0 0
};
\addplot [line width=0.32pt, black, forget plot]
table {%
0 0
0 0
};
\addplot [line width=0.32pt, black, forget plot]
table {%
-0.1375 0
0.1375 0
};
\addplot [line width=0.32pt, black, forget plot]
table {%
-0.1375 0
0.1375 0
};
\path [draw=black, fill=chocolate213940, line width=0.32pt]
(axis cs:0.475,0.027913)
--(axis cs:1.025,0.027913)
--(axis cs:1.025,0.319777944444444)
--(axis cs:0.475,0.319777944444444)
--(axis cs:0.475,0.027913)
--cycle;
\addplot [line width=0.32pt, black, forget plot]
table {%
0.75 0.027913
0.75 0.00569761111111111
};
\addplot [line width=0.32pt, black, forget plot]
table {%
0.75 0.319777944444444
0.75 0.744173694444444
};
\addplot [line width=0.32pt, black, forget plot]
table {%
0.6125 0.00569761111111111
0.8875 0.00569761111111111
};
\addplot [line width=0.32pt, black, forget plot]
table {%
0.6125 0.744173694444444
0.8875 0.744173694444444
};
\addplot [black, mark=*, mark size=1.5, mark options={solid,fill=gray}, only marks, forget plot]
table {%
0.75 0.848758888888889
};
\path [draw=black, fill=cornflowerblue86180233, line width=0.32pt]
(axis cs:1.225,0.011392526661449)
--(axis cs:1.775,0.011392526661449)
--(axis cs:1.775,0.329585982181019)
--(axis cs:1.225,0.329585982181019)
--(axis cs:1.225,0.011392526661449)
--cycle;
\addplot [line width=0.32pt, black, forget plot]
table {%
1.5 0.011392526661449
1.5 0.000181397535112172
};
\addplot [line width=0.32pt, black, forget plot]
table {%
1.5 0.329585982181019
1.5 0.8017450758596
};
\addplot [line width=0.32pt, black, forget plot]
table {%
1.3625 0.000181397535112172
1.6375 0.000181397535112172
};
\addplot [line width=0.32pt, black, forget plot]
table {%
1.3625 0.8017450758596
1.6375 0.8017450758596
};
\addplot [black, mark=*, mark size=1.5, mark options={solid,fill=gray}, only marks, forget plot]
table {%
1.5 0.811694163045989
1.5 0.908485076955795
1.5 1.23452361972131
1.5 0.896602891161495
1.5 1.25706184500196
1.5 1.02219038642714
1.5 1.12954214609729
1.5 0.872489257533815
1.5 0.83043807421165
1.5 0.913895681142807
1.5 1.29694236833562
};
\path [draw=black, fill=darkcyan0114178, line width=0.32pt]
(axis cs:2.425,0.00657152777777778)
--(axis cs:2.975,0.00657152777777778)
--(axis cs:2.975,0.0553189583333333)
--(axis cs:2.425,0.0553189583333333)
--(axis cs:2.425,0.00657152777777778)
--cycle;
\addplot [line width=0.32pt, black, forget plot]
table {%
2.7 0.00657152777777778
2.7 0.000200722222222222
};
\addplot [line width=0.32pt, black, forget plot]
table {%
2.7 0.0553189583333333
2.7 0.122670486111111
};
\addplot [line width=0.32pt, black, forget plot]
table {%
2.5625 0.000200722222222222
2.8375 0.000200722222222222
};
\addplot [line width=0.32pt, black, forget plot]
table {%
2.5625 0.122670486111111
2.8375 0.122670486111111
};
\addplot [black, mark=*, mark size=1.5, mark options={solid,fill=gray}, only marks, forget plot]
table {%
2.7 0.287904986111111
2.7 0.304822388888889
2.7 0.221589069444444
2.7 0.33857425
2.7 0.129259708333333
};
\path [draw=black, fill=orange2301590, line width=0.32pt]
(axis cs:3.175,0.00275033333333333)
--(axis cs:3.725,0.00275033333333333)
--(axis cs:3.725,0.0238449166666667)
--(axis cs:3.175,0.0238449166666667)
--(axis cs:3.175,0.00275033333333333)
--cycle;
\addplot [line width=0.32pt, black, forget plot]
table {%
3.45 0.00275033333333333
3.45 7.63611111111111e-05
};
\addplot [line width=0.32pt, black, forget plot]
table {%
3.45 0.0238449166666667
3.45 0.0547240277777778
};
\addplot [line width=0.32pt, black, forget plot]
table {%
3.3125 7.63611111111111e-05
3.5875 7.63611111111111e-05
};
\addplot [line width=0.32pt, black, forget plot]
table {%
3.3125 0.0547240277777778
3.5875 0.0547240277777778
};
\addplot [black, mark=*, mark size=1.5, mark options={solid,fill=gray}, only marks, forget plot]
table {%
3.45 0.135296958333333
3.45 0.0560892916666667
3.45 0.148277277777778
3.45 0.120972375
3.45 0.0707921805555556
3.45 0.201066527777778
3.45 0.0794882916666667
3.45 0.0708873333333333
3.45 0.0641700694444444
};
\path [draw=black, fill=darkcyan0158115, line width=0.32pt]
(axis cs:3.925,0)
--(axis cs:4.475,0)
--(axis cs:4.475,0)
--(axis cs:3.925,0)
--(axis cs:3.925,0)
--cycle;
\addplot [line width=0.32pt, black, forget plot]
table {%
4.2 0
4.2 0
};
\addplot [line width=0.32pt, black, forget plot]
table {%
4.2 0
4.2 0
};
\addplot [line width=0.32pt, black, forget plot]
table {%
4.0625 0
4.3375 0
};
\addplot [line width=0.32pt, black, forget plot]
table {%
4.0625 0
4.3375 0
};
\path [draw=black, fill=chocolate213940, line width=0.32pt]
(axis cs:4.675,0.00879358333333333)
--(axis cs:5.225,0.00879358333333333)
--(axis cs:5.225,0.0648679444444444)
--(axis cs:4.675,0.0648679444444444)
--(axis cs:4.675,0.00879358333333333)
--cycle;
\addplot [line width=0.32pt, black, forget plot]
table {%
4.95 0.00879358333333333
4.95 0.00243583333333333
};
\addplot [line width=0.32pt, black, forget plot]
table {%
4.95 0.0648679444444444
4.95 0.136364444444444
};
\addplot [line width=0.32pt, black, forget plot]
table {%
4.8125 0.00243583333333333
5.0875 0.00243583333333333
};
\addplot [line width=0.32pt, black, forget plot]
table {%
4.8125 0.136364444444444
5.0875 0.136364444444444
};
\addplot [black, mark=*, mark size=1.5, mark options={solid,fill=gray}, only marks, forget plot]
table {%
4.95 0.150901638888889
4.95 0.215750333333333
4.95 0.233893305555556
4.95 0.199968027777778
4.95 0.276820166666667
};
\path [draw=black, fill=cornflowerblue86180233, line width=0.32pt]
(axis cs:5.425,0.0143159653023614)
--(axis cs:5.975,0.0143159653023614)
--(axis cs:5.975,0.385994056698057)
--(axis cs:5.425,0.385994056698057)
--(axis cs:5.425,0.0143159653023614)
--cycle;
\addplot [line width=0.32pt, black, forget plot]
table {%
5.7 0.0143159653023614
5.7 0.000194444431304932
};
\addplot [line width=0.32pt, black, forget plot]
table {%
5.7 0.385994056698057
5.7 0.93836492776627
};
\addplot [line width=0.32pt, black, forget plot]
table {%
5.5625 0.000194444431304932
5.8375 0.000194444431304932
};
\addplot [line width=0.32pt, black, forget plot]
table {%
5.5625 0.93836492776627
5.8375 0.93836492776627
};
\addplot [black, mark=*, mark size=1.5, mark options={solid,fill=gray}, only marks, forget plot]
table {%
5.7 1.38929788275592
5.7 1.39941729415226
5.7 2.41797178671275
5.7 0.997729107246081
5.7 2.34203862389533
5.7 2.54992467054696
5.7 1.00232757696703
5.7 1.1348817352829
5.7 2.29596308744897
5.7 0.971534359170914
5.7 1.01667029530938
5.7 1.28754457333332
};
\path [draw=black, fill=darkcyan0114178, line width=0.32pt]
(axis cs:6.625,0.000290361111111111)
--(axis cs:7.175,0.000290361111111111)
--(axis cs:7.175,0.00228005555555556)
--(axis cs:6.625,0.00228005555555556)
--(axis cs:6.625,0.000290361111111111)
--cycle;
\addplot [line width=0.32pt, black, forget plot]
table {%
6.9 0.000290361111111111
6.9 3.20833333333333e-05
};
\addplot [line width=0.32pt, black, forget plot]
table {%
6.9 0.00228005555555556
6.9 0.00455530555555556
};
\addplot [line width=0.32pt, black, forget plot]
table {%
6.7625 3.20833333333333e-05
7.0375 3.20833333333333e-05
};
\addplot [line width=0.32pt, black, forget plot]
table {%
6.7625 0.00455530555555556
7.0375 0.00455530555555556
};
\addplot [black, mark=*, mark size=1.5, mark options={solid,fill=gray}, only marks, forget plot]
table {%
6.9 0.0113614444444444
6.9 0.0129329166666667
6.9 0.0108510833333333
6.9 0.0294470694444444
6.9 0.0232310138888889
6.9 0.0527193611111111
6.9 0.0211826111111111
6.9 0.0262084444444444
6.9 0.00574386111111111
6.9 0.00814530555555555
6.9 0.0132544027777778
6.9 0.0140620277777778
6.9 0.00760243055555556
6.9 0.0103960555555556
};
\path [draw=black, fill=orange2301590, line width=0.32pt]
(axis cs:7.375,0.000263861111111111)
--(axis cs:7.925,0.000263861111111111)
--(axis cs:7.925,0.002904125)
--(axis cs:7.375,0.002904125)
--(axis cs:7.375,0.000263861111111111)
--cycle;
\addplot [line width=0.32pt, black, forget plot]
table {%
7.65 0.000263861111111111
7.65 1.80555555555556e-05
};
\addplot [line width=0.32pt, black, forget plot]
table {%
7.65 0.002904125
7.65 0.00664684722222222
};
\addplot [line width=0.32pt, black, forget plot]
table {%
7.5125 1.80555555555556e-05
7.7875 1.80555555555556e-05
};
\addplot [line width=0.32pt, black, forget plot]
table {%
7.5125 0.00664684722222222
7.7875 0.00664684722222222
};
\addplot [black, mark=*, mark size=1.5, mark options={solid,fill=gray}, only marks, forget plot]
table {%
7.65 0.00742891666666667
7.65 0.0162446944444444
7.65 0.00820291666666667
7.65 0.0226783194444444
7.65 0.0110729722222222
7.65 0.0168625416666667
7.65 0.0363288055555556
7.65 0.0118983611111111
7.65 0.00742451388888889
7.65 0.00764827777777778
7.65 0.0142931388888889
7.65 0.0173290555555556
7.65 0.0189276388888889
7.65 0.00900116666666667
};
\path [draw=black, fill=darkcyan0158115, line width=0.32pt]
(axis cs:8.125,0.000278305555555556)
--(axis cs:8.675,0.000278305555555556)
--(axis cs:8.675,0.00327697222222222)
--(axis cs:8.125,0.00327697222222222)
--(axis cs:8.125,0.000278305555555556)
--cycle;
\addplot [line width=0.32pt, black, forget plot]
table {%
8.4 0.000278305555555556
8.4 1.8e-05
};
\addplot [line width=0.32pt, black, forget plot]
table {%
8.4 0.00327697222222222
8.4 0.00776608333333333
};
\addplot [line width=0.32pt, black, forget plot]
table {%
8.2625 1.8e-05
8.5375 1.8e-05
};
\addplot [line width=0.32pt, black, forget plot]
table {%
8.2625 0.00776608333333333
8.5375 0.00776608333333333
};
\addplot [black, mark=*, mark size=1.5, mark options={solid,fill=gray}, only marks, forget plot]
table {%
8.4 0.00873022222222222
8.4 0.0186932222222222
8.4 0.0222955
8.4 0.0110446666666667
8.4 0.00821852777777778
8.4 0.0148298888888889
8.4 0.0279690277777778
8.4 0.00999041666666667
8.4 0.00917002777777778
8.4 0.0129777222222222
8.4 0.0195088888888889
8.4 0.0101206388888889
8.4 0.0161365
8.4 0.0251425277777778
8.4 0.0101374722222222
};
\path [draw=black, fill=chocolate213940, line width=0.32pt]
(axis cs:8.875,0.00741333333333333)
--(axis cs:9.425,0.00741333333333333)
--(axis cs:9.425,0.0643855833333333)
--(axis cs:8.875,0.0643855833333333)
--(axis cs:8.875,0.00741333333333333)
--cycle;
\addplot [line width=0.32pt, black, forget plot]
table {%
9.15 0.00741333333333333
9.15 0.000810888888888889
};
\addplot [line width=0.32pt, black, forget plot]
table {%
9.15 0.0643855833333333
9.15 0.132180777777778
};
\addplot [line width=0.32pt, black, forget plot]
table {%
9.0125 0.000810888888888889
9.2875 0.000810888888888889
};
\addplot [line width=0.32pt, black, forget plot]
table {%
9.0125 0.132180777777778
9.2875 0.132180777777778
};
\addplot [black, mark=*, mark size=1.5, mark options={solid,fill=gray}, only marks, forget plot]
table {%
9.15 0.175814277777778
9.15 0.162009222222222
9.15 0.168085944444444
9.15 0.702642527777778
9.15 0.552504416666667
9.15 0.170776388888889
9.15 0.550465194444444
9.15 0.389030416666667
9.15 0.261029888888889
9.15 0.165207944444444
9.15 0.172122111111111
9.15 0.203406833333333
9.15 0.289649444444444
9.15 0.162531638888889
};
\path [draw=black, fill=cornflowerblue86180233, line width=0.32pt]
(axis cs:9.625,0.0189949995010164)
--(axis cs:10.175,0.0189949995010164)
--(axis cs:10.175,0.41900037250381)
--(axis cs:9.625,0.41900037250381)
--(axis cs:9.625,0.0189949995010164)
--cycle;
\addplot [line width=0.32pt, black, forget plot]
table {%
9.9 0.0189949995010164
9.9 0.000142124428431193
};
\addplot [line width=0.32pt, black, forget plot]
table {%
9.9 0.41900037250381
9.9 1.01677400192653
};
\addplot [line width=0.32pt, black, forget plot]
table {%
9.7625 0.000142124428431193
10.0375 0.000142124428431193
};
\addplot [line width=0.32pt, black, forget plot]
table {%
9.7625 1.01677400192653
10.0375 1.01677400192653
};
\addplot [black, mark=*, mark size=1.5, mark options={solid,fill=gray}, only marks, forget plot]
table {%
9.9 1.27489075085502
9.9 1.34455875865131
9.9 1.47478890162638
9.9 2.37418126309776
9.9 2.40934579780261
9.9 2.37778823004479
9.9 1.62879537694931
9.9 1.4039973686391
9.9 1.0596122416723
9.9 1.28727600891738
9.9 1.03877192752986
9.9 1.37598904717424
9.9 1.57991653973707
};
\addplot [line width=0.48pt, black, forget plot]
table {%
-1.775 0.179478805555556
-1.225 0.179478805555556
};
\addplot [line width=0.48pt, black, forget plot]
table {%
-1.025 2.43888888888889e-05
-0.475 2.43888888888889e-05
};
\addplot [line width=0.48pt, black, forget plot]
table {%
-0.275 0
0.275 0
};
\addplot [line width=0.48pt, black, forget plot]
table {%
0.475 0.10279775
1.025 0.10279775
};
\addplot [line width=0.48pt, black, forget plot]
table {%
1.225 0.0768332277951772
1.775 0.0768332277951772
};
\addplot [line width=0.48pt, black, forget plot]
table {%
2.425 0.0191396805555556
2.975 0.0191396805555556
};
\addplot [line width=0.48pt, black, forget plot]
table {%
3.175 0.0100780694444444
3.725 0.0100780694444444
};
\addplot [line width=0.48pt, black, forget plot]
table {%
3.925 0
4.475 0
};
\addplot [line width=0.48pt, black, forget plot]
table {%
4.675 0.0256523333333333
5.225 0.0256523333333333
};
\addplot [line width=0.48pt, black, forget plot]
table {%
5.425 0.0930915408346388
5.975 0.0930915408346388
};
\addplot [line width=0.48pt, black, forget plot]
table {%
6.625 0.000857194444444444
7.175 0.000857194444444444
};
\addplot [line width=0.48pt, black, forget plot]
table {%
7.375 0.000846
7.925 0.000846
};
\addplot [line width=0.48pt, black, forget plot]
table {%
8.125 0.00103186111111111
8.675 0.00103186111111111
};
\addplot [line width=0.48pt, black, forget plot]
table {%
8.875 0.0251693333333333
9.425 0.0251693333333333
};
\addplot [line width=0.48pt, black, forget plot]
table {%
9.625 0.0904667630755107
10.175 0.0904667630755107
};

\addlegendimage{only marks, mark=square*, mark size=2.5pt, draw=black, fill={rgb,255:red,0;green,114;blue,178}}
\addlegendentry{Primal}
\addlegendimage{only marks, mark=square*, mark size=2.5pt, draw=black, fill={rgb,255:red,230;green,159;blue,0}}
\addlegendentry{Dual}
\addlegendimage{only marks, mark=square*, mark size=2.5pt, draw=black, fill={rgb,255:red,0;green,158;blue,115}}
\addlegendentry{Update}
\addlegendimage{only marks, mark=square*, mark size=2.5pt, draw=black, fill={rgb,255:red,213;green,94;blue,0}}
\addlegendentry{Misc. SP}
\addlegendimage{only marks, mark=square*, mark size=2.5pt, draw=black, fill={rgb,255:red,86;green,180;blue,233}}
\addlegendentry{MP}
\end{axis}

\end{tikzpicture}

%% file: data/final/performance_profiles/solve_time_scaled/scatter_sorted_by_sle_runtime.tikz
\begin{tikzpicture}
\providecommand{\mathdefault}[1]{#1}

\definecolor{darkcyan0114178}{RGB}{0,114,178}
\definecolor{darkgray176}{RGB}{176,176,176}
\definecolor{dimgray102}{RGB}{102,102,102}
\definecolor{silver200}{RGB}{200,200,200}

\begin{axis}[
height=5.5cm,
legend cell align={left},
legend style={
  font=\footnotesize,
  fill=none,
  fill opacity=1,
  draw opacity=1,
  text opacity=1,
  at={(0.97,0.03)},
  anchor=south east,
  draw=none
},
log basis y={10},
tick align=outside,
tick pos=left,
width=7.3cm,
xlabel={Instance (sorted by running time)},
xmajorgrids,
xmin=-4, xmax=152,
y grid style={gray,opacity=0.3},
ylabel={Runtime (s)},
ymajorgrids,
ymin=0.150938271604938, ymax=1514.54205,
ymode=log,
ytick style={color=black},
ytick={0.01,0.1,1,10,100,1000,10000,100000},
yticklabels={
  $\mathdefault{10^{-2}}$,
  $\mathdefault{10^{-1}}$,
  $\mathdefault{10^{0}}$,
  $\mathdefault{10^{1}}$,
  $\mathdefault{10^{2}}$,
  $\mathdefault{10^{3}}$,
  $\mathdefault{10^{4}}$,
  $\mathdefault{10^{5}}$
}
]
\addplot [draw=darkcyan0114178, fill=darkcyan0114178, mark=*, only marks, opacity=0.9]
table{%
x  y
-0.1 0.6199
0.9 0.9605
1.9 1.0339
2.9 1.0584
3.9 0.8229
4.9 2.5907
5.9 1.3059
6.9 2.2973
7.9 0.6113
8.9 1.4591
9.9 1.3709
10.9 0.6246
11.9 0.6712
12.9 0.9167
13.9 1.9676
14.9 2.1302
15.9 1.0316
16.9 1.5017
17.9 1.7956
18.9 1.8958
19.9 1.5288
20.9 1.7904
21.9 1.1692
22.9 1.1279
23.9 2.3326
24.9 2.1797
25.9 1.9781
26.9 0.8789
27.9 1.8491
28.9 1.403
29.9 1.5676
30.9 5.293
31.9 1.1516
32.9 4.0366
33.9 3.1
34.9 2.9236
35.9 2.3182
36.9 2.2457
37.9 6.0634
38.9 1.9039
39.9 4.3274
40.9 2.9249
41.9 5.4051
42.9 6.0017
43.9 3.7332
44.9 7.1074
45.9 3.7917
46.9 3.1314
47.9 14.2923
48.9 5.8525
49.9 11.0187
50.9 4.0235
51.9 4.2704
52.9 10.118
53.9 20.6607
54.9 7.8599
55.9 5.8227
56.9 5.5845
57.9 5.1549
58.9 9.8558
59.9 6.5986
60.9 7.7735
61.9 16.086
62.9 4.9076
63.9 9.3229
64.9 13.9138
65.9 15.2052
66.9 6.7927
67.9 6.3729
68.9 8.8712
69.9 6.289
70.9 25.7879
71.9 9.1236
72.9 9.5067
73.9 13.7675
74.9 88.4353
75.9 10.1456
76.9 10.8147
77.9 21.199
78.9 20.9164
79.9 12.3868
80.9 13.0676
81.9 8.4872
82.9 8.5383
83.9 10.6178
84.9 7.8834
85.9 6.4156
86.9 8.4204
87.9 12.0166
88.9 21.4235
89.9 26.3918
90.9 12.1456
91.9 12.5761
92.9 12.8933
93.9 26.388
94.9 10.4499
95.9 20.6738
96.9 65.1166
97.9 23.986
98.9 51.7178
99.9 224.6006
100.9 70.7954
101.9 254.947
102.9 22.2202
103.9 16.765
104.9 20.2613
105.9 27.8824
106.9 22.001
107.9 18.1536
108.9 26.7909
109.9 42.6008
110.9 373.961
111.9 49.581
112.9 32.0588
113.9 37.5362
114.9 26.9003
115.9 19.7703
116.9 23.102
117.9 20.5707
118.9 29.0552
119.9 41.4333
120.9 100.0232
121.9 28.5578
122.9 28.7596
123.9 33.0201
124.9 53.0299
125.9 98.0381
126.9 44.9463
127.9 36.4493
128.9 77.9979
129.9 118.6827
130.9 39.1546
131.9 34.0705
132.9 95.4668
133.9 42.2473
134.9 52.7811
135.9 42.5008
136.9 30.6628
137.9 52.0292
138.9 98.0708
139.9 77.6802
140.9 157.8497
141.9 129.7473
142.9 104.6566
143.9 116.323
144.9 31.2991
145.9 50.4329
146.9 219.78
147.9 117.9048
};
\addlegendentry{\texttt{SLE-U}}
\addplot [draw=dimgray102, fill=dimgray102, mark=x, only marks, opacity=0.9]
table{%
x  y
0 1.23412259358691
1 1.71930918715772
2 2.30469120891086
3 2.88792119475639
4 2.89681950800901
5 3.53923097369241
6 3.63354061221492
7 3.89382562255148
8 4.00799600532402
9 4.01338880249036
10 4.10155807449505
11 4.50136962573042
12 4.70800502659842
13 4.7746882610734
14 4.96785626947122
15 5.1080790448149
16 5.66149430520219
17 5.76694235301267
18 5.89754324996975
19 6.26090191397872
20 6.4022498629594
21 6.61322604070225
22 6.78858518720533
23 6.79221805301187
24 6.93224080007325
25 6.95774040799196
26 6.97410836191688
27 6.97741498140046
28 7.17104426653987
29 7.34530265679217
30 8.05162587038193
31 8.92814527800085
32 8.94660962899012
33 8.98399288044777
34 8.99511044694738
35 9.45574491831898
36 9.89793982233814
37 9.96159907444901
38 10.3804983348141
39 10.5023887543242
40 11.1263784603352
41 11.1541420681329
42 11.189724814668
43 11.268328866317
44 12.132900663001
45 13.7447920592257
46 13.7923674324067
47 13.8468488897474
48 14.1923031462104
49 14.3493840050746
50 14.7789831960444
51 14.7791045862401
52 14.8948231058786
53 15.3179454556984
54 15.9300401361168
55 15.958906937031
56 16.2104759871206
57 16.4751562256311
58 16.5968840745784
59 17.1397539449039
60 17.9131043841364
61 17.9322846394206
62 17.9331190153117
63 18.1327493870981
64 18.2312945001755
65 18.6821065270425
66 18.8667846383069
67 19.33050391391
68 19.4819898510824
69 20.3241691640568
70 20.367562868101
71 21.4847714842152
72 22.0090889852769
73 22.1264247050542
74 22.3239208228576
75 23.2897549339553
76 23.4446864708846
77 23.526537627428
78 23.852749046359
79 24.0360201625022
80 24.1929309175155
81 24.3636127647568
82 24.5947264015834
83 24.8156463852554
84 25.257657263306
85 25.6731948894883
86 25.8695225710535
87 26.5154218578346
88 26.8694744551469
89 28.2485187083479
90 29.4959846550109
91 29.6746739100618
92 29.7604381163078
93 30.0550187993736
94 30.3414581786093
95 31.0944011267994
96 31.5811462981718
97 32.1460471442343
98 33.3289811561316
99 33.4154063224106
100 34.7748951027409
101 34.7995472378696
102 35.9888102415473
103 37.7557415360935
104 38.7671222402727
105 38.8089440939559
106 39.6207528561212
107 40.5624714017286
108 43.8284655653287
109 45.1266336986239
110 45.5592332747785
111 45.6806974507389
112 45.7830062613847
113 46.1757129445088
114 46.673339259097
115 46.8449667462555
116 47.1743254420772
117 47.6142800970439
118 48.2159120229148
119 49.1173202794539
120 49.5794995203265
121 51.2450396516238
122 54.2174376285069
123 54.3210380487324
124 54.4980016833149
125 54.6480959636825
126 54.8517453667667
127 54.9258564592708
128 55.7850378285099
129 56.1491903750309
130 57.6505270597602
131 58.185253322247
132 58.3971060093447
133 59.8555280456847
134 60.6720293098753
135 61.3529544964667
136 61.4730909282338
137 63.3479865384306
138 67.5195227088488
139 71.2653907660933
140 72.0128484255786
141 72.090955578722
142 72.1121145844076
143 77.0982135069731
144 78.5628486610884
145 114.998979905997
146 132.403694485026
147 138.568387220516
148 155.693867519386
};
\addlegendentry{\texttt{SLE}$^{2/3}$}
\addplot [draw=silver200, fill=silver200, mark=x, only marks, opacity=0.9]
table{%
x  y
0.1 1.17089709197692
1.1 1.50146595032988
2.1 1.87050795240223
3.1 2.21533293208944
4.1 2.22045040476026
5.1 2.58036819078208
6.1 2.63176746693168
7.1 2.77193073506536
8.1 2.83266658821683
9.1 2.83552464281304
10.1 2.88211727728071
11.1 3.09035596655143
12.1 3.19615393872072
13.1 3.23004643929464
14.1 3.32756667852051
15.1 3.39776397061362
16.1 3.67027246944965
17.1 3.72142445845673
18.1 3.78445504663485
19.1 3.95801718035685
20.1 4.02484782321021
21.1 4.12391803992271
22.1 4.20566284906434
23.1 4.20735071036395
24.1 4.2722359485403
25.1 4.28401680668972
26.1 4.29157313813944
27.1 4.29309911369398
28.1 4.38214559319975
29.1 4.4617709488498
30.1 4.77983263305317
31.1 5.16500726040148
32.1 5.17301652810041
33.1 5.18921959450552
34.1 5.19403504031307
35.1 5.39227224831981
36.1 5.58031361125878
37.1 5.60720964473418
38.1 5.78313928588963
39.1 5.83399520054653
40.1 6.09207682157735
41.1 6.10347442036091
42.1 6.11807159160466
43.1 6.15027641655235
44.1 6.50089993770093
45.1 7.13844520886726
46.1 7.15696863203969
47.1 7.17816132446186
48.1 7.31205853368256
49.1 7.37267251408877
50.1 7.53760572065162
51.1 7.53765215435151
52.1 7.58187311948703
53.1 7.74284185554632
54.1 7.97375068584415
55.1 7.98458514889784
56.1 8.07879941575479
57.1 8.17753018948876
58.1 8.22280365812051
59.1 8.42371058382231
60.1 8.70719242925066
61.1 8.71418384015394
62.1 8.7144879367637
63.1 8.78714401839414
64.1 8.82293601926252
65.1 8.98606142867942
66.1 9.05260183593645
67.1 9.21896957365627
68.1 9.27310088373894
69.1 9.57215231805261
70.1 9.58747620596787
71.1 9.97925347909352
72.1 10.1613532563335
73.1 10.2019556948656
74.1 10.2701752662747
75.1 10.6016555310951
76.1 10.6545060889747
77.1 10.6823920542171
78.1 10.7932895819579
79.1 10.8554272140713
80.1 10.9085333569642
81.1 10.9662026244275
82.1 11.0441296624044
83.1 11.1184486327905
84.1 11.2666499013682
85.1 11.4053846932052
86.1 11.4707366807891
87.1 11.6848705598308
88.1 11.8016947935455
89.1 12.2531220511346
90.1 12.6567491876864
91.1 12.7142125198535
92.1 12.7417620445526
93.1 12.8362377665732
94.1 12.9278807234597
95.1 13.1677522759201
96.1 13.3220456387148
97.1 13.50037036529
98.1 13.8712796814137
99.1 13.8982480910365
100.1 14.3202094956743
101.1 14.3278225840495
102.1 14.6935121737453
103.1 15.2313098583149
104.1 15.5363026489574
105.1 15.5488713416762
106.1 15.792178443774
107.1 16.0728653326033
108.1 17.0340218386616
109.1 17.4110396013564
110.1 17.5360713958401
111.1 17.5711240391729
112.1 17.6006306705186
113.1 17.7137376067277
114.1 17.8567186235322
115.1 17.9059431474581
116.1 18.0002805533692
117.1 18.1260392805489
118.1 18.2975435509797
119.1 18.5535064071458
120.1 18.6842901925655
121.1 19.1530911343313
122.1 19.9804028988407
123.1 20.0090304612692
124.1 20.0578986935322
125.1 20.0993159087567
126.1 20.1554657599372
127.1 20.1758865976194
128.1 20.4121287473894
129.1 20.5119818642666
130.1 20.9219645349092
131.1 21.0673396516978
132.1 21.1248431946843
133.1 21.5193029626891
134.1 21.7390915173565
135.1 21.9218201799029
136.1 21.9540064680686
137.1 22.4543047097878
138.1 23.5543987399382
139.1 24.5278229771825
140.1 24.7205137487068
141.1 24.7406204449282
142.1 24.7460663540693
143.1 26.0185587610075
144.1 26.3883913871232
145.1 35.1172493228043
146.1 39.0323839907326
147.1 40.3876070100718
148.1 44.0761477445568
};
\addlegendentry{\texttt{SLE}$^{1/2}$}
\end{axis}

\end{tikzpicture}%

%% file: chapters/6-conclusion.tex
Potential-based flow problems cover a range of different network flow related problems.
Due to their application to several infrastructure domains, they are of significant practical importance.
Affine potential functions cover, for example, applications in lossless DC power systems.
Benders decomposition is a common exact algorithm in this area, particularly when considering different scenarios simultaneously.
Even though the resulting subproblems are linear programs, typically solved with a simplex algorithm,
they can account for a major share of the total running time.

Our contribution builds on the classical fact that in a passive, potential-based network with fixed topology and nodal demands, the resulting flow is unique \citep[as observed by][]{birkhoff1956potential}. 
We conclude that in the subproblems ``there is nothing to optimize,'' and thus the general simplex method is a considerable overhead.
Our work can be interpreted as exposing this special problem structure to the simplex algorithm.
Primal solutions are obtained from solving a single, carefully crafted system of linear equations (SLE).
More importantly, dual solutions are required to construct Benders cuts.
It is only natural that dual solutions arise from solving the same SLE (``transposed''),
and we develop that this is indeed the case (and deriving a dual direction of unboundedness requires only a little more work).
Overall, we need to solve at most two SLEs per subproblem.
Just like in the simplex method, we accelerate solving the SLEs by factorizing the coefficient matrices.
Over multiple subproblems, factorizations are updated, thereby reducing the theoretical complexity and practical running times in certain, natural cases. 
It should be noted that our SLEs have smaller coefficient matrices than those that would be solved in the simplex method.
Overall, this gives us an efficient ``simplex-free,'' but ``simplex-inspired'' subproblem solving algorithm for affine potential-based flow problems with switchable arcs and arc-based contingency scenarios.

Using our approach, we were able to obtain significant improvements over Gurobi as a baseline on instances synthesized from the PGLib-OPF \citep{pglib}.
The shifted geometric mean of the overall running time is reduced by \qty[round-mode=figures]{71}{\percent}.
The subproblem solving time share decreases from \qty[round-mode=figures]{49}{\percent} for Gurobi to a little more than \qty[round-mode=figures, round-precision=1]{2}{\percent} for our approach.
Empirical results also match theoretical results regarding factorization update complexity.

Finally, we would like to mention two potential further developments of our work, one easier, one more challenging:
\begin{itemize}
    \item The Benders algorithm presented in this work, where we solve all subproblems in every iteration and obtain at most one cut per subproblem, is only one option for handling multiple scenarios.
    It might be worthwhile to compute additional cuts, for example by choosing multiple violated components for a feasibility cut or by solving additional, related subproblems in a given Benders iteration.
    Adding additional cuts can lead to a reduced number of iterations overall, and thus less time spent in the master problem.
    Our algorithm allows for exploring many additional subproblems at very low computational cost.
    \item A natural extension would be to consider flow problems with non-affine potential functions.
    However, our algorithm relies on the subproblems being linear to construct linear systems.
    Nevertheless, it suffices that potential functions $\PotFSym$ are continuous, strictly increasing and $\PotFSym(0) = 0$ for the flow to be unique \citep{pfetsch_potential-based_2026}, so there is room for improvement also in these cases.
\end{itemize}